\documentclass[11pt]{amsart}

\usepackage{amsfonts}
\usepackage{amscd}
\usepackage{amssymb}
\usepackage{amsthm}
\usepackage{amsmath} 
\usepackage{bbm} 
\usepackage{epsfig}
\usepackage{enumerate}

\input xy
\xyoption{all}

\usepackage{color}
\definecolor{dred}{rgb}{.65, 0, 0.15}
\definecolor{dblue}{rgb}{0.15, 0, 0.65}
\definecolor{grey}{rgb}{.7, .7, .7}
\definecolor{dgrey}{rgb}{.1, .1, .1}
\definecolor{orange}{rgb}{1, .3, 0}
\definecolor{brick}{rgb}{.5, 0, 0}

\theoremstyle{plain}
\newtheorem{thm}{Theorem}[section]
\newtheorem{thm*}{Theorem}
\newtheorem{lemma}[thm]{Lemma}
\newtheorem{prop}[thm]{Proposition}
\newtheorem{cor}[thm]{Corollary}
\theoremstyle{definition}

\newtheorem{remark}[thm]{Remark}
\theoremstyle{example}
\newtheorem{example}{Example}[section]
\theoremstyle{remark}

\numberwithin{equation}{section}

\def\cA{\mathcal{A}}
\def\cB{\mathcal{B}}
\def\cC{\mathcal{C}}

\def\cG{\mathcal{G}}
\def\cH{\mathcal{H}}

\def\cL{\mathcal{L}}

\def\cO{\mathcal{O}}
\def\cP{\mathcal{P}}

\def\cR{\mathcal{R}}

\def\cT{\mathcal{T}}
\def\cU{\mathcal{U}}

\def\CC{\mathbb{C}}

\def\ZZ{\mathbb{Z}}

\def\ww{\mathrm{w}}

\def\fg{\mathfrak{g}} 
\def\fh{\mathfrak{h}}

\def\fn{\mathfrak{n}}

\def\fgl{\mathfrak{gl}}  
\def\fsl{\mathfrak{sl}}  
\def\fso{\mathfrak{so}}  
\def\fsp{\mathfrak{sp}}  

\def\ad{\mathrm{ad}}

\def\dim{\mathrm{dim}} 
\def\ext{\mathrm{ext}} 
\def\End{\mathrm{End}} 
 
\def\Hom{\mathrm{Hom}} 
\def\half{\hbox{$\frac12$}}
\def\id{\mathrm{id}}  
\def\Ind{\mathrm{Ind}}

\def\Res{\mathrm{Res}}

\def\tr{\mathrm{tr}} 
 
\def\vep{\varepsilon}

\def\tensorspace{$M\otimes N \otimes V^{\otimes k}$}

\def\<{\langle}
\def\>{\rangle}

\def\VVsym{L\left(~\!\begin{picture}(20, 10)(0,12)
	\color{black}
	\thicklines
		\multiput(0,10)(10,0){3} {\line(0,1){10}}
		\multiput(0,10)(0,10){2} {\line(1,0){20}}
	\end{picture}~\!\right) }

\def\VVantisym{  L\left(~\!\begin{picture}(10, 15)(0,17)
	\color{black}
	\thicklines
		\multiput(0,10)(10,0){2} {\line(0,1){20}}
		\multiput(0,10)(0,10){3} {\line(1,0){10}}
	\end{picture}~\!\right)  }

\newcommand{\comment}[1]{}

\usepackage[pdftex,bookmarks]{hyperref}

\title{Degenerate two-boundary centralizer algebras}
\author{Zajj Daugherty}
\address{Department of Mathematics\\
Dartmouth College\\
6188 Kemeny Hall\\
Hanover, NH 03755
}
\email{zajj.b.daugherty@dartmouth.edu}
\urladdr{http://www.math.dartmouth.edu/~zdaugherty}

\date{\today}

\thanks{This paper is a revised and streamlined version of thesis work \cite{Da}, and is the result of research and teaching initiatives supported by the National Science Foundation under Grant No. 0353038. The author would like to thank, among others, Georgia Benkart, Matthew Davis, and especially Arun Ram for many helpful conversations. The author would also like to thank the referee for very thorough and thoughtful guidance.}


\begin{document}

\maketitle

\begin{abstract} 
Diagram algebras (e.g. graded braid groups, Hecke algebras, Brauer algebras) arise as tensor power centralizer algebras, algebras of commuting operators for a Lie algebra action on a tensor space. This work explores centralizers of the action of a complex reductive Lie algebra $\mathfrak{g}$ on tensor space of the form $M \otimes N \otimes V^{\otimes k}$. We define the degenerate two-boundary braid algebra $\mathcal{G}_k$ and show that centralizer algebras contain quotients of this algebra in a general setting. As an example, we study in detail the combinatorics of special cases corresponding to Lie algebras $\mathfrak{gl}_n$ and $\mathfrak{sl}_n$ and modules $M$ and $N$ indexed by rectangular partitions. For this setting, we define the degenerate extended two-boundary Hecke algebra $\mathcal{H}_k^{\mathrm{ext}}$ as a quotient of $\mathcal{G}_k$, and show that a quotient of $\mathcal{H}_k^{\mathrm{ext}}$ is isomorphic to a large subalgebra of the centralizer. We further study the representation theory of $\mathcal{H}_k^{\mathrm{ext}}$ to find that the seminormal representations are indexed by a known family of partitions. The bases for the resulting modules are given by paths in a lattice of partitions, and the action of $\mathcal{H}_k^{\mathrm{ext}}$ is given by combinatorial formulas. 
\end{abstract}

\tableofcontents

\pagebreak

\section{Introduction}
\label{sec:intro}
The phenomenon now known as {Schur-Weyl duality} was first studied by Frobenius and Schur in their work connecting the representation theory of the symmetric groups and the general linear groups. It has since stimulated many advances in the development of \emph{tensor power centralizer algebras}, algebras of operators which preserve symmetries in a tensor space. Striking examples include:  
\begin{enumerate}
	\item the \emph{Brauer algebras} in \cite{Br}
		centralize the action of symplectic and orthogonal groups on tensor space $\left(\CC^n\right)^{\otimes k}$;
	 \item the \emph{graded Hecke algebra of type A} 
	 	centralizes the action of $\fsl_n$ on $L(\lambda) \otimes (\CC^n)^{\otimes k}$, where $L(\lambda)$ is the irreducible $\fsl_n$ module indexed by a partition $\lambda$ (see \cite{AS});
	 \item the \emph{degenerate affine Wenzl algebra} in \cite{Naz}
		 centralizes the action of symplectic and orthogonal groups on $L(\lambda) \otimes (\CC^n)^{\otimes k}$.
	\end{enumerate}
The paper of Orellana and Ram \cite{OR} provides a unified approach to studying tensor power centralizer algebras, including the \emph{affine and cyclotomic Hecke and  Birman-Murakami-Wenzl algebras}.

Recent work in the study of loop models and spin chains in statistical mechanics uncovered yet another potential use of Schur-Weyl duality in \cite{GN}. Specifically, a connection was discovered between the two-boundary Temperley-Lieb algebra and a quotient of the affine Hecke algebra of type C. The Temperley-Lieb algebra is the centralizer of the quantum group $\cU_q \fsl_2$ on tensor space  $M \otimes N \otimes (\CC^2)^{\otimes k}$, where $M$ and $N$ are simple $\cU_q\fsl_2$-modules, which suggested the possibility of constructing affine Hecke algebra type C modules explicitly using Schur-Weyl duality tools.

In Section \ref{sec:braid_group}, we begin the study of the centralizer of the action of $\fg$ on \tensorspace, where $\fg$ is a finite-dimensional complex reductive Lie algebra and  $M$, $N$, and $V$ are finite-dimensional irreducible $\fg$-modules. The new definition is that of the  \emph{degenerate two-boundary braid algebra} $\cG_k$, an associative algebra over the complex numbers. This braid algebra can be pictured as the degeneration of the quantum group analog, group algebra of the braid group in a space with two punctures, a generalization of the affine braid group studied in \cite{OR}. The algebra $\cG_k$ is generated by
	$$\CC[x_1, \dots, x_k] , ~\CC[y_1, \dots, y_k], ~\CC[z_0, z_1, \dots, z_k], \text{ and }  \CC S_k,$$
with relations twisting the polynomial rings and the symmetric group together. The first main theorem, Theorem \ref{thm:braid_group_rep}, is that $\cG_k$ acts on $M \otimes N \otimes V^{\otimes k}$ and that this action commutes with the action of $\fg$. In many cases, both historic and new, this action will produce $\End_\fg(M \otimes N \otimes V^{\otimes k})$. For example, 
\begin{enumerate}[\qquad (1)]
\item when $\fg = \fgl_n$ or $\fsl_n$, $V$ is the standard representation, and
\begin{enumerate}[(a)]
\item $M$ and $N$ are trivial, the image of $\cG_k$ in $\End_{\fg}(M \otimes N \otimes V^{\otimes k})$ is the same as that of the symmetric group $S_k$;
\item $M$ is trivial and $N$ is a simple highest weight module, the image of  $\cG_k$ in $\End_{\fg}(M \otimes N \otimes V^{\otimes k})$ is the same as that of the graded Hecke algebra of type A;
\end{enumerate}
\item when $\fg = \fso_n$ or $\fsp_{2n}$,  $V$ is the standard representation, and 
\begin{enumerate}
\item $M$ and $N$ are trivial,  the image of  $\cG_k$ in $\End_{\fg}(M \otimes N \otimes V^{\otimes k})$ is the same as that of the Brauer algebras;
\item $M$ is trivial and $N$  is a simple highest weight module, the  image of  $\cG_k$ in $\End_{\fg}(M \otimes N \otimes V^{\otimes k})$ is the same as that of the degenerate affine Wenzl algebra.
\end{enumerate}
\end{enumerate}
We discuss the specifics of these examples in Remark \ref{rk: other examples}.

In Section \ref{sec:Hecke}, we consider the new cases where $\fg = \fsl_n$ or $\fgl_n$, $M = L((a^p))$ (the finite-dimensional irreducible $\fg$-module indexed by the rectangular partition with $p$ parts of length $a$), $N= ((b^q))$, and $V$ is the standard representation. Theorem \ref{thm:Hecke-action} states that a twist of the representation given in Theorem \ref{thm:braid_group_rep} factors through this quotient. We call this quotient of $\cG_k$ the \emph{extended degenerate two-boundary Hecke algebra} $\cH_k^\ext$.

We further study the representation theory of $\cH_k^\ext$ throughout Sections \ref{sec:reps-of-H} and \ref{sec:seminormal_reps}, classifying the seminormal representations. Using the combinatorics of Young tableaux, we describe these representations explicitly in   Sections \ref{sec:seminormal_bases} and \ref{sec:seminormal_reps}. The basis elements for the resulting modules are given by paths in a lattice of partitions, and the action of $\cH_k^\ext$ is given in terms of contents of boxes in those partitions.

This work may proceed in a number of directions. Firstly,  an analogous theory may also be developed for centralizers of type B, C, and D, which will parallel that of the degenerate affine Wenzl algebra as studied in Nazarov in \cite{Naz} and \cite{AMR}.  Also, functorial techniques developed in \cite{OR} may be used to promote the study of calibrated $\cH_k^\ext$-modules, given in Section \ref{sec:seminormal_reps}, to that of all standard modules. This should extend to the study of standard modules for types B, C, and D as well.

Finally, one subalgebra of $\cH_k^\ext$, the \emph{degenerate two-boundary Hecke algebra} $\cH_k$, is of particular interest as it is strikingly similar to the graded Hecke algebra of type C. This can be seen through the combinatorics presented throughout Sections \ref{sec:seminormal_bases} and \ref{sec:seminormal_reps} and in the action of the type C Weyl group in the final proof of the paper. This observation suggests the possibility of studying representations of type C Hecke algebras using Schur-Weyl duality techniques, a study which is further developed in forthcoming papers.

\section{The degenerate two-boundary braid algebra}
\label{sec:braid_group}
%

Fix $k \in \ZZ_{\geq 0}$.  Let $S_k$ be the symmetric group, which is generated by simple transpositions $s_i = (i~ i+1)$ and braid relations
\begin{equation}
s_is_{i+1}{s_i} = {s_{i+1}} {s_i}{s_{i+1}} \qquad \text{ and } \qquad  {s_i}{s_j} = {s_j}{s_i}, \quad \text{ if } j \neq i \pm1. 
\label{eq:S_braidA}
\end{equation}
The main object of study in this section is the \emph{degenerate two-boundary braid algebra}, denoted $\cG_k$, and is a two-boundary analog to the degenerate one-boundary braid algebra in \cite{DRV1}. By design, we will see in the Section \ref{sec:braid-action} that $\cG_k$ acts on tensor space of a specific form for a finite-dimensional reductive Lie algebra. The algebra $\cG_k$ is generated over $\CC$ by 
\begin{equation}\label{generators}
t_\ww \quad \text{for } \ww\in S_k, \qquad x_1, \dots, x_k, \quad y_1, \dots, y_k, 
\quad \text{and} \quad  z_0, z_1, \dots, z_k, 
\end{equation}
subject to relations as follows. Let $m_1 = 0$ and, for $j>1$, define $m_{j} = \sum_{1\leq i <j}m_{i,j}$, where  
\begin{equation}\label{eq:m}
\begin{array}{ll}m_{j-1,j} = x_{j} - t_{s_{j-1}}x_{j-1}t_{s_{j-1}}, & \text{and} \\
 m_{i,j} = t_{(i~j-1)}m_{j-1, j} t_{(i~j-1)} & \text{for } 1\leq i < j-1.
 \end{array}\end{equation}
Then $\cG_k$ is the associative algebra generated over $\CC$ by elements \eqref{generators} with relations



%
%
%
%
\begin{align}
\label{rel:zdefn}
	&	z_i = x_i + y_i - m_i,& \qquad 1 \leq i \leq k,\\
	&	t_\ww t_{\ww'} = t_{\ww\ww'} & \quad \text{for } \ww \in S_k\\
	&	x_ix_j = x_j x_i, \quad y_iy_j=y_jy_i, \quad z_i z_j = z_j z_i
	& \text{ for all admissible } i,j\\
	&
		 t_{s_i} x_j = x_jt_{s_i} , \quad
		  t_{s_i}y_j = y_jt_{s_i}, \quad
		  t_{s_i}z_j = z_jt_{s_i},
  	& \mbox{for } j \neq i, i+1, \label{rel:graded_braid3}\\
	&\!\!\!\!		\begin{array}{c}
		 \left(z_0 + \cdots + z_i \right)x_j =  x_j \left(z_0 + \cdots + z_i\right),\\
		 \left(z_0 + \cdots + z_i \right)y_j =  y_j \left(z_0 + \cdots + z_i\right),
		 \end{array}
  	& \mbox{for }i \geq j, \label{rel:graded_braid_zsum}\\
	&
		  t_{s_i}(x_i + x_{i+1}) = (x_i + x_{i+1})t_{s_i},\quad
		  t_{s_i}(y_i + y_{i+1}) = (y_i + y_{i+1})t_{s_i},
	& \text{for } 1 \leq i \leq k-1, \label{rel:graded_braid4}\\
	&\!\!\!\!\begin{array}{l}
		  (t_{s_i} t_{s_{i+1}})\left(x_{i+1} - t_{s_i}x_it_{s_i}\right)(t_{s_{i+i}} t_{s_{i}}) = x_{i+2} - t_{s_{i+1}}x_{i+1}t_{s_{i+1}},\\ 
		  (t_{s_i} t_{s_{i+1}})\left(y_{i+1} - t_{s_i}y_it_{s_i}\right)(t_{s_{i+i}} t_{s_{i}}) = y_{i+2} - t_{s_{i+1}}y_{i+1}t_{s_{i+1}},\\ 
 	\end{array}
	& \text{for } 1 \leq i \leq k-2, \label{rel:graded_braid5}\\
	&x_{i + 1} - t_{s_i}x_it_{s_i} = y_{i + 1} - t_{s_i}y_it_{s_i}, & \text{for } 1 \leq i \leq k-1. \label{rel:graded_braid6}
\end{align}


\subsection{Action on tensor space}
\label{sec:braid-action}
Let $\fg$ be a finite-dimensional complex reductive Lie algebra. We fix a triangular decomposition 
	\begin{equation}\label{triangular}\qquad\fg = \fn^- \oplus \fh \oplus \fn^+, 
			 \quad \text{ with }  
			\fn^{+} = \bigoplus_{\alpha \in R^{+}} \fg_{\alpha}, 
			\end{equation}
and $R^+$ is a fixed set of positive roots for $\fg$.  
A \emph{weight} is an element of $\fh^* = \Hom(\fh, \CC)$.

The trace form $\<,\>: \fg \otimes \fg \to \CC$ associated to a faithful representation $\theta$ of $\fg$ is defined by 
	$$\< x, y \> = \tr(\theta(x) \theta(y)).$$
This is an $\ad$-invarient, symmetric, bilinear form which is nondegenerate on both $\fg$ and $\fh$. 
Therefore the map
 \begin{equation*}
 \qquad \qquad 
	\begin{array}{rcc}
		\fh &\longrightarrow& \fh^*\\
		h &\mapsto&\!\! \< h, \cdot \> \\
		h_{\mu} &\mapsto& \mu
	\end{array}
	\qquad \qquad \text{ is an isomorphism},
\end{equation*}
where $h_\mu$ is the unique element of $\fh$ such that $\< h_\mu, h\> = \mu(h)$ for all $h \in \fh$. 
We define the symmetric, bilinear, nondegenerate form $\<, \> : \fh^* \otimes \fh^* \to \CC$ by
$\< \lambda, \mu \> = \< h_{\lambda}, h_{\mu} \>$.

Let $M$, $N$, and $V$ be finite-dimensional simple $\fg$-modules, and consider  the action of $\fg$ on the tensor space $M \otimes N \otimes V^{\otimes k}$. Denote the \emph{centralizer} of the action of $\fg$ on a $\fg$-module $U$ by
	$$\End_{\fg}(U) = \{ \varphi \in \End(U) ~|~ x\varphi = \varphi x \text{ for all } x  \in \fg  \}.$$
We will construct a homorphism $\Phi: \cG_k \to \End_\fg(M \otimes N \otimes V^{\otimes k})$ using the observation that the map given by
\begin{equation}\label{eq:embedded_centralizers}
\begin{array}{rl}
\End_{\fg}(U) &\to \End_{\fg}(U \otimes U')\\
\varphi & \mapsto \varphi \otimes \id_{U'}
\end{array}
\end{equation}
is an injective algebra homomorphism for any $\fg$-modules $U$ and $U'$.

Fix a basis $\{ b_i \}$ for $\fg$ and let $\{b_i^*\}$ be the dual basis to $\{b_i\}$ with respect to $\<,\>$. The \emph{Casimir element} of the enveloping algebra $\cU\fg$ is
	\begin{equation}\label{eq:kappadefn}
	\kappa = \sum_{i} b_{i}b_{i}^*,
	\end{equation}
and is central in $\cU\fg$. If $U$ and $U'$ are $\fg$-modules, $\kappa$ acts on $U \otimes U'$ by 
	\begin{equation}\label{gammadefn}\kappa \otimes \id_{U'}+ \id_U \otimes \kappa + 2\gamma , \qquad\hbox{where}\quad \gamma = \sum_i b_i\otimes b_i^*.
	\end{equation}
Since $\kappa$ and $\gamma$ are independent of the choice of the basis, we have, in particular,
\begin{equation}\tag{\ref*{gammadefn}b}\label{gammadefnb}\gamma = \sum_i b_i^* \otimes b_i.\end{equation}
Let $\gamma_{j,j'}$ be the operator given by the action of $\gamma$ on the $j$ and $j'$ factors of $V^{\otimes k}$ (acting by the identity on all other factors). Note that $\gamma_{j,j'} = \gamma_{j',j}$ because of \eqref{gammadefnb}.
Similarly, for a factor $X$ ($X=M, N,$ $M\otimes N$, or $V$, where applied), denote by
	\begin{align*}
		\gamma_{M,N} & \qquad \text{$\gamma$ acting on factors $M$ and $N$ in a tensor space},\\
		\gamma_{X,i} & \qquad \text{$\gamma$ acting on factor $X$ and the $i^{\text{th}}$ copy of $V$ in a tensor space},\\
		\kappa_X & \qquad \text{$\kappa$ acting on the factor $X$ in a tensor space,}\\
		\kappa_{X,\leq j} &\qquad \text{$\kappa$ acting on the factor $X$ and the first $j$ factors of $V$},\\
				& \qquad \text{where $\kappa_{X,\leq 0} = \kappa_X$.}
	\end{align*}
Using  \eqref{gammadefn}  to apply $\kappa$ iteratively to $M \otimes V^{\otimes k}$,  $N\otimes V^{\otimes k}$, and  $M \otimes N\otimes V^{\otimes k}$,  we find that as operators on $M \otimes N \otimes V^{\otimes k}$, for $X = M, N$ or $M\otimes N$,
	\begin{equation}\label{eq:iterative-kappa-expansion}
	\kappa_{X,\leq j} = \kappa_X  + j\kappa_V +  2\left(\sum_{1 \leq i \leq j} \gamma_{X, i} + \sum_{1\leq r<s \leq j} \gamma_{r,s}\right).\end{equation}

\begin{thm}
\label{thm:braid_group_rep}
There is an algebra homomorphism 
$$\Phi\colon  \cG_k\to \End_{\fg}(M\otimes N\otimes V^{\otimes k})$$
defined by
$$\Phi(x_i) = \frac{1}{2}(\kappa_{M,\leq i} -\kappa_{M,\leq i-1}), \qquad 
\Phi(y_i) = \frac{1}{2}( \kappa_{N,\leq i} - \kappa_{N,\leq i-1}),$$
$$\text{ and } \quad \Phi(z_i) = \frac{1}{2}(\kappa_{M\otimes N,\leq i} -\kappa_{M\otimes N,\leq i-1} + \kappa_V),
\qquad \text{ for } 1 \leq i \leq k,$$
$$\Phi(z_0) = \frac{1}{2}(\kappa_{M \otimes N} - \kappa_M - \kappa_N) = \gamma_{M,N}, \qquad \text{and} $$
$$\Phi(t_{s_i}) = \id_M\otimes \id_N\otimes \id_V^{\otimes(j-1)}\otimes s \otimes \id_V^{\otimes (k-j-1)} \qquad \text{ for $1 \leq i \leq k-1$,}$$
where $s \cdot (u \otimes v) =  v \otimes u.$
\end{thm}

\begin{proof}

The $t_{s_i}$ act by simple transpositions, so they generate an action of $\CC S_k$ on $V^{\otimes k}$. Since the coproduct is cocommutative, the action of $\CC S_k$ commutes with the $\fg$-action.

Since $\kappa$ is central, $\kappa_{M,\leq i} \in \End_{\fg}(M \otimes V^{\otimes i})$. By  \eqref{eq:embedded_centralizers}, this means $\kappa_{M,\leq i}\otimes\id_V^{j-i}$ is an element of $\End_{\fg}(M \otimes V^{\otimes j})$ for $i<j$. So the actions of $\kappa_{M,\leq i}$, $i=1, 2, \dots, k$, and therefore the actions of $x_1, \dots, x_k$, pairwise commute. Similarly,  $\{y_1, \dots, y_k\}$ and $\{z_0, \dots, z_k\}$ each act commutatively on $M \otimes N \otimes V^{\otimes k}$. Again by \eqref{eq:embedded_centralizers}, these operators are also all contained in $\End_\fg(M \otimes N \otimes V^{\otimes k}).$ Moreover, since $M$, $N$, and $V$ are simple, $\kappa_M$, $\kappa_N$, and $\kappa_V$ act as constants. So
$$\Phi(z_0 + \cdots + z_i) = \frac{1}{2}( \kappa_{M\otimes N, \leq i} + i\kappa_V - \kappa_M - \kappa_N)$$
commutes with $\kappa_{M,\leq j}$ and $\kappa_{N,\leq j}$ for $j \leq i$, verifying \eqref{rel:graded_braid_zsum}. 
				
The relations in (\ref{rel:graded_braid4}) follow from 
\begin{align*} \Phi(t_{s_i} (x_i + x_{i+1}))
		&=   \half  t_{s_{i}}(\kappa_{M,\leq i+1} - \kappa_{M,\leq i-1})\\
		&=   \half t_{s_{i}}\left( \gamma_{M,i} +\gamma_{M,i+1} + 2 \kappa_V + 2 \sum_{\ell=1}^{i-1} (\gamma_{\ell,i} + \gamma_{\ell, i+1}) + 2 \gamma_{i, i+1} 
				 \right) & \text{ by \eqref{eq:iterative-kappa-expansion}, }\\
		&=   {\frac{1}{2}}\left( \gamma_{M,i+1} +\gamma_{M,i} + 2 \kappa_V + 2 \sum_{\ell=1}^{i-1} (\gamma_{\ell,i+1} + \gamma_{\ell, i}) + 2 \gamma_{i+1, i} 
				 \right)  t_{s_{i}} \\
		&= \Phi((x_i + x_{i+1})t_{s_i})
\end{align*}
(a similar computation confirms $\Phi(t_{s_i} (y_i + y_{i+1})) = \Phi((y_i + y_{i+1})t_{s_i})$). The action of the symmetric group commutes with the action of $\fg$, \and if $j<i$, $\kappa_{X,\leq j}$ acts by the identity on the $i$ and $i + 1$ factors of $V^{\otimes k}$. Thus
 	\begin{equation}\label{rel:graded_braid3_alt}
	t_{s_i} \kappa_{X,\leq j} = \kappa_{X,\leq j} t_{s_i}, \textrm{ if $j \neq i$, and $X = M$, $N$, or $M \otimes N$},\end{equation}
is satisfied for all $i \neq j$, which implies \eqref{rel:graded_braid3}.

Finally, as operators on $M \otimes N \otimes V^{\otimes k}$ via $\Phi$,
\begin{align}
x_i 	= \half \left( \kappa_{M, \leq i} - \kappa_{M, \leq i-1} \right)\nonumber
	=\half\kappa_V+ \gamma_{M,i} + \sum_{1\le \ell<i} \gamma_{\ell, i}, \nonumber
	&& \text{ by \eqref{eq:iterative-kappa-expansion}}
\end{align}
and similarly
$$
y_i =\half \kappa_i + \gamma_{N,i} + \sum_{1\le \ell<i} \gamma_{\ell, i},
	\quad \text{ and } \quad
z_i = \kappa_V + \gamma_{N,i} + \gamma_{M,i} + \sum_{1\le \ell<i} \gamma_{\ell, i}.
$$
So
\begin{equation}\label{tonfactors}
m_{i, i+1} = x_{i+1} - t_{s_i}x_it_{s_i} = y_{i+1} - t_{s_i}y_it_{s_i} = z_{i+1} - t_{s_i}z_it_{s_i}  = \gamma_{i,i+1}.
\end{equation}
So  \eqref{rel:graded_braid6} and \eqref{rel:zdefn} are satisfied. Since $t_{s_i}t_{s_{i+1}} \gamma_{i,i+1}t_{s_{i+1}}t_{s_i} = t_{s_i}\gamma_{i,i+2}t_{s_i} = \gamma_{i+1,i+2}$, relation 
\eqref{rel:graded_braid5} follows from   \eqref{tonfactors}.

\end{proof}

\begin{remark}\label{rk: other examples}
As discussed in the introduction, the degenerate two-boundary braid algebra is meant to be the degeneration of the group algebra of the two-boundary braid group, which is the braid group in a space with two punctures, or ``flag poles.'' The two-boundary braid group is the generalization of the affine braid group used in \cite{OR}, and just like the affine braid group, has many centralizer algebras  for quantum groups as quotients. Analogously, $\cG_k$ has many familiar centralizer algebras for Lie algebras as quotients, and the map in Theorem \ref{thm:braid_group_rep} factors through these quotients (in some cases after applying an automorphism). For example:
\begin{enumerate}[\qquad (1)]
\item When $\fg = \fgl_n$ or $\fsl_n$, $V$ is the standard representation, 
the action of $t_{s_i}$ on $V \otimes V$ is the same as that of $\gamma$ on $V \otimes V$. So
\begin{enumerate}[(a)]
\item when $M$ and $N$ are trivial, the images of $x_i$ and $y_i$ are linear combinations of the images of $t_{s_j}$, $j=1, \dots, i-1$, and so the image of $\cG_k$ in $\End(M \otimes N \otimes V^{\otimes k})$ is the same as that of the symmetric group $S_k$; and
\item when $M$ is trivial and $N$ is a simple highest weight module, the image of $x_i$ is redundant as above, and the image of  $\cG_k$ in $\End(M \otimes N \otimes V^{\otimes k})$ (after a version of the automorphism in Lemma \ref{lem:G-auto} when $\fg=\fsl_n$) is the same as that of the graded Hecke algebra of type A. 
\end{enumerate}
\item When $\fg = \fso_n$ or $\fsp_{2n}$,  $V$ is the standard representation, and $M$ is trivial, we hope to see the Brauer algebra in \cite{Br} and the degenerate affine Wenzl algebra in \cite{Naz} (when $N$ is trivial or not, respectively), in which case we expect to see elements which are diagrammatically represented by 
$$s_i = 
\begin{picture}(73,12)(0,6) 
	\put(0,0){\line(0,1){15}}
	\put(20,0){\line(0,1){15}}
	\put(55,0){\line(0,1){15}}
	\put(75,0){\line(0,1){15}}
	\put(30,0){\line(1,1){15}}
	\put(45,0){\line(-1,1){15}}
	\multiput(-2,-2)(0,15){2}{\tiny $\bullet$}
	\multiput(18,-2)(0,15){2}{\tiny $\bullet$}
	\multiput(28,-2)(0,15){2}{\tiny $\bullet$}
	\multiput(43,-2)(0,15){2}{\tiny $\bullet$}
	\multiput(53,-2)(0,15){2}{\tiny $\bullet$}
	\multiput(72.5,-2)(0,15){2}{\tiny $\bullet$}	
	\multiput(-2,-10)(0,29.5){2}{\tiny $1$}
	\multiput(18,-10)(0,29.5){2}{\tiny $$}
	\multiput(28,-10)(0,29.5){2}{\tiny $i$}
	\multiput(40,-10)(0,29.5){2}{\tiny $i\!+\!1$}
	\multiput(53,-10)(0,29.5){2}{\tiny $$}
	\multiput(72,-10)(0,29.5){2}{\tiny $k$}
	\multiput(3,5)(56,0){2}{\tiny$\cdots$}
\end{picture}
\qquad \text{ and } \qquad \bar s_i = 
\begin{picture}(76,12)(0,6) 
	\put(0,0){\line(0,1){15}}
	\put(20,0){\line(0,1){15}}
	\put(55,0){\line(0,1){15}}
	\put(75,0){\line(0,1){15}}
	\put(37.5,0){\oval(15,8)[t]}
	\put(37.5,15){\oval(15,8)[b]}
	\multiput(-2,-2)(0,15){2}{\tiny $\bullet$}
	\multiput(18,-2)(0,15){2}{\tiny $\bullet$}
	\multiput(28,-2)(0,15){2}{\tiny $\bullet$}
	\multiput(43,-2)(0,15){2}{\tiny $\bullet$}
	\multiput(53,-2)(0,15){2}{\tiny $\bullet$}
	\multiput(72.5,-2)(0,15){2}{\tiny $\bullet$}	
	\multiput(-2,-10)(0,29.5){2}{\tiny $1$}
	\multiput(18,-10)(0,29.5){2}{\tiny $$}
	\multiput(28,-10)(0,29.5){2}{\tiny $i$}
	\multiput(40,-10)(0,29.5){2}{\tiny $i\!+\!1$}
	\multiput(53,-10)(0,29.5){2}{\tiny $$}
	\multiput(72,-10)(0,29.5){2}{\tiny $k$}
	\multiput(3,5)(56,0){2}{\tiny$\cdots$}
\end{picture}~.\phantom{\Bigg|}
$$
The diagram $s_i$ corresponds to $t_{s_i}$ in $\cG_k$, and the diagram $\bar s_i$ corresponds to the element
$$e_i = t_{s_i} y_i - y_{i+1} t_{s_i} -1.$$ The map $\Phi$ factors through the quotient of $\cG_k$ by the relations in \cite[\S 4]{Naz} (after an automorphism in the case where $\fg = \fsp_{2n}$) or by relations in \cite[\S 2.2]{DRV1} (with no automorphism), as is shown in \cite{DRV2}. So, 
\begin{enumerate}
\item when $M$ and $N$ are trivial,  the image of  $\cG_k$ in $\End(M \otimes N \otimes V^{\otimes k})$ is the same as that of the Brauer algebras, and 
\item when $M$ is trivial and $N$  is a simple highest weight module, the  image of  $\cG_k$ in $\End(M \otimes N \otimes V^{\otimes k})$ is the same as that of the degenerate affine Wenzl algebra.
\end{enumerate}
\end{enumerate}

\end{remark}


\section{The degenerate two-boundary Hecke algebra}
\label{sec:Hecke}

Our next goal is to consider the case where $\fg$ is of type  $\fgl_n$ or $\fsl_n$, and  $M$, $N$, and $V$ are  three specific $\fg$-modules. In general, even if we specify  $V$ to be the standard representation as usual, the decomposition of  $M \otimes N$ is not in general multiplicity free, and so the method of considering quotients of the braid algebra in studying centralizers $\End_\fg(M \otimes N \otimes V^{\otimes k})$ is ineffective. However, in the case where $M$ and $N$ are indexed by rectangular partitions, it is an amazing consequence of the Littlewood-Ricardson rule that the decomposition of $M \otimes N$ is multiplicity free.  Furthermore, when constructing Hecke algebras in the quantum case, one places quadratic relations on all generators corresponding to the $\cR$-matrices. In the degenerate case, these generators are specifically $t_{s_1}, \dots, t_{s_{k-1}}$, $x_1$ and $y_1$.  By choosing $M$ and $N$ to be indexed by rectangular partitions, we will force quadratic relations on $x_1$ and $y_1$ as desired.

In this way, we use the representations of $\cG_k$ in Theorem \ref{thm:braid_group_rep} to motivate the construction of a new algebra, the degenerate extended two-boundary Hecke algebra. In Section \ref{sec:reps-of-H} we will carefully lay out the combinatorics behind this construction and explore this motivation further. This section is devoted to the definition and two presentations of $\cH^\ext_k$.

~

Fix $a,b,p,q \in \ZZ_{>0}$. The \emph{degenerate extended two-boundary Hecke algebra} $\cH_k^\text{ext}$ is the quotient of the degenerate two-boundary braid algebra  by the relations
	\begin{equation}\label{rel:hecke}
	t_{s_i} x_i = x_{i+1} t_{s_i}  -1, \quad t_{s_i} y_i = y_{i+1} t_{s_i} -1, \quad i=1, \dots, k-1,
	\end{equation}	
	\begin{equation}\label{rel:hecke2}
	(x_1-a)(x_1+p) = 0, \qquad (y_1 -b)(y_1 + q) = 0.
	\end{equation}
The \emph{degenerate two-boundary Hecke algebra} $\cH_k$ is the subalgebra of $\cH_k^\text{ext}$ generated by $x_1, \ldots, x_k$, $y_1, \ldots, y_k$, $z_1, \ldots, z_k$, $t_{s_1}, \ldots, t_{s_{k-1}}$.

Proposition \ref{thm:hecke_ext_presentation} provides a presentation of $\cH^\ext_k$ which is a consolidation of the presentation of $\cG_k$ using the quotient in \eqref{rel:hecke} and \eqref{rel:hecke2}. We follow this up with Theorem \ref{thm:hecke_ext_presentation-short}, which provides a much more efficient presentation which we will make use of in Section \ref{sec:seminormal_reps}. 

\begin{prop}\label{thm:hecke_ext_presentation}
Define 
$$x_{i} = t_{s_{i-1}}x_{i-1}t_{s_{i-1}} + t_{s_{i-1}}, \quad \quad z_{i} = t_{s_{i-1}}z_{i-1}t_{s_{i-1}} + t_{s_{i-1}},\quad \text{ for $i=2, \ldots, k$,}$$	
$$m_1 = 0,\qquad  m_i = \sum_{j=1}^{i-1} t_{(j~i)} \text{ for $i>0$}
\quad \text{ and} \quad 
y_i = z_i - x_i + m_i \quad \text{for $i=1, \ldots, k$.}$$
Then $\cH_k^\text{ext}$ is generated as an algebra over $\CC$ by $x_1, z_0, z_1$ and $t_\ww$ for $\ww \in S_k$ with relations $$t_{\ww} t_{\ww'} = t_{\ww \ww'}\qquad \text{ for } \ww, \ww' \in S_k,$$ and
\begin{enumerate}
\item[] Quadratic relations:
	$$(x_1-a)(x_1+p) = 0, \qquad (y_1 -b)(y_1 + q) = 0, \qquad a,b,p,q \in \ZZ_{>0},$$
\item[] Commutation relations:
	$$t_{s_i} x_j = x_jt_{s_i} , ~
		  t_{s_i}z_j = z_jt_{s_i},
 	\quad \mbox{for } j \neq i, i+1, \qquad $$	
	$$x_i x_j = x_jx_i , ~
		  y_iy_j = y_j y_i, ~
		  z_iz_j = z_j z_i, ~
		  z_0 z_i = z_i z_0, ~
 	\quad \mbox{for } 1 \leq i,j \leq k, \qquad $$
	$$x_jz _i = z_i x_j, \quad \text{ for } i>j.$$

\item[] Twisting relations:
$$\begin{array}{c}
		 x_i (z_0 + \cdots + z_i) = (z_0 + \cdots + z_i)x_i,\\
		 y_i (z_0 + \cdots + z_i) = (z_0 + \cdots + z_i)y_i,
		 \end{array}
 	\quad \mbox{for }i =1 \dots k. $$
\end{enumerate}
\end{prop}

\begin{proof}
Equation \eqref{eq:m} can be rewritten as
	$$m_{j,j+1} = x_{j+1} - t_{s_j} x_j t_{s_j} = t_{s_j} $$ $$\text{and} \qquad
		m_{i,j} = t_{(i~j-1)}m_{j-1, j} t_{(i~j-1)}= t_{(i~j)} \quad  \text{ if }  i <  j-1.$$
So $$m_1 = 0, \quad m_i = \sum_{1<j<i} t_{(i~j)}.$$ Therefore \eqref{rel:zdefn} implies
	\begin{align*}
		t_{s_i} z_i t_{s_i} &= t_{s_i} (x_i + y_i - m_i) t_{s_i} \\
			&= x_{i+1} - t_{s_i} + y_{i+1} - t_{s_i} - t_{s_i}\left( \sum_{1<j<i} t_{(i~j)}\right) t_{s_i}\\
			&=x_{i+i} + y_{i+1} - t_{s_i} - t_{s_i} -  \sum_{1<j<i} t_{(i+1~j)}\\
			 &= x_{i+i} + y_{i+1} - m_{i} - t_{s_i}\\
			 &= z_{i+1} - t_{s_i}.
	\end{align*}
Similarly, any two of
	$$x_{i+1} - t_{s_i} x_i t_{s_i}=t_{s_i}, \quad
		y_{i+1} - t_{s_i} y_i t_{s_i}=t_{s_i}, \quad \text{ and } \quad
		z_{i+1} - t_{s_i} z_i t_{s_i}=t_{s_i}, \quad \quad i=1, \dots, k-1, $$
		imply the third.
So we use \eqref{rel:hecke} to discard the generators $x_2, \dots, x_k$, $y_1, \dots, y_k$, and $z_2, \dots, z_k$,  by defining
$$x_{i} = t_{s_{i-1}}x_{i-1}t_{s_{i-1}} + t_{s_{i-1}}, \quad z_{i} = t_{s_{i-1}}z_{i-1}t_{s_{i-1}} + t_{s_{i-1}},
	\quad i=2, \ldots, k,$$
$$\text{and } \quad y_i = z_i - x_i + m_i, \quad i=1, \ldots, k$$

Relation \eqref{rel:graded_braid4} can be rewritten as
	 $$t_{s_i}x_i -  x_{i+1}t_{s_i}= t_{s_i} (t_{s_i}x_i -  x_{i+1}t_{s_i})t_{s_i}  \quad \text{and} \quad  t_{s_i}y_i -  y_{i+1}t_{s_i} = t_{s_i}(t_{s_i}y_i -  y_{i+1}t_{s_i})t_{s_i}$$
for $1 \leq i \leq k-1$, which is equivalent to $-1 = -1$. Relation \eqref{rel:graded_braid5} is equivalent to
	$$\begin{array}{l}
		  (t_{s_i} t_{s_{i+1}})\left(t_{s_i}\right)(t_{s_{i+i}} t_{s_{i}}) = t_{s_{i+1}}
	\end{array} \quad \text{for } 1 \leq i \leq k-2,$$
which is redundant with relation $t_{s_i}^2 = 1$ and the second relation in \eqref{eq:S_braidA}. 
Relation \eqref{rel:graded_braid6} is equivalent to $t_{s_i} = t_{s_i}$. So by introducing \eqref{rel:hecke}, we can discard relations \eqref{rel:graded_braid4} - \eqref{rel:graded_braid6}. The second relation in \eqref{rel:graded_braid3} can also be discarded since for $j \neq i, i+1$,
	\begin{align*}
		 t_{s_i}y_j &=  t_{s_i}\left(z_j - x_j + \sum_{\ell=1}^{j-1} t_{(\ell~j)}\right)
		 	= \left(z_j - x_j + \sum_{\ell=1}^{j-1} t_{(\ell~j)}\right) t_{s_i}
			= y_j  t_{s_i},
	\end{align*}

Finally, independent of  \eqref{rel:hecke}, we rewrite relation \eqref{rel:graded_braid_zsum} as 
	$$\begin{array}{c}
		 x_i z_0 = z_0 x_i + \big((z_1 + \cdots + z_i)x_i - x_i (z_1 + \cdots + z_i)\big),\\
		 y_i z_0 = z_0 y_i + \big((z_1 + \cdots + z_i)y_i - y_i (z_1 + \cdots + z_i)\big),
		 \end{array}
 	\quad \mbox{for }i =1, \dots, k, $$
and 
	\begin{align*}
	x_j z_i &= x_j(z_0 + \cdots + z_i) - x_j (z_0 + \cdots + z_{i-1})\\
		&= (z_0 + \cdots + z_i)x_j - (z_0 + \cdots + z_{i-1})x_j\\
		&= z_i x_j, \quad \text{ and } & 	\text{ for $i>j$}.\\
  	y_j z_i &= z_i y_j.
		\end{align*}
\end{proof}

The following is a streamlined version of Proposition \ref{thm:hecke_ext_presentation}, which will be our favorite presentation for calculating representations in Section \ref{sec:seminormal_reps}.

\begin{thm}\label{thm:hecke_ext_presentation-short} 
Let $$w_i = z_i - \frac12(a-p+b-q).$$
$\cH_k^\ext$ is generated as an algebra over $\CC$ by  $w_0$, $w_1$, \dots, $w_k$, $x_1$, $t_{s_1}, \dots, t_{s_{k-1}}$ with relations 
\begin{enumerate}
\item[]Braid relations:
	\begin{equation}\label{eq:S_braid} t_{s_i}^2 = 1, \qquad t_{s_i}t_{s_{i+1}}t_{s_i} = t_{s_{i+1}}t_{s_{i}}t_{s_{i+1}}, \qquad
t_{s_i}t_{s_j} = t_{s_j} t_{s_i} \quad \text{ for } |i-j|>1,\end{equation}
	\begin{align}
	x_1(t_{s_1}x_1t_{s_1} + t_{s_1}) &= (t_{s_1}x_1t_{s_1} + t_{s_1}) x_1,\label{eq:braid 4}
	\end{align}
\item[]Quadratic relation:
	\begin{equation}
	(x_1 - a)(x_1 + p) = 0, \label{eq:quadratic}
	\end{equation}
\item[]Commutation relations:
 \begin{align}
 	t_{s_i}w_j &= w_jt_{s_i},& j \neq i, i+1, 	\label{eq:comm-tw}\\
	x_1 w_i &= w_i x_1,& i=2, \dots, k,		\label{eq:comm-xw}\\
	x_1 t_{s_i} &= t_{s_i} x_1,& i=2, \dots, k-1,	\label{eq:comm-xt}\\
	w_iw_j &= w_j w_i,& i,j = 0, \dots, k,		\label{eq:comm-ww}
	\end{align}
 \item[] Twisting relations:
	\begin{align}
	t_{s_i} w_i &= w_{i+1} t_{s_i}  -1,& \quad i=1, \dots, k-1,\label{eq:twist-tw}\\
	x_1w_0 &= w_0x_1 - (x_1 w_1 - w_1 x_1),\label{eq:twist-xw0}
	\end{align}
	\begin{equation}\label{eq:xw}
	x_1 w_1 =  -w_1 x_1+ (a-p)w_1 + w_1^2 + \left(\frac{a+p + b+ q}{2}\right) \left(\frac{a+p - (b+q)}{2}\right).
	\end{equation}
\end{enumerate}
\end{thm}
\begin{proof}
With the exception of  $$y_i = z_i - x_i + m_i = w_i - x_i + m_i + \frac12(a-p+b-q),$$
every substitution of $z_i = w_i + \half(a-p+b-q)$ in the presentation in Proposition  \ref{thm:hecke_ext_presentation} results in a cancelation,  i.e. 
		$$  t_{s_i}w_j = w_jt_{s_i},
 	\quad \mbox{for } j \neq i, i+1, \qquad\qquad	  w_iw_j = w_j w_i, ~
		  w_0 w_i = w_i w_0, ~
 	\quad \mbox{for } 1 \leq i,j \leq k, \qquad $$
	$$x_jw _i = w_i x_j \quad \text{ for } i>j \qquad\qquad \text{ and } \qquad\qquad \begin{array}{c}
		 x_i (w_0 + \cdots + w_i) = (w_0 + \cdots + w_i)x_i,\\
		 y_i (w_0 + \cdots + w_i) = (w_0 + \cdots + w_i)y_i,
	\end{array}$$
are immediate.

\noindent Next, we address \eqref{eq:xw} by proving the following claim: 
\begin{enumerate}
\item[]\textbf{Claim 1:} The set of relations \\
(A): $\quad (x_1 -a)(x_1 + p) = 0,$ 
$(y_1 - b)(y_1 +q) = 0,$  and  
$w_1 = x_1 + y_1 - \half(a-p+b-q)$\\~\\
are equivalent to the set of relations \\~\\
(B): $\qquad \begin{matrix} (x_1 -a)(x_1 + p) = 0 \quad \textrm{ and } \quad \\
x_1 w_1 =  -w_1 x_1+ (a-p)w_1 + w_1^2 + \left(\frac{a+p + b+ q}{2}\right) \left(\frac{a+p - (b+q)}{2}\right)\end{matrix}
$
\\~\\
\emph{Proof:}
\begin{enumerate}
\item[]\textbf{(A)$\implies$(B):} First notice that
	$$x_1^2 = (a-p)x_1 + ap, \quad y_1^2 = (b-q) y_1 + bq,$$
	 $$\quad \text{ and } \quad
	z_1^2 = (x_1 + y_1)^2 = x_1y_1 + y_1 x_1  + (a-p)x_1 + (b-q)y_1 + ap+bq.$$
So
	\begin{align*}
	x_1 w_1 + w_1 x_1 &= x_1 ( x_1 + y_1  - (a-p + b-q)/2) \\
	&~~+ ( x_1 + y_1  - (a-p + b-q)/2) x_1\\
		&= 2 x_1^2 + (x_1 y_1 + y_1 x_1) - (a - p + b - q)x_1\\
		&= (a-p - (b-q))x_1 + 2ap + (x_1 y_1 + y_1 x_1)  .
	\end{align*}
Since
	\begin{align*}
		w_1^2 &= z_1^2 - (a-p+b-q)z_1 + \frac{1}{4}(a-p+b-q)^2\\
			&=x_1y_1 + y_1 x_1  + (a-p)x_1 + (b-q)y_1 + ap+bq\\
			&\qquad  -  (a-p+b-q)(x_1 + y_1)  + \frac{1}{4}(a-p+b-q)^2\\
			&= (x_1y_1 + y_1 x_1) - (b-q)x_1 - (a-p) (w_1 -x_1 + (a-p+b-q)/2) \\
			&\qquad + ap+bq + \frac{1}{4}(a-p+b-q)^2\\
			&= (x_1y_1 + y_1 x_1) + (a-p - (b-q)) x_1 - (a-p)w_1\\
			&\qquad  + ap+bq  - (a-p)^2/4 + (b-q)^2/4
	\end{align*}
we have
	\begin{align*}
	x_1 w_1 + w_1 x_1
		&= (a-p - (b-q))x_1 +2 ap \\
		&\qquad	+ (w_1^2 - \big( (a-p - (b-q)) x_1 - (a-p)w_1 \\
		&\qquad+ ap+bq  - (a-p)^2/4 + (b-q)^2/4)\big)\\
		&=  w_1^2 + (a-p)w_1+ \left(\frac{a+p + b+ q}{2}\right) \left(\frac{a+p - (b+q)}{2}\right).\\
	\end{align*}
\item[] \textbf{(B)$\implies$(A):} If $y_1 = w_1 - x_1 +\half(a-p+b-q)$, then using both relations in (B) to expand $x_1^2$ and $w_1 x_1 + x_1 w_1$, direct calculation yields 
\begin{align*}
	(y_1 - b)& (y_1 + q)\\
	& = (w_1 - x_1+ \half(a-p+b-q) - b)(w_1 - x_1+ \half(a-p+b-q)  + q)\\
	 &= 0.
\end{align*}
\end{enumerate}
\end{enumerate}

\noindent The remainder is showing that the relations in Proposition \ref{thm:hecke_ext_presentation} follow from relations \eqref{eq:S_braid} and  \eqref{eq:braid 4}-\eqref{eq:xw}. As in Proposition \ref{thm:hecke_ext_presentation}, define $x_{i+1} = t_{s_i}x_it_{s_i} + t_{s_i}$. By induction on $\ell$, 	
\begin{align}x_{i+1} &= t_{s_{i}}\cdots t_{s_{\ell+1}}(t_{s_{\ell}}x_\ell t_{s_{\ell}}+ t_{s_{\ell}})t_{s_{\ell +1}}\cdots t_{s_{i}}\nonumber\\
		&\qquad + \sum_{r = \ell +1}^{i} t_{s_{i}}\cdots t_{s_{r+1}}t_{s_{r}}t_{s_{r+1}}\cdots t_{s_{i}}. \label{eq:x_induct}
	\end{align}

\begin{enumerate}
\item[]\textbf{Claim 2:} $t_{s_{i}}x_j = x_j t_{s_{i}}$ for $i > j$.\\
\noindent\emph{Proof:} If $i>j$, then $t_{s_i}$ commutes with $t_{s_\ell}$ for all $\ell < j$, so by \eqref{eq:comm-xt} and \ref{eq:x_induct}
	\begin{align*}
		t_{s_i} x_j &= t_{s_i}(t_{s_{j-1}}\cdots t_{s_{2}})(t_{s_{1}}x_1t_{s_{1}}+ t_{s_{1}})(t_{s_{2}}\cdots t_{s_{j-1}})\\
		&\qquad + t_{s_i}\sum_{\ell = 2}^{j-1} t_{s_{j-1}}\cdots t_{s_{\ell+1}}t_{s_{\ell}}t_{s_{\ell+1}}\cdots t_{s_{j-1}} \\		&= (t_{s_{j-1}}\cdots t_{s_{2}})(t_{s_{1}}x_1t_{s_{1}}+ t_{s_{1}})(t_{s_{2}}\cdots t_{s_{j-1}})t_{s_i}\\
		&\qquad + \left(\sum_{\ell = 2}^{j-1} t_{s_{j-1}}\cdots t_{s_{\ell+1}}t_{s_{\ell}}t_{s_{\ell+1}}\cdots t_{s_{j-1}} \right)t_{s_i}\\
		&= x_jt_{s_i} .
	\end{align*}
\item[]\textbf{Claim 3:} $t_{s_{i}}x_j = x_j t_{s_{i}}$ for $i<j-1$.\\
\noindent\emph{Proof:} 
By \eqref{eq:x_induct}, 
\begin{align*}
	t_{s_{i}}x_j &= t_{s_i}(t_{s_{j-1}}\cdots t_{s_{i+2}} t_{s_{i+1}} )(t_{s_{i}}x_it_{s_{i}}+ t_{s_{i}})(t_{s_{i+1}}t_{s_{i+2}}\cdots t_{s_{j-1}})\\
		&\qquad + t_{s_i} t_{s_{j-1}}\cdots t_{s_{i+2}}t_{s_{i+1}}t_{s_{i+2}}\cdots t_{s_{j-1}} \\
		&\qquad + t_{s_i}\sum_{\ell = i+2}^{j-1} t_{s_{j-1}}\cdots t_{s_{\ell+1}}t_{s_{\ell}}t_{s_{\ell+1}}\cdots t_{s_{j-1}} \\
		&= (t_{s_{j-1}}\cdots t_{s_{i+2}})(t_{s_i})\left( t_{s_{i+1}} (t_{s_{i}}x_it_{s_{i}}+ t_{s_{i}})t_{s_{i+1}}+t_{s_{i+1}}\right)(t_{s_{i+2}}\cdots t_{s_{j-1}})\\
		&\qquad + \left(\sum_{\ell = i+2}^{j-1} t_{s_{j-1}}\cdots t_{s_{\ell+1}}t_{s_{\ell}}t_{s_{\ell+1}}\cdots t_{s_{j-1}}\right)t_{s_i} \\
\end{align*}
But, by Claim 2, since $i+1 > i$, 
\begin{align*}
t_{s_i}( t_{s_{i+1}} (t_{s_{i}}x_it_{s_{i}}+& t_{s_{i}})t_{s_{i+1}}+t_{s_{i+1}}) \\
	&= 
t_{s_i} t_{s_{i+1}} t_{s_{i}}x_it_{s_{i}}t_{s_{i+1}}+ t_{s_i} t_{s_{i+1}}t_{s_{i}}t_{s_{i+1}}+t_{s_{i}}t_{s_{i+1}}\\
&=
t_{s_{i+1}} t_{s_{i}} t_{s_{i+1}}x_it_{s_{i}}t_{s_{i+1}}+ t_{s_{i+1}}t_{s_{i}}t_{s_{i+1}}^2+t_{s_{i}}t_{s_{i+1}}t_{s_i}^2\\
&=
t_{s_{i+1}} t_{s_{i}} x_i t_{s_{i+1}}t_{s_{i}}t_{s_{i+1}}+ t_{s_{i+1}}t_{s_{i}}+t_{s_{i+1}}t_{s_{i}}t_{s_{i+1}}t_{s_i}\\
&=
t_{s_{i+1}} t_{s_{i}} x_i t_{s_{i}}t_{s_{i+1}}t_{s_{i}}+ t_{s_{i+1}}t_{s_{i}}+t_{s_{i+1}}t_{s_{i}}t_{s_{i+1}}t_{s_i}\\
&=( t_{s_{i+1}} (t_{s_{i}}x_it_{s_{i}}+ t_{s_{i}})t_{s_{i+1}}+t_{s_{i+1}}) t_{s_i}.
\end{align*}
So $t_{s_{i}}x_j  = x_j t_{s_{i}}$.

\item[]\textbf{Claim 4:} $x_i x_{j} = x_{j}x_i$ for $i,j = 1, \dots, k$.\\
\noindent\emph{Proof:} First, $x_1 x_2 = x_2 x_1$ by \ref{eq:braid 4}. Next, we induct on  $i$ to show $x_i x_{i+1} = x_{i+1}x_i$ for $i = 1, \dots, k-1$: 
	\begin{align*}
			x_i &x_{i+1}
				 =x_it_{s_{i}}x_it_{s_{i}} + x_it_{s_{i}}\\
				&= \big((t_{s_{i-1}}x_{i-1}t_{s_{i-1}} + t_{s_{i-1}})t_{s_{i}}(t_{s_{i-1}}x_{i-1}t_{s_{i-1}} + t_{s_{i-1}}) + (t_{s_{i-1}}x_{i-1}t_{s_{i-1}} + t_{s_{i-1}})\big)t_{s_{i}}\\
				&= (t_{s_{i-1}}x_{i-1}t_{s_{i-1}}t_{s_{i}} t_{s_{i-1}}x_{i-1}t_{s_{i-1}} 
					+ t_{s_{i-1}}x_{i-1}t_{s_{i-1}} )t_{s_{i}} 
				\\&\qquad
					 +  (t_{s_{i-1}}t_{s_{i}}t_{s_{i-1}}x_{i-1}t_{s_{i-1}} 
					+ t_{s_{i-1}}x_{i-1}t_{s_{i-1}}t_{s_{i}} t_{s_{i-1}})t_{s_{i}} 
				\\
				&\qquad 					+( t_{s_{i-1}}t_{s_{i}}t_{s_{i-1}}
					+ t_{s_{i-1}} )t_{s_{i}}.
	\end{align*}
But
\begin{align*}
 (t_{s_{i-1}}x_{i-1}t_{s_{i-1}}t_{s_{i}} &t_{s_{i-1}}x_{i-1}t_{s_{i-1}} 
					+ t_{s_{i-1}}x_{i-1}t_{s_{i-1}} )t_{s_{i}}
				\\&= t_{s_{i-1}}x_{i-1}t_{s_{i}}t_{s_{i-1}} t_{s_{i}}x_{i-1}t_{s_{i-1}}t_{s_{i}} 
					+ t_{s_{i-1}}x_{i-1}t_{s_i}^2t_{s_{i-1}} t_{s_{i}} 
\\
				&= t_{s_{i-1}}t_{s_{i}}x_{i-1}t_{s_{i-1}} x_{i-1} t_{s_{i}}t_{s_{i-1}}t_{s_{i}}
					+ t_{s_{i-1}}t_{s_i}x_{i-1}t_{s_i}t_{s_{i-1}} t_{s_{i}} 
\\
				&= t_{s_{i-1}}t_{s_{i}}x_{i-1}t_{s_{i-1}} x_{i-1} t_{s_{i-1}}t_{s_{i}}t_{s_{i-1}}
					+ t_{s_{i-1}}t_{s_i}x_{i-1}t_{s_{i-1}}t_{s_{i}} t_{s_{i-1}} 
\\
				&= t_{s_{i-1}}t_{s_{i}}(x_{i-1}t_{s_{i-1}} x_{i-1} t_{s_{i-1}} + x_{i-1}t_{s_{i-1}})t_{s_{i}}t_{s_{i-1}} 
\\
				&= t_{s_{i-1}}t_{s_{i}}(t_{s_{i-1}} x_{i-1} t_{s_{i-1}}x_{i-1} + t_{s_{i-1}}x_{i-1})t_{s_{i}}t_{s_{i-1}} 
				\\
				&= t_{s_{i}}(t_{s_{i-1}}x_{i-1}t_{s_{i-1}}t_{s_{i}} t_{s_{i-1}}x_{i-1}t_{s_{i-1}} 
					+ t_{s_{i-1}}x_{i-1}t_{s_{i-1}} ),
\end{align*}
\begin{align*}
  (t_{s_{i-1}}t_{s_{i}}t_{s_{i-1}}x_{i-1}t_{s_{i-1}} 
					&+ t_{s_{i-1}}x_{i-1}t_{s_{i-1}}t_{s_{i}} t_{s_{i-1}})t_{s_{i}} 
					\\
	&=   t_{s_{i}}t_{s_{i-1}}x_{i-1}t_{s_{i}}t_{s_{i-1}} + t_{s_{i-1}}x_{i-1}t_{s_{i}}t_{s_{i-1}}\\
	&=   t_{s_{i}}t_{s_{i-1}}x_{i-1}t_{s_{i}}t_{s_{i-1}} + t_{s_i}^2t_{s_{i-1}}t_{s_{i}}x_{i-1}t_{s_{i-1}}\\
	&=  t_{s_{i}}( t_{s_{i-1}}x_{i-1}t_{s_{i-1}}t_{s_{i}} t_{s_{i-1}}
		t_{s_{i-1}}t_{s_{i}}t_{s_{i-1}}x_{i-1}t_{s_{i-1}}) ,
	\end{align*}
and
$$
( t_{s_{i-1}}t_{s_{i}}t_{s_{i-1}}
					+ t_{s_{i-1}} )t_{s_{i}} = t_{s_{i}}t_{s_{i-1}}t_{s_{i}}^2 + t_{s_{i}}^2t_{s_{i-1}}t_{s_{i}}
		= t_{s_{i}}( t_{s_{i-1}}+ t_{s_{i-1}}t_{s_{i}}t_{s_{i-1}}).$$
So 
	\begin{align*}
			x_i x_{i+1}
				&=t_{s_{i}} \big((t_{s_{i-1}}x_{i-1}t_{s_{i-1}} + t_{s_{i-1}})t_{s_{i}}(t_{s_{i-1}}x_{i-1}t_{s_{i-1}} + t_{s_{i-1}}) + (t_{s_{i-1}}x_{i-1}t_{s_{i-1}} + t_{s_{i-1}})\big)\\
				 &=t_{s_{i}}x_it_{s_{i}}x_i + t_{s_{i}}x_i
				= x_{i+1}x_i .
		\end{align*}
Finally, assume, without loss of generality, that $i < j$. Then Claim 4 follows by \eqref{eq:x_induct} and Claim 2.


\item[]\textbf{Claim 5:} $w_j x_i= x_i w_j$ for $j>i$.\\
\noindent\emph{Proof:} By \eqref{eq:x_induct}, \\
	\centerline{$\displaystyle x_i = t_{s_{i-1}}\cdots t_{s_{1}}x_{1}t_{s_{1}}\cdots t_{s_{i-1}}
		 + \sum_{\ell = 1}^{i-1} t_{s_{i-1}}\cdots t_{s_{\ell+1}}t_{s_{\ell}}t_{s_{\ell+1}}\cdots t_{s_{i-1}}. $}
So \eqref{eq:comm-tw} and \eqref{eq:comm-xw} imply $w_j x_i= x_i w_j$ for $j>i$.

\item[]\textbf{Claim 6:} $ x_i (w_0  + \cdots + w_i)= (w_0  + \cdots + w_i) x_i $ for $i = 1, \dots k$. 
\\
\noindent\emph{Proof:} This follows by induction on $i$, with $i=1$ satisfied by \eqref{eq:twist-xw0}. Rewrite \eqref{eq:twist-tw} as $t_{s_i}(w_i + w_{i+1}) = (w_i + w_{i+1}) t_{s_i}$, so 
	\begin{align*} 
		x_i (w_0 + \cdots + w_i) &= (t_{s_{i-1}} x_{i-1} t_{s_{i-1}} + t_{s_{i-1}}) (w_0 + \cdots + w_i) \\
				&=  (w_0 + \cdots + w_i) (t_{s_{i-1}} x_{i-1} t_{s_{i-1}} + t_{s_{i-1}})\\
				&= (w_0  + \cdots + w_i) x_i 
	\end{align*}
since $x_{i-1} w_i =  w_ix_{i-1}$ by Claim 5, and $t_{s_\ell} w_j = w_j t_{s_\ell}$ for $\ell < j$ by \eqref{eq:comm-tw}.

\item[]\textbf{Claim 7:} If  $y_1 = w_1 - x_1 +\half(a-p+b-q)$ and  $y_2 = w_2 - x_2 + t_{s_1} +\half(a-p+b-q)$, then $y_1 y_2 = y_2 y_1$ and $y_1 t_{s_i} = t_{s_i} y_1$ for $i > 1$. \\
\noindent\emph{Proof:} Let $K = \half(a-p+b-q)$. So
	\begin{align*}
	y_1 y_2 &= (w_1 - x_1 +K)(w_2 -(t_{s_1}x_1 t_{s_1} + t_{s_1})  + t_{s_1} +K)\\
		&=  (w_2+K)(w_1 - x_1 +K) -(w_1 - x_1 +K)t_{s_1}x_1 t_{s_1} \\
		&=  (w_2+K)(w_1 - x_1 +K) - (t_{s_1}x_1 t_{s_1}) K - w_1t_{s_1}x_1 t_{s_1} + x_1t_{s_1}x_1 t_{s_1} \\
		&=  (w_2+K)(w_1 - x_1 +K) - (t_{s_1}x_1 t_{s_1}) K + x_1t_{s_1}x_1 t_{s_1} 
		 - (t_{s_1}w_2 - 1)x_1 t_{s_1} \\
		&=  (w_2+K)(w_1 - x_1 +K) - (t_{s_1}x_1 t_{s_1}) K 
		  +x_1t_{s_1}x_1 t_{s_1}+x_1 t_{s_1} - t_{s_1}x_1w_2 t_{s_1}\\
		&=  (w_2+K)(w_1 - x_1 +K) - (t_{s_1}x_1 t_{s_1}) K
		  +t_{s_1}x_1 t_{s_1}x_1+ t_{s_1}x_1 - t_{s_1}x_1(t_{s_1}w_1 +1)\\
		&=  (w_2+K)(w_1 - x_1 +K) - (t_{s_1}x_1 t_{s_1}) K 
		  + t_{s_1}x_1 t_{s_1}x_1- t_{s_1}x_1t_{s_1}w_1 \\
		&=(w_2 -t_{s_1}x_1 t_{s_1}  +K) (w_1 - x_1 +K)\\
		&= y_2 y_1.
	\end{align*}
The latter is simply $t_{s_i}$ commuting with $w_1$, $x_1$, and $K$ for $i>1$.

\item[]\textbf{Claim 8:} Let $ m_i = \sum_{j=1}^{i-1} t_{(j~i)} $ and $K = \half(a-p+b-q)$. If $y_1 = w_1 -x_1 + K$, then
$$y_i = w_i - x_i + m_i + K\qquad  \text{ and } \qquad y_i = t_{s_{i-1}}y_{i-1}t_{s_{i-1}} + t_{s_{i-1}}$$
 for $i=2, \ldots, k$ are equivalent definitions of $y_i$. \\
\noindent\emph{Proof:} 
Since  $t_{s_{j}}t_{(i~j)}t_{s_{j}} = t_{(i~j+1)}$, we have $t_{s_j}m_j t_{s_j} = m_{j+1} - t_{s_j}$,
and so
\begin{align*}
	t_{s_i}y_it_{s_i} + t_{s_i} &= t_{s_i}(w_i - x_i + m_i + K)t_{s_i} + t_{s_i}\\
		&= (w_{i+1} - t_{s_i}) - (x_{i+1} - t_{s_i}) + (m_{i+1} - t_{s_i}) + K + t_{s_i}\\
		&= w_{i+1} - x_{i+1} + m_{i+1}  + K = y_{i+1}.
\end{align*}
The other direction follows by induction. 

\item[]\textbf{Claim 9:} If $y_i$ is as in Claim 8, then 
	$$y_i y_j = y_j y_i \textrm{ for $i,j = 1, \dots, k$}, \quad t_{s_i} y_j = y_j t_{s_i} \textrm{ for $j \neq i, i+1$,}$$ 
	$$\textrm{ and } y_i w_0 = w_0 y_i + \big((w_1 + \cdots + w_i)y_i - y_i (w_1 + \cdots + w_i)\big)\text{ for $i = 1, \dots k$. }$$
\emph{Proof:}
These follow from Claims 7 and 8 analogously to the  $x_i$-valued relations above.

\end{enumerate}

\end{proof} 

As a final remark, Theorem \ref{thm:hecke_ext_presentation-short} implies
$\cH_k^\ext \cong \CC[w_0] \otimes \cH_k$ as vector spaces.


\section{Tensor space as a $\cH_k^\ext$-module} 
\label{sec:reps-of-H}

Now we  fix $\fg = \fgl_n$ or $\fsl_n$, and show that for special choices of  $\fg$-modules $M$, $N$, and $V$, the algebra $\cH_k^\ext$ acts on tensor space $M\otimes N \otimes V^{\otimes k}$ by a twist of the representation $\Phi$ in Theorem \ref{thm:braid_group_rep} via an automorphism of $\cG_k$. We go on to explore seminormal representations arising from this representation.

\subsection{Preliminaries on $\fgl_n$ and $\fsl_n$}
\label{sec:prelims}

Let $V = \CC^{n}$ with orthonormal basis $\{ v_1, \dots, v_{n} \}$. We consider the Lie algebras
	$$\fgl_n = \End(V)  \quad \text{and} \quad \fsl_n = \{ x \in \End(V) ~|~ \tr(x) = 0 \}.$$
 Let $\varepsilon_1, \dots, \varepsilon_{n}$ be the orthonormal basis of the weight space $\fh^*$, where if $E_{i,j}\in \End(V)$ is given by $E_{i,j}v_k = \delta_{j,k}v_i$, then $\vep_i(E_{j,j}) = \delta_{i,j}$. The set of positive roots is given by 
$$R^+ = 
\{ \varepsilon_i-\varepsilon_j\ |\ 1\le i<j\le n\}. 
$$
The set roots is $R= R^+ \cup R^-,$ where $R^- = \{ -\alpha ~|~ \alpha \in R^+\}$, and has basis $\{\alpha_i = \varepsilon_i - \varepsilon_{i+1} ~|~  i = 1,\dots, n-1 \}$.

\begin{figure}[!tp]
$$
\begin{matrix}~
\text{(I)}&\quad&\text{(II)}&\quad&\text{(III)}\\
	\begin{picture}(30,60)(2,-12)
		\setlength{\unitlength}{1.1pt}
		\thicklines		
	\color{black}
		\put(0,10) {\line(0,1){30}}
		\put(10,10) {\line(0,1){30}}
		\multiput(20,10)(10,0) {1}  {\line(0,1){30}}
		\multiput(30,30)(10,0) {1}  {\line(0,1){10}}
		\multiput(0,10)(0,10) {2}  {\line(1,0){20}}
		\multiput(0,30)(0,10) {2}  {\line(1,0){30}}
	\end{picture}~
	&&
	\begin{picture}(70,60)(-33.5,-12)
	\setlength{\unitlength}{1.1pt}
	\color{dgrey}
		\put(0,-5) {\line(0,1){50}}
		\put(-2,-13) {\small$0$}
		\thicklines		
	\color{black}
		\put(0,10) {\line(0,1){30}}
		\put(-10,0) {\line(0,1){40}}
		\put(-20,0) {\line(0,1){40}}
		\put(10,10) {\line(0,1){30}}
		\multiput(20,10)(10,0) {1}  {\line(0,1){30}}
		\multiput(30,30)(10,0) {1}  {\line(0,1){10}}
		\multiput(-27,10)(0,10) {2}  {\line(1,0){47}}
		\multiput(-27,0)(0,10) {1}  {\line(1,0){17}}
		\multiput(-27,30)(0,10) {2}  {\line(1,0){57}}
		\multiput(-42,-2.5)(0,10){5} {$\cdots$}
	\end{picture}~
	&&~
		\begin{picture}(70,60)(-32,-12)
	\setlength{\unitlength}{1.1pt}
		\color{dgrey}
		\put(0,-5) {\line(0,1){50}}
		\put(-2,-13) {\small$0$}
		\thicklines		
	\color{black}
		\put(0,10) {\line(0,1){30}}
		\put(-10,0) {\line(0,1){40}}
		\put(-20,0) {\line(0,1){40}}
		\put(10,10) {\line(0,1){30}}
		\multiput(20,10)(10,0) {1}  {\line(0,1){30}}
		\multiput(30,30)(10,0) {1}  {\line(0,1){10}}
		\multiput(-27,10)(0,10) {2}  {\line(1,0){47}}
		\multiput(-27,0)(0,10) {1}  {\line(1,0){17}}
		\multiput(-27,30)(0,10) {2}  {\line(1,0){57}}
		\multiput(-42,-2.5)(0,10){5} {$\cdots$}
		 \put(1,12){\small -2} \put(1,22){\small -1}  \put(3.2,32){\small 0} 
		 \put(-8.5,32){\small -1}
		 \put(-8.5,22){\small -2}\put(-8.5,12){\small -3}
		\put(-18.5,32){\small -2} \put(-18.5,22){\small -3}\put(-18.5,12){\small -4}\put(-18.5,2){\small -5}
		 \put(11,12){\small -1} \put(13.2,22){\small 0}  \put(13.2,32){\small 1}
		  \put(23.2,32){\small 2}
	\end{picture}~
	\\
	\fg=\fsl_n && \fg=\fgl_n &&\\
	\text{Partition assoc. to}&&\text{Multiseg. assoc. to} && \text{Multisegment from (II)}\\
	\lambda = 3 \varepsilon_1 + 2 \varepsilon_2 + 2 \varepsilon_3 &&\lambda = 3 \varepsilon_1 + 2 \varepsilon_2 &&\text{ filled in with contents,}\\
	- \frac{7}{n}(\vep_1 + \cdots + \vep_n)&&+ 2 \varepsilon_3 - \varepsilon_4&&\text{as defined in \eqref{content}.}
	\end{matrix}
$$
\caption{Weights represented as multisegments. }
\label{fig:weights-as-partitions}
\end{figure}
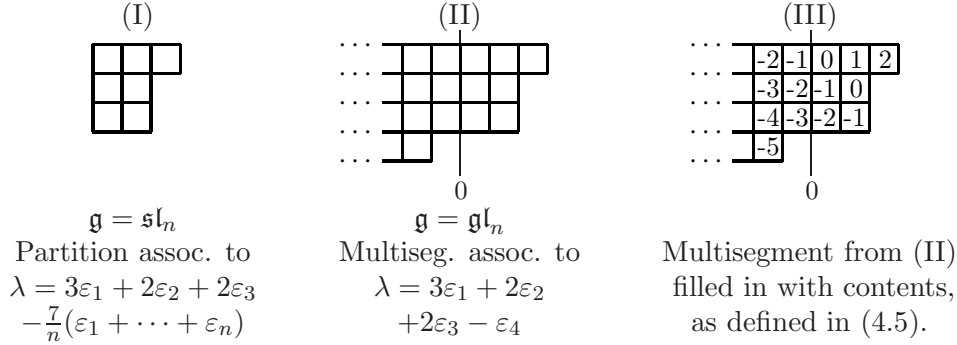

 The finite-dimensional irreducible $\fsl_n$-modules are indexed by the dominant integral weights 
$$P^+ = \left\{
		\lambda = \lambda_1 \varepsilon_1 + \cdots + \lambda_{n-1}\varepsilon_{n-1}
				 - \frac{|\lambda|}{n}(\varepsilon_1 + \cdots + \varepsilon_{n})
~~\Bigg|~
			\begin{array}{c}
				\lambda_i \in \ZZ_{\geq 0},\\
				\lambda_1 \geq \cdots \geq \lambda_{n-1} \geq 0,\\
				|\lambda| = \lambda_1 + \cdots + \lambda_{n-1},
			\end{array} \right\},
		$$
and we identify each weight $\lambda$ with the partition with $\lambda_i$ boxes in row $i$. The finite-dimensional irreducible $\fgl_n$-modules are indexed by the dominant integral weights 
 $$P^+ = \left\{ \lambda = \lambda_1 \varepsilon_1 + \cdots + \lambda_{n} \varepsilon_{n}, 
	~|~  \lambda_i \in \ZZ, ~ \lambda_1 \geq \cdots \geq \lambda_{n}\right\},$$
and we identify each weight $\lambda$ the partition which extends infinitely to the left, and ends on the right in column $\lambda_i$. For examples of each, see Figure \ref{fig:weights-as-partitions} parts (I) and (II). In the case where $\fg=\fgl_n$ and $\lambda_i \geq 0$ for all $1\leq i \leq n$, we often represent $\lambda$ as a finite partition, leaving off boxes to the left of 0. In either case, the first fundamental weight is indexed by a single box, i.e.\ is given by 
	$$\omega_1 = \begin{cases} \vep_1, & \fg=\fgl_n,\\
			\vep_1 - \frac{1}{n}( \vep_1 + \dots + \vep_n),& \fg=\fsl_n.
			\end{cases}$$

Let $L(\lambda)$ be the finite-dimensional irreducible highest weight $\fg$-module of weight $\lambda$, i.e. the irreducible $\fg$-module generated by highest weight vector $v_{\lambda}^+$ of weight $\lambda$ with action
	$$hv_{\lambda}^+ = \lambda(h) v_{\lambda}^+ \quad
		\text{ and } \quad
	   x v_{\lambda}^+ =0,  \qquad \text{ for }  h\in \fh, ~x \in \fn^+.$$
In particular, when $\fg= \fgl_n$ or $\fsl_n$, the \emph{standard representation} is $L(\omega_1)$.

The decomposition numbers for the tensor product of two highest weight modules can be calculated using the Littlewood-Richardson rule (see \cite[I,1]{Mac}). The two special cases we require are as follows. 

\begin{example}[Adding a box]
\label{ex:LR-add-a-box}
For $\fg=\fgl_n$ or $\fsl_{n+1}$ and $\mu \in P^+$, 
	\begin{equation}\label{eq:omega_1_tensor}
	L(\mu) \otimes L(\omega_1) =\bigoplus_{\lambda \in \mu^+ } L(\lambda),  \
	\quad
\text{where } \quad
\mu^+ = 
	\left\{\begin{array}{c}\text{partitions of height $\leq n$}\\ \text{obtained by adding a box to $\mu$ }\end{array} \right\} .
	\end{equation}
	\end{example}

\begin{example}[Rectangles] (See \cite[Lem. 3.3]{St}, \cite[Thm 2.4]{Ok}) 	\label{ex:rectangles} Let $p \geq q$ and $a, b$ be non-negative integers. Denote the rectangular partition with $p$ rows of length $a$ by $(a^p)$. 
Then each $L(\lambda)$ has multiplicity 1 in $L((a^p)) \otimes L((b^q))$ if $\lambda \in  \cP$, and is zero otherwise, where
$\cP = \cP((a^p), (b^q))$ is the set of partitions $\lambda$ with height $\leq p+q$ such that 
	\begin{equation}\label{eq:rectangle_reqs}
	\begin{matrix}
		\lambda_{q+1} = \lambda_{q+2} = \cdots = \lambda_p = a,\\
		\lambda_q \geq \mathrm{max}(a,b),\\
		\lambda_i + \lambda_{p+q - i + 1} = a+b, \quad {i = 1, \dots, q}.
	\end{matrix}
	\end{equation}
In other words, $\cP$ is the set of partitions made by placing $(b^q)$ to the right of $(a^p)$, carving a corner out of $(b^q)$, rotating it $180^{\circ}$ and gluing it to the bottom of $(a^p)$. For example, 
{\def\UNIT{0.7pt}
\begin{align*}\setlength{\unitlength}{\UNIT}
	\begin{picture}(22,12)(0,-2)
		\thicklines		
		\color{dred}
		\multiput(0,-10)(10,0) {3}  {\line(0,1){20}}
		\multiput(0,-10)(0,10) {3}  {\line(1,0){20}} 		
	\end{picture}
	\times
	\begin{picture}(35,12)(-1,-2)
		\thicklines
		\multiput(0,-10)(10,0) {4}  {\line(0,1){20}}
		\multiput(0,-10)(0,10) {3}  {\line(1,0){30}} 		
	\end{picture}
	=	&
		\begin{picture}(40, 20)(0,0)\setlength{\unitlength}{\UNIT}
		\color{grey}
		\multiput(0,0)(10,0) {3}  {\line(0,-1){20}}
		\multiput(0,-10)(0,-10) {2}  {\line(1,0){20}} 	
		\color{grey}
		\multiput(40,20)(10,0) {2}  {\line(0,-1){20}}
		\multiput(30,20)(0,-10) {3}  {\line(1,0){20}} 	
		\color{black}
		\thicklines
		\multiput(0, 0)(10,0) {4}  {\line(0,1){20}}
		\multiput(0, 0)(0,10) {3}  {\line(1,0){30}} 		
		\color{dred}
		\multiput(40,10)(10,0) {2}  {\line(0,-1){10}}
		\multiput(40,10)(0,-10) {2}  {\line(1,0){10}} 	
		\multiput(40,20)(10,0) {2}  {\line(0,-1){10}}
		\multiput(40,20)(0,-10) {2}  {\line(1,0){10}} 	
		\multiput(30,10)(10,0) {2}  {\line(0,-1){10}}
		\multiput(30,10)(0,-10) {2}  {\line(1,0){10}} 	
		\multiput(30,20)(10,0) {2}  {\line(0,-1){10}}
		\multiput(30,20)(0,-10) {2}  {\line(1,0){10}} 	
	\end{picture} 
	~+~
	\begin{picture}(40, 20)(0,0)\setlength{\unitlength}{\UNIT}
		\color{grey}
		\multiput(0,0)(10,0) {3}  {\line(0,-1){20}}
		\multiput(0,-10)(0,-10) {2}  {\line(1,0){20}} 	
		\color{grey}
		\multiput(40,20)(10,0) {2}  {\line(0,-1){20}}
		\multiput(30,20)(0,-10) {3}  {\line(1,0){20}} 	
		\color{black}
		\thicklines
		\multiput(0, 0)(10,0) {4}  {\line(0,1){20}}
		\multiput(0, 0)(0,10) {3}  {\line(1,0){30}} 		
		\color{dred}
		\multiput(0,0)(10,0) {2}  {\line(0,-1){10}}
		\multiput(0,0)(0,-10) {2}  {\line(1,0){10}} 	
		\multiput(40,20)(10,0) {2}  {\line(0,-1){10}}
		\multiput(40,20)(0,-10) {2}  {\line(1,0){10}} 	
		\multiput(30,10)(10,0) {2}  {\line(0,-1){10}}
		\multiput(30,10)(0,-10) {2}  {\line(1,0){10}} 	
		\multiput(30,20)(10,0) {2}  {\line(0,-1){10}}
		\multiput(30,20)(0,-10) {2}  {\line(1,0){10}} 	
	\end{picture}~+~
	\begin{picture}(40, 20)(0,0)\setlength{\unitlength}{\UNIT}
		\color{grey}
		\multiput(0,0)(10,0) {3}  {\line(0,-1){20}}
		\multiput(0,-10)(0,-10) {2}  {\line(1,0){20}} 	
		\color{grey}
		\multiput(40,20)(10,0) {2}  {\line(0,-1){20}}
		\multiput(30,20)(0,-10) {3}  {\line(1,0){20}} 	
		\color{black}
		\thicklines
		\multiput(0, 0)(10,0) {4}  {\line(0,1){20}}
		\multiput(0, 0)(0,10) {3}  {\line(1,0){30}} 		
		\color{dred}
		\multiput(0,0)(10,0) {2}  {\line(0,-1){10}}
		\multiput(0,0)(0,-10) {2}  {\line(1,0){10}} 	
		\multiput(0,-10)(10,0) {2}  {\line(0,-1){10}}
		\multiput(0,-10)(0,-10) {2}  {\line(1,0){10}} 	
		\multiput(30,10)(10,0) {2}  {\line(0,-1){10}}
		\multiput(30,10)(0,-10) {2}  {\line(1,0){10}} 	
		\multiput(30,20)(10,0) {2}  {\line(0,-1){10}}
		\multiput(30,20)(0,-10) {2}  {\line(1,0){10}} 	
	\end{picture}
\\ &\\
	&\quad ~+~
	\begin{picture}(40, 20)(0,0)\setlength{\unitlength}{\UNIT}
		\color{grey}
		\multiput(0,0)(10,0) {3}  {\line(0,-1){20}}
		\multiput(0,-10)(0,-10) {2}  {\line(1,0){20}} 	
		\color{grey}
		\multiput(40,20)(10,0) {2}  {\line(0,-1){20}}
		\multiput(30,20)(0,-10) {3}  {\line(1,0){20}} 	
		\color{black}
		\thicklines
		\multiput(0, 0)(10,0) {4}  {\line(0,1){20}}
		\multiput(0, 0)(0,10) {3}  {\line(1,0){30}} 		
		\color{dred}
		\multiput(0,0)(10,0) {2}  {\line(0,-1){10}}
		\multiput(0,0)(0,-10) {2}  {\line(1,0){10}} 	
		\multiput(40,20)(10,0) {2}  {\line(0,-1){10}}
		\multiput(40,20)(0,-10) {2}  {\line(1,0){10}} 	
		\multiput(10,0)(10,0) {2}  {\line(0,-1){10}}
		\multiput(10,0)(0,-10) {2}  {\line(1,0){10}} 	
		\multiput(30,20)(10,0) {2}  {\line(0,-1){10}}
		\multiput(30,20)(0,-10) {2}  {\line(1,0){10}} 	
	\end{picture}	 ~+~
	\begin{picture}(40, 20)(0,0)\setlength{\unitlength}{\UNIT}
		\color{grey}
		\multiput(0,0)(10,0) {3}  {\line(0,-1){20}}
		\multiput(0,-10)(0,-10) {2}  {\line(1,0){20}} 	
		\color{grey}
		\multiput(40,20)(10,0) {2}  {\line(0,-1){20}}
		\multiput(30,20)(0,-10) {3}  {\line(1,0){20}} 	
		\color{black}
		\thicklines
		\multiput(0, 0)(10,0) {4}  {\line(0,1){20}}
		\multiput(0, 0)(0,10) {3}  {\line(1,0){30}} 		
		\color{dred}
		\multiput(0,0)(10,0) {2}  {\line(0,-1){10}}
		\multiput(0,0)(0,-10) {2}  {\line(1,0){10}} 	
		\multiput(0,-10)(10,0) {2}  {\line(0,-1){10}}
		\multiput(0,-10)(0,-10) {2}  {\line(1,0){10}} 	
		\multiput(10,0)(10,0) {2}  {\line(0,-1){10}}
		\multiput(10,0)(0,-10) {2}  {\line(1,0){10}} 	
		\multiput(30,20)(10,0) {2}  {\line(0,-1){10}}
		\multiput(30,20)(0,-10) {2}  {\line(1,0){10}} 	
	\end{picture}~+~
	\begin{picture}(40, 20)(0,0)\setlength{\unitlength}{\UNIT}
		\color{grey}
		\multiput(0,0)(10,0) {3}  {\line(0,-1){20}}
		\multiput(0,-10)(0,-10) {2}  {\line(1,0){20}} 	
		\color{grey}
		\multiput(40,20)(10,0) {2}  {\line(0,-1){20}}
		\multiput(30,20)(0,-10) {3}  {\line(1,0){20}} 	
		\color{black}
		\thicklines
		\multiput(0, 0)(10,0) {4}  {\line(0,1){20}}
		\multiput(0, 0)(0,10) {3}  {\line(1,0){30}} 		
		\color{dred}
		\multiput(0,0)(10,0) {2}  {\line(0,-1){10}}
		\multiput(0,0)(0,-10) {2}  {\line(1,0){10}} 	
		\multiput(0,-10)(10,0) {2}  {\line(0,-1){10}}
		\multiput(0,-10)(0,-10) {2}  {\line(1,0){10}} 	
		\multiput(10,0)(10,0) {2}  {\line(0,-1){10}}
		\multiput(10,0)(0,-10) {2}  {\line(1,0){10}} 	
		\multiput(10,-10)(10,0) {2}  {\line(0,-1){10}}
		\multiput(10,-10)(0,-10) {2}  {\line(1,0){10}} 	
	\end{picture}
\\&
%
%
\end{align*}}
A useful visualization of these partitions is given in Figure \ref{fig:vis-of-rectangles}. 

\begin{figure}[!tp]
$$
\begin{picture}(165,130)(-45,-60)
		\put(0,-30){\dashbox{2}(25,30)[c]{$\mu'$}} 
		\put(50,10){\dashbox{2}(25,30)[c]{$\mu$}} 
	\thicklines
		\multiput(0,0)(50,0){2} {\line(0,1){40}} 
		\multiput(0,0)(0,40){2} {\line(1,0){50}} 
		\put(22,42){$a$} 
		\put(60,42){$b$} 
		\put(-8,17){$p$} 
		\put(78,22){$q$} 
		\put(-8,-17){$q$} 
		\put(10,-40){$b$} 
	\put(0,65){\makebox(75,10)[c]{$a > b$ :} }
	\put(0,-60){\makebox(75,10)[c]{$\mu$ is a partition in a $b \times q$ box} } 
	\put(0,-72){\makebox(75,10)[c]{$\mu'$ is the $180^\circ$ rotation of $(b^q)/\mu$} } 
\end{picture}~~~~
\begin{picture}(165,110)(-45,-60)
		\put(0,0){\dashbox{1}(25,40)[c]{}} 
		\put(0,10){\dashbox{1}(50,30)[c]{}} 
		\put(0,-30){\dashbox{2}(25,30)[c]{$\mu'$}} 
		\put(50,10){\dashbox{2}(25,30)[c]{$\mu$}} 
	\thicklines
		\put(0,0) {\line(0,1){40}} 
		\put(0,40) {\line(1,0){50}} 
		\put(50,10) {\line(0,1){30}} 
		\put(25,10) {\line(1,0){25}} 
		\put(25,0) {\line(0,1){10}} 
		\put(0,0) {\line(1,0){25}} 
		\put(22,42){$b$} 
		\put(60,42){$a$} 
		\put(-8,17){$p$} 
		\put(78,22){$q$} 
		\put(-8,-17){$q$} 
		\put(10,-40){$a$}
	\put(0,65){\makebox(75,10)[c]{$a < b$ :} } 
	\put(0,-60){\makebox(75,10)[c]{$\mu$ is a partition in an $a \times q$ box} } 
	\put(0,-72){\makebox(75,10)[c]{$\mu'$ is the $180^\circ$ rotation of $(a^q)/\mu$} } 
\end{picture}$$
\caption{An illustration of partitions in $\cP = \cP((a^p),(b^q))$. Outlined sections are filled full with boxes, and dashed regions are filled with complementary partitions.} 
\label{fig:vis-of-rectangles}
\end{figure}
\end{example}

\subsubsection{The Casimir element and the operator $\gamma$}
\label{sec:Casimir}
When $\fg = \fsl_n$, we distinguish the weight
	\begin{equation}\label{rho_def}
	\rho = \frac{1}{2} \sum_{\alpha \in R^+} \alpha = \frac{1}{2}\sum_{i=1}^{n} (n + 1 - 2i)\varepsilon_i.
\end{equation}
When $\fg=\fgl_n$, we choose the analogous weight 
	 \begin{equation}\label{eq:delta}
	 \delta = (n-1)\varepsilon_1 + (n-2)\varepsilon_2 + \cdots +\varepsilon_{n-1} = \sum_{i=1}^{n} (n-i)\varepsilon_i,
	 \end{equation}
which matches \cite[I,1]{Mac}. Keeping $\fg = \fgl_n$ or $\fsl_n$, recall from \eqref{eq:kappadefn} and \eqref{gammadefn} that 
	$$\kappa = \sum_{i} b_{i}b_{i}^*\qquad\hbox{and}\qquad \gamma = \sum_i b_i\otimes b_i^*.$$ 

\begin{lemma}~
\label{lem:Casimir_constant}
The Casimir element $\kappa$ acts on  $L(\lambda)$ by the constant 
		$$\kappa_{L(\lambda)} = \begin{cases}
			\< \lambda, \lambda + 2 \delta \>  - (n-1)|\lambda| &
			\text{when } \fg=\fgl_n, \\
			\< \lambda, \lambda + 2 \rho \> & \text{when } \fg=\fsl_n.
			\end{cases}$$ 
It follows that if $L(\lambda)$ is a submodule of $L(\mu)\otimes L(\nu)$, then $\gamma$ acts on the $L(\lambda)$ isotypic component of $L(\mu)\otimes L(\nu)$ by the constant 
\begin{equation*}\label{tvalue}
\gamma_{\mu \nu}^\lambda = \begin{cases}
	\half\big(\< \lambda,\lambda+2\delta\>
- \< \mu,\mu+2 \delta\> - \< \nu,\nu+2 \delta\>\big)& \text{when } \fg=\fgl_n,\\
	\half(\<\lambda,\lambda+2\rho\>
- \< \mu,\mu+2\rho\> - \< \nu,\nu+2\rho\>)& \text{when } \fg=\fsl_n.
\end{cases}
\end{equation*}
\end{lemma}
\begin{proof}
Both cases are classical results. We include here an argument for $\fg=\fgl_n$, as it is illustrative of both. For the action of $\kappa$ when $\fg=\fsl_n$, see also \cite[\S 8.2]{Jac}. The elementary matrices $\{E_{ij} ~|~ 1 \leq i,j \leq n \}$ form a basis of $\fgl_n$ with dual basis $\{ E_{ji} ~|~  1 \leq i,j \leq n\}$ with respect to the trace form. So
	$$
	 \kappa	=  \sum_{1 \leq i,j \leq n} E_{ij}E_{ji}
= \sum_{i=1}^n E_{ii}E_{ii}
				+ \sum_{1 \leq i<j \leq n} (E_{ii} - E_{jj} + 2E_{ji}E_{ij} ),$$
and therefore
	 \begin{align*}
	 \kappa v_{\lambda}^+ 
	 			&=\left(\sum_{i=1}^n \lambda_i^2
				+ \sum_{1 \leq i<j \leq n}  \lambda_i - \lambda_j + 0
	\right)v_{\lambda}^+\\			
				&= \left(\< \lambda, \lambda \>
				+ \sum_{i=1}^n \big((n-i) - (i-1)\big)\lambda_i
	\right)v_{\lambda}^+\\
				&= \left(\< \lambda, \lambda \>
					+ \sum_{i=1}^n \big( 2n-2i \big)\lambda_i - (n-1)\lambda_i
	\right)v_{\lambda}^+\\	
	&= \left(\< \lambda, \lambda \>
					+ \<\lambda, 2 \delta \> - (n-1)|\lambda|
	\right)v_{\lambda}^+.
	 \end{align*}
Since $\kappa$ acts on $L(\mu)\otimes L(\nu)$ by $(\kappa \otimes \id_{L(\nu)}) + (\id_{L(\mu)} \otimes \kappa) + 2\gamma$, 
$$	\gamma_{\mu \nu}^\lambda = \half\big(\< \lambda,\lambda+2\delta\>
- \< \mu,\mu+2 \delta\> - \< \nu,\nu+2 \delta\>\big) - \hbox{$\frac{n-1}{2}$}\big( | \lambda| - |\mu| - |\nu|\big).$$
But if $L( \lambda) \subseteq L(\mu) \otimes L(\nu)$, then $|\lambda| = |\mu| + |\nu|$, so the desired action of  $\gamma$ follows.
\end{proof}

If $B$ is the box in  column $c$ and row $r$ of a partition $\lambda$, the \emph{content} of $B$ is
	\begin{equation}\label{content}
	c(B) = c-r.
	\end{equation}
See Figure \ref{fig:weights-as-partitions} part (III) for an example of a filling of boxes in a multisegment with their respective contents. We can now give a combinatorial description of the values $\gamma$ takes on tensor products in the special cases described in Examples \ref{ex:LR-add-a-box} and \ref{ex:rectangles}.

\begin{lemma} 
\label{thm:gamma_contents}
If $L(\lambda)$ is a submodule of $L(\mu)\otimes L(\omega_1)$, then $\gamma$ acts on the $L(\lambda)$ isotypic component of $L(\mu)\otimes L(\omega_1)$ by the constant
$$\gamma_{\mu \omega_1}^{\lambda} = \begin{cases}
	 c(\lambda/\mu),& \text{if } \fg=\fgl_n,\displaystyle\\
	 c(\lambda/\mu) - \frac{|\mu|}{n}, & \text{if } \fg=\fsl_n,
	 \end{cases} $$
		where $\lambda/\mu$ is the box added to $\mu$ to obtain $\lambda$.
\end{lemma}
\begin{proof}These values are also known in the literature, but we give an illustrative calculation. Let $\fg = \fgl_n$ and write 
$\mu = \mu_1 \varepsilon_1 + \cdots + \mu_{n} \varepsilon_{n} $.
Adding a box to $\mu$ in the $i^{\text{th}}$ row is equivalent to adding $\varepsilon_i$ to $\mu$. So, since $\omega_1 = \varepsilon_1$,  by Lemma \ref{lem:Casimir_constant},
	\begin{align*}
		2 \gamma_{\mu \omega_1}^{\lambda} &= \big( \<\mu + \varepsilon_i , \mu + \varepsilon_i + 2\delta \> -  \< \mu, \mu + 2 \delta \>  - \< \omega_1, \omega_1 + 2 \delta \> \big)\\
			&= 2\< \mu, \vep_1 \> + 2 \< \varepsilon_i - \varepsilon_1, \mu \> + 2 \< \varepsilon_i - \varepsilon_1, \vep_1 \>
					+ \<\varepsilon_i - \varepsilon_1, \varepsilon_i - \varepsilon_1 + 2 \delta \>\\
			&= 2\big( \mu_1 +\mu_i - \mu_1 - 1 + 1+ (n-i) - (n-1)\big) = 2( \mu_i +1-i).
		\end{align*}
A box added to row $i$ of $\mu$ is in position $(i, \mu_i +1)$ and has content $(\mu_i + 1) - i$, so $\gamma_{\mu \omega_1}^{\lambda} = c(\lambda/\mu).$

The case where $\fg=\fsl_n$ follows analogously, since adding a box to $\mu$ in the $i^{\text{th}}$ row is equivalent to adding $\varepsilon_i - \frac{1}{n}(\varepsilon_1 + \cdots + \varepsilon_{n})$ to $\mu$, and so $\lambda = \mu + \varepsilon_i - \varepsilon_1 + \omega_1.$

\end{proof}


Finally, we will need to understand the the action of $\gamma$ on  $M\otimes N$, where $M$ and $N$ are indexed by rectangular partitions. Specifically, fix $a, b, p, q$ positive integers with $p \geq q$ and 
$$p+q \leq \begin{cases} n, & \text{if } \fg = \fgl_n,\\ n-1, & \text{if } \fg = \fsl_n.\end{cases} $$
Let $M = L((a^p))$, $N= L((b^q))$, and recall from Example \ref{ex:rectangles} that nontrivial submodules of $M\otimes N$ have multiplicity 1 and are indexed by partitions in $\cP$. 

\begin{remark}\label{rk:move-a-box}
As a consequence of the description of $\cP$ in \eqref{eq:rectangle_reqs}, if a box in $\lambda \in \cP((a^b),(p^q))$ is moved from position $(i,j)$ to form another partition in $\cP$, it must be moved to position $(a+b+1-i, p+q+1-j)$. 
\end{remark}

Any partition in  $\cP$ can be built iteratively by beginning with the partition $$(a^p) + (b^q) =  \begin{cases}
a(\vep_1 + \cdots + \vep_p) + b(\vep_{1} + \cdots + \vep_q),  & \text{when }\fg = \fgl_n,\\
a(\vep_1 + \cdots + \vep_p) + b(\vep_{1} + \cdots + \vep_q)  - \frac{ap + bq}{n}(\vep_1 + \cdots + \vep_n), 
 & \text{when }\fg = \fsl_n,\end{cases}
$$
and moving successive boxes down. An example of this process is given in Figure \ref{fig:move-a-box}.
%

%
%
%

\begin{lemma} 
\label{lem:move-a-box}
Let $\mu$  and $\lambda$ index distinct non-trivial components of $M \otimes N$, assume $\lambda$ differs from $\mu$ by moving one box from position $(\mu_i,i)$.  Denote the constant by which $\kappa$ acts on an irreducible component $L(\nu)$ as $\kappa_{L(\nu)}$.
Then
	$$\kappa_{L(\lambda)} = \kappa_{L(\mu)} - 4( (\mu_i - i) - \half (a-p+b-q) ).$$
\end{lemma}
\begin{proof}
If $\fg = \fgl_{n}$ and $\lambda = \mu - \varepsilon_i + \varepsilon_j$ is obtained from $\mu$ by moving a box from row $i$ into row $j$, then
	\begin{align*}
		\kappa_{L(\lambda)}	&=  \< \lambda, \lambda + 2 \delta \> - (n-1)|\lambda|\\
			&= \<\mu, \mu + 2 \delta \>  - (n-1)|\mu| + 2 \< \mu, \vep_j - \vep_i\> + \<\vep_j- \vep_i, \vep_j - \vep_i + 2 \delta \> \\
			&= \kappa_{L(\mu)} - 2\Big(  (\mu_i - i) - (\lambda_j - j)\Big)\\
			&= \kappa_{L(\mu)} - 2 ( \mathrm{ content~of~old~box } -  \mathrm{ content~of~new~box }).
	\end{align*}
A similar computation will show the same for $\fg = \fsl_{n}$.	
Now, if $\lambda$ and $\mu$ are both elements of $\cP$, then 
	$j = p+q +1 -i \text{ and } \lambda_j = a+b + 1 - \mu_i.$
So 
	\begin{align*}\kappa_{L(\lambda)} 
		&= \kappa_{L(\mu)} - 2\Big( (\mu_i - i) - \big((a-p) + (b-q) -(\mu_i-i)\big)\Big)\\
		&= \kappa_{L(\mu)} - 4( (\mu_i-i) - \half (a-p+b-q) ).
	\end{align*}
\end{proof}

%

%

\begin{lemma}
\label{lem:z_0_on_lambda}
Let $\lambda \in \cP$ and define $\cB_{\lambda}$ to be the set of boxes in $\lambda$ in rows $p+1$ and below. Then $\gamma$ acts on an irreducible component $L(\lambda)$ of $L((a^p)) \otimes L((b^q))$ by the constant 
\begin{equation*}
	\gamma_{(a^p)(b^q)}^{\lambda}= 
	\begin{cases} \displaystyle
	abq + 2\sum_{B \in \cB_\lambda} \left(c(B) - \half(a-p + b-q)\right), &
	\text{if } \fg= \fgl_n,\\ \displaystyle
	abq - \frac{abpq}{n} + 2\sum_{B \in \cB_\lambda} \left(c(B) - \half(a-p + b-q)\right),&
	\text{if } \fg=\fsl_{n}.
	\end{cases}
\end{equation*}
\end{lemma}
\begin{proof}
Both cases proceed similarly by a direct calculation. 

Let $\fg = \fgl_{n}$. By Lemma \ref{lem:Casimir_constant}, 
	\begin{align*}
		\kappa_{L((a^p) + (b^q))}
			&= \left\< (a^p) + (b^q), (a^p) + (b^q) + 2 \delta \right\> - (n-1)(ap + bq) \\
			& = \< (a^p), (a^p) +  2 \delta \> - (n-1)ap\\
			&\qquad
				+\<(b^q), (b^q) + 2 \delta \>  - (n-1)bq
				+2 \< (a^p), (b^q) \>\\
			&= \kappa_M + \kappa_N + 2 \< a(\vep_1 + \cdots + \vep_p), b(\vep_1 + \cdots + \vep_q)\>\\
			&=  \kappa_M + \kappa_N + 2abq.
	\end{align*}

Since any partition indexing a component of $M \otimes N$ can be arrived at recursively by beginning with $\lambda_0 = (a^p) + (b^q)$ and moving boxes down, iteratively applying Lemma \ref{lem:move-a-box} implies
	\begin{align*}
	\kappa_{L(\lambda)} &= \kappa_{L(\lambda_0)}  + 4\sum_{B \in \cB_\lambda} \left(c(B)- \half(a-p+b-q)\right)\\
	&= \kappa_M + \kappa_N +2abq + 4\sum_{B \in \cB_\lambda}\left(c(B)- \half(a-p+b-q)\right).
	\end{align*}
So $\gamma$ acts on the $L(\lambda)$ component of $M \otimes N$ by 
$$\gamma_{(a^p), (b^q)}^{\lambda} = abq  + 2\sum_{B \in \cB_\lambda} \left(c(B)- \half(a-p+b-q)\right).$$

In the case where $\fg=\fsl_n$, 
	$$\kappa_{L((a^p) + (b^q))} = \kappa_M + \kappa_N + 2abq - 2 \frac{apbq}{n},$$
and so the desired result
follows analagously.

\end{proof}

{An illustration of Lemmas \ref{lem:move-a-box} and \ref{lem:z_0_on_lambda} is given in Figure \ref{fig:move-a-box}. We now have all of the machinery needed to rework the representation of the braid algebra $\cG_k$ from Section \ref{sec:braid-action} into a representation of $\cH_k^\ext$, which we will do in Section \ref{sec:hecke-action}. However, the recursive process in Figure \ref{fig:move-a-box} suggests something further about those partitions obtained by adding a box to a partition in $\cP$, as we explore in the following two lemmas. These results will prove useful later in Section \ref{sec:seminormal_reps}. Let $\cP_1$ be the set of partitions which are obtained by adding a box to an element of $\cP$. }
\begin{figure}
$$
{\def\UNIT{.75pt}
\xymatrix{
&&
\setlength{\unitlength}{\UNIT}
\begin{picture}(50, 65)(0,-5)
	\color{black}
	\thicklines
		\multiput(0,20)(42,0){2} {\line(0,1){35}}
		\multiput(0,20)(0,35){2} {\line(1,0){42}}
		\multiput(0,20)(10,0){2}{\line(0,-1){20}}
		\multiput(0,0)(0,10){2}{\line(1,0){10}}
		\put(52,55){\line(0,-1){20}}
		\multiput(42,55)(0,-10){3}{\line(1,0){10}}
		\put(-2,-7){\scriptsize-$p$-$1$}
		\put(3,3){$\cdot$}
		\put(54,37){\scriptsize $a$-$1$}
		\put(45,37){$\cdot$}
		\put(3,36){\tiny$2(\!a~\!$-$p~\!$-$1\!)$}
		%
\end{picture} 
\ar[dr]|(.4){-(a-1+p-1)}
&&\\
\setlength{\unitlength}{\UNIT}
\begin{picture}(55, 45)(0,25)
	\color{black}
	\thicklines
		\multiput(0,20)(40,0){2} {\line(0,1){35}}
		\multiput(0,20)(0,35){2} {\line(1,0){40}}
		\multiput(40,55)(0,-10){3}  {\line(1,0){20}}
		\multiput(50,55)(10,0){2}{\line(0,-1){20}} 
		\put(52,38){\scriptsize $a$}
		\put(15,35){\scriptsize $4a$}
		%
		\end{picture} 
		\ar[r]^{-(a+p)~~~}
~~&~~
\setlength{\unitlength}{\UNIT}
\begin{picture}(65, 45)(0,10)
	\color{black}
	\thicklines
		\multiput(0,20)(40,0){2} {\line(0,1){35}}
		\multiput(0,20)(0,35){2} {\line(1,0){40}}
		\multiput(40,55)(0,-10) {2}  {\line(1,0){20}}
		\put(40,35) {\line(1,0){10}} 
		\put(50,55) {\line(0,-1){20}} 
		\put(60,55) {\line(0,-1){10}}
		\multiput(0,10)(10,0) {2}  {\line(0,1){10}}
		\multiput(0,10)(0,10) {1}  {\line(1,0){10}} 
		\put(1,13){\scriptsize-$\!p$}
		\put(41,28){\scriptsize $a$-$1$}
		\put(43,37){$\cdot$}
		\put(61,48){\scriptsize $a\!\!+\!\!1$}
		\put(53,47){$\cdot$}
		\put(5,38){\scriptsize $3a\!-\!p$}
		%
		\end{picture} 
		\ar[dr]|(.3){-(a-1+p-1)}	\ar[ur]|(.7){-(a+1+p+1)}
~~&~&~~
\setlength{\unitlength}{\UNIT}
\begin{picture}(50, 65)(0,-5)
	\color{black}
	\thicklines
		\multiput(0,20)(40,0){2} {\line(0,1){35}}
		\multiput(0,20)(0,35){2} {\line(1,0){40}}
		\multiput(0,20)(10,0){2}{\line(0,-1){20}}
		\put(20,20){\line(0,-1){10}}
		\put(0,10){\line(1,0){20}}
		\put(0,0){\line(1,0){10}}
		\multiput(40,55)(0,-10){2}{\line(1,0){10}}
		\put(50,55){\line(0,-1){10}}
		\put(5,38){\scriptsize $a\!-\!3p$}
		\put(42,48){\scriptsize $a$}
		\put(-2,-7){\scriptsize-$p$-$1$}
		\put(3,2){$\cdot$}
		\put(23,13){\scriptsize-$p\!+\!1$}
		\put(13,12){$\cdot$}
		%
		\end{picture} 
\ar[r]^{-(a+p)}
&~~
\setlength{\unitlength}{\UNIT}
\begin{picture}(40, 65)(0,5)
	\color{black}
	\thicklines
		\multiput(0,20)(40,0){2} {\line(0,1){35}}
		\multiput(0,20)(0,35){2} {\line(1,0){40}}
		\multiput(0,20)(10,0){3}{\line(0,-1){20}}
		\multiput(0,0)(0,10){2}{\line(1,0){20}}
		\put(11,3){\scriptsize-$\!p$}
		\put(9,35){\scriptsize $-4p$}
		%
\end{picture} 
\\
&&
\setlength{\unitlength}{\UNIT}
\begin{picture}(60, 50)(0,0)
	\color{black}
	\thicklines
		\multiput(0,10)(42,0){2} {\line(0,1){35}}
		\multiput(0,10)(0,35){2} {\line(1,0){42}}
		\multiput(0,10)(10,0){3}{\line(0,-1){10}}
		\put(0,0){\line(1,0){20}}
		\multiput(42,45)(0,-10){2}{\line(1,0){20}}
		\multiput(52,45)(10,0){2}{\line(0,-1){10}}
		\put(3,-7){\scriptsize -$p\!+\!1$}
		\put(13,2){$\cdot$}
		\put(47,27){\scriptsize $a\!+\!1$}
		\put(55,37){$\cdot$}
		\put(2,25){\tiny $2(\!a~\!$-$p\!+\!\!1\!)$}
		%
\end{picture} 
\ar[ur]|(.7){-(a+1+p+1)}
&&\\
}
}
$$
\caption{An illustration of Lemmas \ref{lem:move-a-box} and \ref{lem:z_0_on_lambda}: the process of constructing partitions in $\cP$, those partitions indexing nontrivial components of $L((a^p)) \otimes L((b^q))$. In this example, $a,p \geq 2$ and $b=q=2$. The leftmost partition is $(a^p) + (2^2)$. The larger outlined area represents $a \times p$ boxes. Partitions are labeled with the action of $\gamma$ in the case where $\fg = \fgl_n$. Edges represent a box in the leftmost partition being moved down to its lower complementary position (as described in Remark \ref{rk:move-a-box}) to form the rightmost partition, and are labeled by the change this presents in the value of $\gamma$. Boxes are marked if they are a change to the left or right, and are labeled by their contents. }
\label{fig:move-a-box}
\end{figure}

\begin{lemma} 
\label{lem:one_or_two}
If $\mu \in \cP_1((a^p),(b^q))$, then there are exactly one or two $\lambda \in \cP$ for which $\lambda \subseteq \mu$. 
\end{lemma}
\begin{proof}
 As described in  Example \ref{ex:rectangles},  $\cP$ is the set of partitions $\lambda$ with height $\leq p+q$ such that 
	\begin{equation}\label{eq:rectangle_reqs_tensor_space}
	\begin{matrix}
		\lambda_{q+1} = \lambda_{q+2} = \cdots = \lambda_p = a,\\
		\lambda_q \geq \mathrm{max}(a,b),\\
		\lambda_i + \lambda_{p+q - i + 1} = a+b, \quad {i = 1, \dots, q}.
	\end{matrix}
	\end{equation}
Again, a useful visualization of these partitions is provided in Figure \ref{fig:vis-of-rectangles}. As stated in Remark \ref{rk:move-a-box}, if a box is removed from $\lambda \in \cP$ in position $(i,j)$, then a box must be added to position $(a+b+1-i, p+q + 1 - j)$ to get another partition in $\cP$. 
Consider a partition $\mu \in \cP_1((a^p),(b^q))$.  Assume, in addition to having $p \geq q$, that if $p=q$ then $a \geq b$. By moving through the criteria in \eqref{eq:rectangle_reqs_tensor_space} and considering addable boxes for a partition which meets these criteria, we can see that this partition falls into one of the following categories.
	\begin{enumerate}
	\item $\mu$ has height $p+q+1$: In this case, exactly one box can be removed to form a partition which satisfies \eqref{eq:rectangle_reqs_tensor_space}, the box in position $(1, p+q+1)$. This partition $\mu$ looks like the partition in Figure \ref{fig:added-boxes-for-1Ds} with only box 1 added. 
	For example
		$$
	\text{if }
	(a^p) = \begin{picture}(25,10)(0,7)
	\thicklines
		\multiput(0,0)(5,0){6} {\line(0,1){20}} 
		\multiput(0,0)(0,5){5} {\line(1,0){25}} 
	\end{picture}~, \quad 
	(b^q) = \begin{picture}(10,5)(0,2)
	\thicklines
		\multiput(0,0)(5,0){3} {\line(0,1){10}} 
		\multiput(0,0)(0,5){3} {\line(1,0){10}} 
	\end{picture}~, \quad \text{ and } \quad
	\mu=\begin{picture}(25,17.5)(0,0)
	\thicklines
		\multiput(0,0)(5,0){6} {\line(0,1){20}} 
		\multiput(0,0)(0,5){5} {\line(1,0){25}} 
		\put(0,-5) {\line(1,0){10}} 
		\put(0,-10) {\line(1,0){10}} 
		\put(0,-15) {\line(1,0){5}} 
		\put(0,0) {\line(0,-1){15}} 
		\put(5,0) {\line(0,-1){15}} 
		\put(10,0) {\line(0,-1){10}} 
	\end{picture}~,
	\text{ then $\mu$ came from } 
	\lambda = \begin{picture}(25,17.5)(0,0)
	\thicklines
		\multiput(0,0)(5,0){6} {\line(0,1){20}} 
		\multiput(0,0)(0,5){5} {\line(1,0){25}} 
		\put(0,-5) {\line(1,0){10}} 
		\put(0,-10) {\line(1,0){10}} 
		%
		\put(0,0) {\line(0,-1){10}} 
		\put(5,0) {\line(0,-1){10}} 
		\put(10,0) {\line(0,-1){10}} 
	\end{picture}~. \phantom{\Bigg|}$$
	\item $\mu_{q+1} = a+1$: In this case, there is exactly one box which can be removed to obtain a partition which satisfies \eqref{eq:rectangle_reqs_tensor_space}, the box in position $(a+1, q+1)$. This partition $\mu$ looks like the partition in Figure \ref{fig:added-boxes-for-1Ds} with only box 2 added. 

	\item $\mu_1 = a+b+1$: Again, there is exactly one box which can be removed, the box in position $(a+b+1, 1)$. This partition $\mu$ looks like the partition in Figure \ref{fig:added-boxes-for-1Ds} with only box 3 added. 
	
	\item $\mu_{p+1} = b+1$: This is similar to the case above, but is a little more complex. We can only see $\mu_{p+1} = b+1$ when $a>b$ and $\mu_{q} = a$. So the only removable box is the one in position $(b+1, p+1)$. This partition $\mu$ looks like the partition in Figure \ref{fig:added-boxes-for-1Ds} with only box 4 added. 
	
	\item $\mu_j + \mu_{p+q - j + 1} = a+b +1$ for some $1 \leq j \leq p$, but $\mu_j < a+ b +1$ and  $\mu_{p+q - j + 1} < \min(a,b)+1$: This is the case which will yield two partitions. One is the partition in which we remove the box in position $(\mu_j, j)$; the other is the partition in which we remove the box in position $(a+b+1-\mu_j, p+q+1 - j)$.  This partition $\mu$ looks like those in Figure \ref{fig:added-boxes-for-2Ds}, where the boxes marked $x$ and $y$ are corner boxes, one of $x$ or $y$ has position $(i,j)$, and the other has position  $(a+b+1-\mu_j, p+q+1 - j)$. 
	
		For example,
$$
	\text{if }
	(a^p) = \begin{picture}(25,10)(0,7)
	\thicklines
		\multiput(0,0)(5,0){6} {\line(0,1){20}} 
		\multiput(0,0)(0,5){5} {\line(1,0){25}} 
	\end{picture}~, \quad 
	(b^q) = \begin{picture}(10,5)(0,2)
	\thicklines
		\multiput(0,0)(5,0){3} {\line(0,1){10}} 
		\multiput(0,0)(0,5){3} {\line(1,0){10}} 
	\end{picture}~, \quad \text{ and } \quad
	\mu=\begin{picture}(30,20)(0,0)
	\thicklines
		\multiput(0,0)(5,0){6} {\line(0,1){20}} 
		\multiput(0,0)(0,5){5} {\line(1,0){25}} 
		\put(25,20) {\line(1,0){10}} 
		\put(25,15) {\line(1,0){10}} 
		\put(25,10) {\line(1,0){10}} 
		\put(30,20) {\line(0,-1){10}} 
		\put(35,20) {\line(0,-1){10}} 
		\put(0,-5) {\line(1,0){5}} 
		\put(0,0) {\line(0,-1){5}} 
		\put(5,0) {\line(0,-1){5}} 
		\end{picture}
~,$$

	$$\text{ then $\mu$ came from } 
	\lambda = \begin{picture}(30,15)(0,7)
		\thicklines
		\multiput(0,0)(5,0){6} {\line(0,1){20}} 
		\multiput(0,0)(0,5){5} {\line(1,0){25}} 
		\put(25,20) {\line(1,0){10}} 
		\put(25,15) {\line(1,0){10}} 
		\put(25,10) {\line(1,0){10}} 
		\put(30,20) {\line(0,-1){10}} 
		\put(35,20) {\line(0,-1){10}} 
		%
		\end{picture} \quad \text{ or } \quad 
		\lambda = \begin{picture}(30,15)(0,5)
		\thicklines
		\multiput(0,0)(5,0){6} {\line(0,1){20}} 
		\multiput(0,0)(0,5){5} {\line(1,0){25}} 
		\put(25,20) {\line(1,0){10}} 
		\put(25,15) {\line(1,0){10}} 
		\put(25,10) {\line(1,0){5}} 
		\put(30,20) {\line(0,-1){10}} 
		\put(35,20) {\line(0,-1){5}} 
		\put(0,-5) {\line(1,0){5}} 
		\put(0,0) {\line(0,-1){5}} 
		\put(5,0) {\line(0,-1){5}} 
		\end{picture}~.$$
	\end{enumerate}
\end{proof}

\begin{figure}[!t]
{\def\UNIT{1.3pt}
$$\setlength{\unitlength}{\UNIT}
\begin{picture}(100,100)(-10,-40)
	\put(0,-37){\framebox(7,7)[c]{1}}
	\put(50,3){\framebox(7,7)[c]{2}}
	\put(75,33){\framebox(7,7)[c]{3}}
	\put(25,-7){\framebox(7,7)[c]{4}}
		\put(0,-30){\dashbox{2}(25,30)[c]{$\nu'$}} 
		\put(50,10){\dashbox{2}(25,30)[c]{$\nu$}} 
	\thicklines
		\multiput(0,0)(50,0){2} {\line(0,1){40}} 
		\multiput(0,0)(0,40){2} {\line(1,0){50}} 
	\put(0,55){\makebox(75,10)[c]{$a > b$ :} }
\end{picture}\quad \text{ or } \quad 
\begin{picture}(100,100)(-10,-40)
	\put(0,-37){\framebox(7,7)[c]{1}}
	\put(25,3){\framebox(7,7)[c]{2}}
	\put(75,33){\framebox(7,7)[c]{3}}
		\put(0,-30){\dashbox{2}(25,30)[c]{$\nu'$}} 
		\put(50,10){\dashbox{2}(25,30)[c]{$\nu$}} 
	\thicklines	
	\put(50,-7){ \color{dgrey}\dashbox(7,7)[c]{4}}
		\put(0,0) {\line(0,1){40}} 
		\put(0,40) {\line(1,0){50}} 
		\put(50,10) {\line(0,1){30}} 
		\put(25,10) {\line(1,0){25}} 
		\put(25,0) {\line(0,1){10}} 
		\put(0,0) {\line(1,0){25}} 
	\put(0,55){\makebox(75,10)[c]{$a < b$ :} } 
\end{picture}$$}
\caption{Added boxes corresponding to partitions with one parent, as described in cases 1-4 in the proof of  Lemma \ref{lem:one_or_two} (see also Figure \ref{fig:vis-of-rectangles}).}
\label{fig:added-boxes-for-1Ds}
\end{figure}
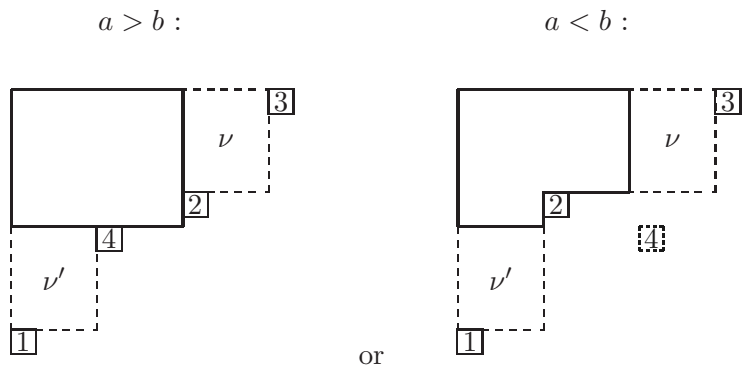
\begin{figure}[!t]
{\def\UNIT{1.3pt}
$$\setlength{\unitlength}{\UNIT}
\begin{picture}(100,115)(-10,-40)
\put(70,45){\color{grey} \line(1,-1){17}}
	\put(87,23){\footnotesize $a+b$}
\put(45,15){\color{grey} \line(1,-1){27}}
	\put(65,-18){\footnotesize$a-q$}
\put(20,5){\color{grey} \line(1,-1){27}}
	\put(40,-28){\footnotesize$b-p$}
\put(-5,-25){\color{grey} \line(1,-1){17}}
	\put(5,-48){\footnotesize$-p-q$}
	\put(63,30){\framebox(7,7)[c]{$x$}}
	\put(5,-27){\framebox(7,7)[c]{$y$}}
		\put(0,-30){\dashbox{2}(25,30)[c]{$~$}} 
		\put(50,10){\dashbox{2}(25,30)[c]{$~$}} 
	\thicklines
		\multiput(0,0)(50,0){2} {\line(0,1){40}} 
		\multiput(0,0)(0,40){2} {\line(1,0){50}} 
	\put(0,55){\makebox(75,10)[c]{$a > b$ :} }
\end{picture}\quad \text{ or } \quad 
\begin{picture}(100,100)(-10,-40)
\put(70,45){\color{grey} \line(1,-1){17}}
	\put(87,23){\footnotesize $a+b$}
\put(45,15){\color{grey} \line(1,-1){27}}
	\put(66,-18){\footnotesize$a-q$}
\put(20,5){\color{grey} \line(1,-1){27}}
	\put(39,-28){\footnotesize$b-p$}
\put(-5,-25){\color{grey} \line(1,-1){17}}
	\put(5,-48){\footnotesize$-p-q$}
	\put(63,30){\framebox(7,7)[c]{$x$}}
	\put(5,-27){\framebox(7,7)[c]{$y$}}
%
		\put(0,-30){\dashbox{2}(25,30)[c]{$~$}} 
		\put(50,10){\dashbox{2}(25,30)[c]{$~$}} 
	\thicklines
		\put(0,0) {\line(0,1){40}} 
		\put(0,40) {\line(1,0){50}} 
		\put(50,10) {\line(0,1){30}} 
		\put(25,10) {\line(1,0){25}} 
		\put(25,0) {\line(0,1){10}} 
		\put(0,0) {\line(1,0){25}} 
	\put(0,55){\makebox(75,10)[c]{$a < b$ :} } 
\end{picture}$$}
\caption{Added boxes corresponding to partitions with two parents, as described in case 5 in the proof of  Lemma \ref{lem:one_or_two} (see also Figure \ref{fig:vis-of-rectangles}). Additionally, critical diagonals are marked with contents $a+b$, $a - q$, $b-p$, and $-p-q$ for Lemma \ref{lem:content-distinction}.}
\label{fig:added-boxes-for-2Ds}
\end{figure}
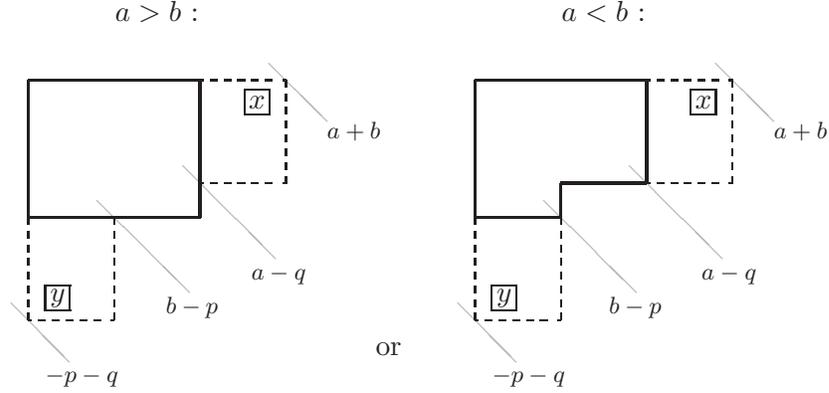

\begin{lemma}\label{lem:content-distinction}
Suppose $\mu \in \cP_1$ and $\lambda \in \cP$ differ by a box, and $c(\mu/\lambda)$ is the content of this box. Then 
\begin{enumerate}
\item there is exactly one such $\lambda$ if and only if $c(\mu/\lambda) = -p-q, ~a-q,~ a+b,~\text{ or } ~ b-p,$ and
\item if $c(\mu/\lambda) \neq -p-q, ~a-q,~ a+b,~\text{ or } ~ b-p,$ then there is exactly one 
	$\lambda' \in \cP$ distinct from $\lambda$ which differs from $\mu$ by a box, and 
		$$c(\mu/\lambda') = a-p+b-q - c(\mu/\lambda) .$$
\end{enumerate}
\end{lemma}
\begin{proof}
If $\mu \in \cP_1$ satisfies cases 1-4 in Lemma \ref{lem:one_or_two},  then 
$$	c(\mu/\lambda) =-p-q, ~a-q,~ a+b,~\text{ or } ~ b-p.$$
The final case yielded two partitions which differ by the movement of one box. If a box in position $(i,j)$ in $\lambda \in \cP$ can be moved to get another partition in $ \cP$, then that box must satisfy either 
\begin{enumerate}
\item[(1)] $\max(a,b) < i \leq a+b \text{ and } 0< j \leq q$, or
\item[(2)] $0< i \leq \min(a,b) \text{ and } p< j \leq p+q$.
\end{enumerate}
If $(i,j)$ satisfies (1), then 	
	$$ \max(a,b) - q < i-j < a+b - q.$$
So since $p \geq q$, 
	$$-p-q < i-j, \quad a-q < i-j, \quad  b-p <i-j , \quad \text{ and } \quad i-j < a+b .$$
If $(i,j)$ satisfies (2), then 
	$$-p-q < i-j < \min(a,b)-p.$$
So, similarly, 
	$$-p-q< i-j, \quad i-j <  a-q, \quad i-j<  b-p \quad \text{ and } \quad i-j< a+b .$$
Thus, if there are two partitions in $ \cP$ which can be obtained by removing a box from $\mu$, then the contents of those boxes are distinct from $-p-q, ~a-q,~ a+b,~\text{ and } ~ b-p.$ See Figure \ref{fig:added-boxes-for-2Ds} for an illustration of these bounds.
\end{proof}

\subsection{Action on tensor space}
\label{sec:hecke-action}
We return now to the representation $\Phi\colon  \cG_k\to \End_{\fg}(M\otimes N\otimes V^{\otimes k})$ in Theorem \ref{thm:braid_group_rep}. Unfortunately, this map does not factor through the quotient defined in \eqref{rel:hecke} and \eqref{rel:hecke2}. However, we can twist by the following automorphism of $\cG_k$ to acquire the desired representations of $\cH_k^\ext$. 

\begin{lemma} \label{lem:G-auto}Fix $c_x$, $c_y$, $c_z$, $d \in \CC$.
	The map $\phi:  \cG_k  \to  \cG_k$ given by 
$$t_{s_i} \mapsto t_{s_i},$$ 
$$x_i  \mapsto 	x_i   + \left((i-1)d + c_x \right), \quad 
y_i  \mapsto 	y_i   + \left((i-1)d + c_y \right), $$
$$z_i  \mapsto 	z_i   + \left((i-1)d + c_x + c_y \right), \quad \text{ and } \quad 
z_0 \mapsto 	z_0   +  c_z,
$$
is an algebra automorphism. 
\end{lemma}
\begin{proof}
Observe that for $i \geq 1$,  
\begin{align*}
\phi(x_{i+1} - t_{s_i} x_i t_{s_i} ) & = x_{i+1} - t_{s_i} x_i t_{s_i} + d,\\
\phi(y_{i+1} - t_{s_i} y_i t_{s_i} ) & = y_{i+1} - t_{s_i} y_i t_{s_i} + d,\quad \text{ and }\\
\phi(m_i) & = m_i + (i-1)d.
\end{align*}
Relations \eqref{rel:graded_braid5}, \eqref{rel:graded_braid6}, and \eqref{rel:zdefn} follow directly.
\end{proof}
Reviewing notation from Section \ref{sec:braid-action}, we denote by $\gamma_{M,N}$ the operator which acts on $M$ and $N$ in $M \otimes N \otimes V^{\otimes k}$ by $\gamma$ and on all other factors by the identity. Similarly, $\gamma_{i,j}$ acts by $\gamma$ on $i^\text{th}$ and $j^\text{th}$ factors of $V$, and for a factor $X$, $\gamma_{X,i}$ acts on $X$ and the $i^\text{th}$ copy of $V$, $\kappa_X$ acts by $\kappa$ on $X$, and $\kappa_{X,\leq i}$ acts by $\kappa$ on $X$ and the first $i$ factors of $V$ (where $\kappa_{X,\leq 0} = \kappa_X$).

Now, define 
\begin{equation} \label{eq:Phi'}
	\Phi' = \Phi \circ \phi: \cG_k \to \End_\fg(M \otimes N \otimes V^{\otimes k}),
	\end{equation}
so that 
$$\Phi'(x_i) = \frac{1}{2}(\kappa_{M,\leq i} -\kappa_{M,\leq i-1})+ \left((i-1)d + c_x \right)\id, \qquad 
\Phi'(y_i) = \frac{1}{2}( \kappa_{N,\leq i} - \kappa_{N,\leq i-1})+ \left((i-1)d + c_y \right)\id,$$
$$\Phi'(z_i) = \frac{1}{2}(\kappa_{M\otimes N,\leq i} -\kappa_{M\otimes N,\leq i-1} + \kappa_V)+ \left((i-1)d + c_x + c_y \right)\id, \qquad \text{for $1 \leq i \leq k$,}$$
$$\Phi'(z_0) = \frac{1}{2}(\kappa_{M \otimes N} - \kappa_M - \kappa_N) +c_z\id= \gamma_{M,N}+c_z\id,\qquad \text{and }$$
$$
\Phi'(t_{s_j}) = \id_M\otimes \id_N\otimes \id_V^{\otimes(j-1)}\otimes s \otimes \id_V^{\otimes (k-j-1)}, \quad \text{where $s \cdot (u\otimes v) =  v \otimes u$}$$
for $1 \leq j \leq k-1$. 


\begin{thm}\label{thm:Hecke-action} Fix $\fg = \fgl_n$ or $\fsl_n$, and let $M = L\left(\left(a^p\right)\right)$, $N= ((b^q))$, and $V = L(\omega_1)$. 
\begin{enumerate}
	\item[(a)] When $\fg = \fgl_n$, fix 
		$c_x = c_y = -\half n,$ and $d=0.$
	\item[(b)] When  $\fg = \fsl_n$, fix 
	$$c_x= \frac{ap}{n}- \frac12\left(n-\frac{1}{n}\right) \quad 
	c_y =  \frac{bq}{n} - \frac{1}{2}\left(n-\frac{1}{n}\right),\quad \text{ and }\quad
	d= \frac{1}{n}.$$
	\end{enumerate}
For this choice of $c_x$, $c_y$, and $d$, and any choice of $c_z$, the map $\Phi'$ in \eqref{eq:Phi'} factors through the quotient by \eqref{rel:hecke} and \eqref{rel:hecke2}, so defines an action of $\cH^\mathrm{ext}_k$ which commutes with the action of $\fg$.

\end{thm}
\begin{proof}~

The relations in \eqref{rel:hecke} can be rewritten as
	\begin{equation*}
		 x_{i+1} - t_{s_i} x_i t_{s_i}  = t_{s_i}, \quad y_{i+1} - t_{s_i} y_i t_{s_i}   = t_{s_i}, \quad i=1, \dots, k-1.
	\end{equation*}
Recall from \eqref{eq:iterative-kappa-expansion} that for $X=M, N,$ or $M\otimes N$,
$$\kappa_{X,\leq j} = \kappa_X  + j\kappa_V +  2\left(\sum_{1 \leq i \leq j} \gamma_{X, i} + \sum_{1\leq r<s \leq j} \gamma_{r,s}\right)$$
and so
	\begin{equation}\label{eq:operator_x}\kappa_{X,\leq i} - \kappa_{X,\leq i-1} = \kappa_V+ 2\gamma_{X,i} + 2\sum_{1\le \ell<i} \gamma_{\ell, i}\end{equation}
as an operator on $X \otimes V^{\otimes k}$. Therefore
	\begin{align*}
	(\kappa_{X, \leq i+1} -& \kappa_{X, \leq i})- s_i  (\kappa_{X, \leq i} - \kappa_{X, \leq i-1}) s_i\\
	&=\kappa_V+ 2\gamma_{X,i+1} + 2\!\!\sum_{1\le \ell<i+1} \gamma_{\ell, i+1}
		- s_i\left( \kappa_V+ 2\gamma_{X,i} + 2\sum_{1\le \ell<i} \gamma_{\ell, i} \right)s_i\\
	&=\kappa_V+ 2\gamma_{X,i+1} + 2\!\!\sum_{1\le \ell<i+1} \gamma_{\ell, i+1}
		- \left( \kappa_V+ 2\gamma_{X,i+1} + 2\sum_{1\le \ell<i} \gamma_{\ell, i+1} \right)\\
	 &=2\gamma_{i,i+1}.
	\end{align*}
	This means that to show \eqref{rel:hecke}, it only remains to be checked that
	\begin{align*}
	\id_M \otimes \id_N  \otimes &\,\id_V^{\otimes i-1} \otimes s \otimes \id_V^{\otimes k-i-1}\\
	=\Phi'(t_{s_i}) &= \Phi'( x_{i+1} - t_{s_i} x_i t_{s_i} ) \\
		&= \half \left( (\kappa_{M, \leq i+1} - \kappa_{M, \leq i}) +2 \left( i~\!d + c_x\right) - s_i  (\kappa_{M, \leq i} - \kappa_{M, \leq i-1}  + 2((i-1)d + c_x)) s_i\right)\\
		&= \gamma_{i,i+1} + d\\
		& = \begin{cases}
		\gamma_{i,i+1} & \text{when $\fg=\fgl_n$,}\\
		\gamma_{i,i+1}  + \frac{1}{n} & \text{when $\fg=\fsl_n$,}
		\end{cases}
	\end{align*}
as operators on $M\otimes N \otimes V^{\otimes k}$ (the check for $\Phi'(t_{s_i}) = \Phi'( y_{i+1} - t_{s_i} y_i t_{s_i} )$ is the same). 

The decomposition of $V\otimes V$ is 
	$$V\otimes V = \VVsym \oplus \VVantisym,$$
	 where if $v_1, \dots, v_n$ is a basis for $V$, then 
	 \begin{equation*}
	 \begin{array}{rl}
	 \VVsym \!\!\!\!&= \mathrm{span}_\CC\{v_i \otimes v_j + v_j \otimes v_i ~|~ 1 \leq i,j \leq n \},
	 \quad \text{and} \\ &\\
	 \VVantisym \!\!\!\!&= \mathrm{span}_\CC\{v_i \otimes v_j - v_j \otimes v_i ~|~ 1 \leq i,j \leq n \}.\end{array}\end{equation*}
	 It follows from this decomposition and Lemma \ref{thm:gamma_contents} that the actions of $s$ and $\gamma$ are given by
$$\begin{array}{l|cc}
\multicolumn{3}{c}{\fg=\fgl_n}\\
& \VVsym & \VVantisym\\
\hline
s&1&-1\\
\gamma&1 &-1 \\
\end{array}
\qquad 
\begin{array}{l|cc}
\multicolumn{3}{c}{\fg=\fsl_n}\\
& \VVsym & \VVantisym\\
\hline
s&1&-1\\
\gamma&1 - \frac{1}{n}&-1 - \frac{1}{n}\\
\end{array}
$$
so \eqref{rel:hecke} is satisfied. 

Next we check $(x_1 - a)(x_1 + p) = 0$. By \eqref{eq:operator_x}, we have 
	$$\Phi'(x_1) = \half\kappa_V+ \gamma_{M,1} + c_x.$$
	The module $M \otimes V$ decomposes as 
	\begin{equation}\label{MVdecomp}M \otimes V  = L\left(
	\begin{picture}(30,15)(0,7)
	\thicklines
		\multiput(0,0)(25,0){2} {\line(0,1){20}} 
		\multiput(0,0)(0,20){2} {\line(1,0){25}} 
		\put(25,20) {\line(1,0){5}} 
		\put(25,15) {\line(1,0){5}} 
		\put(30,20) {\line(0,-1){5}} 
		\put(2,9) {\scriptsize $p$}
		\put(10,22) {\scriptsize $a$}
		\end{picture}
 \right) \oplus L\left( 
 	\begin{picture}(25,20)(0,5)
	\thicklines
		\multiput(0,0)(25,0){2} {\line(0,1){20}} 
		\multiput(0,0)(0,20){2} {\line(1,0){25}} 
		\put(0,-5) {\line(1,0){5}} 
		\put(0,0) {\line(0,-1){5}} 
		\put(5,0) {\line(0,-1){5}} 
		\put(2,9) {\scriptsize $p$}
		\put(10,22) {\scriptsize $a$}
		\end{picture}\right).
	\end{equation}
\begin{enumerate}
\item[]\textbf{Case 1: $\fg = \fgl_n$}\\
By Lemma \ref{lem:Casimir_constant},
	\begin{equation}\label{kappaV-gl}\kappa_V = \< \omega_1, \omega_1 + 2 \delta\> - (n-1)|\omega_1| = 1 + (n-1) - (n-1) = n,\end{equation}
so $\half\kappa_V+  c_x = 0$. By Lemma \ref{thm:gamma_contents} and the decomposition in \eqref{MVdecomp}, $\gamma_{M,1} = a$ or $-p$, so  $\Phi'{(x_1 - a)(x_1 + p)} = 0$ as desired.
\item[]\textbf{Case 2: $\fg = \fsl_n$}\\
By Lemma \ref{lem:Casimir_constant},	
\begin{equation}\label{kappaV-sl}\kappa_V = \< \omega_1, \omega_1 + 2 \rho\>  =n - \frac{1}{n}, \end{equation}
so $\half\kappa_V+  c_x = \frac{ap}{n}$. By Lemma \ref{thm:gamma_contents} and the decomposition in \eqref{MVdecomp}, $\gamma_{M,1} = \left(a - \frac{ap}{n}\right)$ or  $\left(-p - \frac{ap}{n}\right)$ so  $\Phi'{(x_1 - a)(x_1 + p)} = 0$ as desired.
\end{enumerate}
The relation $(y_1 - b)(y_1 + q) =0$ follows analogously, and therefore \eqref{rel:hecke2} is satisfied.

\end{proof}


\subsection{Bratteli diagrams and seminormal bases}
\label{sec:seminormal_bases}
Let 
\begin{equation}\label{chain}
\CC = \cC_0 \subseteq \cC_1 \subseteq  \cC_2 \subseteq \cdots 
\end{equation}
be a chain of semisimple algebras. Let $\hat \cC_k$ be the set of equivalence classes of finite-dimensional irreducible $\cC_k$-modules for $k = 0, 1, \dots$, and write $\cC^\mu$ for a module in the class $\mu \in \hat \cC_k$. Here we will be describing an oriented ranked graph, the \emph{Bratteli diagram} for \eqref{chain}, which encodes the representation theory of $\cC_k$ in terms of the representation theory of $\cC_i$ for $i < k$. In the example where $\cC_k$ is the group algebra of the symmetric group $\CC S_k$, this diagram the same as Young's diagram. This exposition on Bratteli diagrams and seminormal bases closely follows \cite{OV}, where they discuss chains of finite-dimensional semisimple associative algebras in general.

The \emph{Bratteli diagram} associated to a chain \eqref{chain} is an oriented ranked graph, with a rank for each $\cC_i$. The vertices \emph{of rank} or \emph{on level} $k$ are the elements of the set $\hat \cC_k$. Two vertices $\mu \in \hat\cC_{k-1}$ and $\nu \in \hat\cC_{k}$ are joined by $d$ oriented edges from $\mu$ to $\nu$ if 

\begin{equation}\label{multiplicity}d = \dim \Hom_{\cC_{k-1}} (\cC^\mu, \cC^\nu),\end{equation}
i.e.\ $d$ is the multiplicity of $\cC^\mu$ in the restriction of $\cC^\nu$ to a $\cC_{k-1}$-module. Write  
$$\mu \searrow \nu \quad \text{ if $\mu$ and $\nu$ are connected by an edge from $\mu$ to $\nu$.}$$ 
If $\mu \in \hat \cC_i$ and $\lambda \in \hat \cC_k$ with $i<k$, write 
$$\mu \subset \lambda \quad  \text{if there is a path $\mu \searrow \cdots \searrow \lambda$ from $\mu$ to $\lambda$ in the Bratteli diagram}.$$ 
In other words, $\mu \subset \lambda$ if and only if the multiplicity of $\mu$ in $\lambda$ after appropriate restriction is nonzero. 

Our favorite examples are when $\fg = \fsl_n$ or $\fgl_n$, $M$ and $N$ are finite-dimensional simple modules indexed by rectangular partitions, $V$ is the standard representation, and 
\begin{enumerate}[(1)]
\item $\cC_k = \End_\fg(M \otimes V^{\otimes k})$,  
\item $\cC_k = \End_\fg(N \otimes V^{\otimes k})$, or
\item $\cC_0 = \End_\fg(M)$ and $\cC_k = \End_\fg(M \otimes N \otimes V^{\otimes {k-1}}). $
\end{enumerate}
In fact, as we will see in Examples \ref{ex:BratMV} and \ref{ex:BratMNV}, these Bratteli diagrams are all multiplicity free ($d$ in \eqref{multiplicity} is always 0 or 1). With multiplicity free diagrams, the decomposition 
$$\cC^\lambda = \bigoplus_{\mu \in \hat \cC_{k-1} \atop \mu \searrow \lambda} \cC^\mu$$
is canonical. By induction, we obtain a canonical decomposition of the module $\cC^\lambda$ into irreducible one-dimensional $\cC_0$-modules
$$\cC^\lambda = \bigoplus_T \CC v_T$$
indexed by all possible paths
\begin{equation}\label{eq:paths1}
T = (T^{(0)} \searrow T^{(1)} \searrow \dots \searrow T^{(k)} = \lambda),
\end{equation}
where $T^{(i)} \in \hat \cC_i$ for each $0 \leq i \leq k$. 
In particular, $v_T$ is the unique element (up to scalar multiplication) respecting the inductive process, i.e. for each $0 \leq i \leq k$,  after the induction 
$$\Ind_{\cC_0}^{\cC_i} \bigoplus_T \CC v_T = \bigoplus_\nu c_\nu \cC^\nu,$$
each vector $v_T$ lands in the isotypic component $c_{T^{(i)}} \cC^{T^{(i)}}$. 
We call the basis $\{v_T\}$ of $\cC^\lambda$ a (non-normalized) \emph{seminormal basis}. 

Any finite dimensional $\fg$-module $U$ decomposes as a ($\fg$,$\End_\fg(U)$)-bimodule as
\begin{equation}\label{eq:centralizer_thm}
	U \cong  \bigoplus_\lambda L(\lambda) \otimes \cL^\lambda
\end{equation}
where $\cL^\lambda$ are distinct irreducible $\End_\fg(U)$-modules (see \cite[Thm 3.3.7]{GW}), isomorphic to the span of all highest weight vectors of weight $\lambda$ in $U$. So both the irreducible $\fg$-modules and the irreducible $\End_\fg(M \otimes N \otimes V^{\otimes k})$-modules appearing in $M \otimes N \otimes V^{\otimes k}$ are indexed by the same set. 
Therefore, the result of the combinatorics outlined in Section \ref{sec:prelims} is that, for our favorite examples, the paths in \eqref{eq:paths1} are in bijection with specific sets of \emph{tableaux}, which we define now. 

For two partitions $\lambda \subseteq \mu$, the \emph{skew shape} $\mu/\lambda$ is the portion of $\mu$ not contained in $\lambda$. A \emph{(standard) $\mu/\lambda$-tableau} is a filling of the skew shape $\mu/\lambda$ with the integers $1, \dots, |\mu| - |\lambda|$ so that the row fillings increase from left to right and the column fillings increase from top to bottom. For example, if $\lambda = \begin{picture}(25,15)(0,5)
		{\thicklines
			\multiput(0,0)(15,0){2} {\line(0,1){15}} 
		\multiput(0,0)(0,15){2} {\line(1,0){15}} 
		\put(15,15) {\line(1,0){10}} 
		\put(15,10) {\line(1,0){10}} 
		\put(15,5) {\line(1,0){10}} 
		\put(20,15) {\line(0,-1){10}} 
		\put(25,15) {\line(0,-1){10}} 
		}
		\end{picture}
$
and 
$
		\mu = \begin{picture}(35,15)(0,5)
		\thicklines
		\multiput(0,0)(15,0){2} {\line(0,1){15}} 
		\multiput(0,0)(0,15){2} {\line(1,0){15}} 
		\put(15,15) {\line(1,0){10}} 
		\put(15,10) {\line(1,0){10}} 
		\put(15,5) {\line(1,0){10}} 
		\put(20,15) {\line(0,-1){10}} 
		\put(25,15) {\line(0,-1){10}} 
		%
		\put(25,15) {\line(1,0){10}} 
		\put(25,10) {\line(1,0){10}} 
		%
		\put(30,15) {\line(0,-1){5}} 
		\put(35,15) {\line(0,-1){5}} 
		%
		\put(0,-5) {\line(1,0){5}} 
		\put(0,0) {\line(0,-1){5}} 
		\put(5,0) {\line(0,-1){5}} 
		\end{picture}$,
then there are three $\mu/\lambda$-tableaux,
{\def\UNIT{2pt}
\begin{equation}\label{eq:tab_example}
\setlength{\unitlength}{\UNIT}
	\phantom{\Bigg|}T_1 = \begin{picture}(35,15)(0,2.5)
		{\color{grey}
			\multiput(0,0)(15,0){2} {\line(0,1){15}} 
		\multiput(0,0)(0,15){2} {\line(1,0){15}} 
		\put(15,15) {\line(1,0){10}} 
		\put(15,10) {\line(1,0){10}} 
		\put(15,5) {\line(1,0){10}} 
		\put(20,15) {\line(0,-1){10}} 
		\put(25,15) {\line(0,-1){10}} 
		}
		\thicklines
		\put(25,15) {\line(1,0){10}} 
		\put(25,10) {\line(1,0){10}} 
		%
		\put(25,15) {\line(0,-1){5}} 
		\put(30,15) {\line(0,-1){5}} 
		\put(35,15) {\line(0,-1){5}} 
		\put(0,0) {\line(1,0){5}} 
		\put(0,-5) {\line(1,0){5}} 
		\put(0,0) {\line(0,-1){5}} 
		\put(5,0) {\line(0,-1){5}} 
		{\color{dblue}
		\put(1,-4){3}
		\put(26,11){1}
		\put(31,11){2}}
		\end{picture}, \quad 
		T_2 = \begin{picture}(35,15)(0,2.5)
		{\color{grey}
			\multiput(0,0)(15,0){2} {\line(0,1){15}} 
		\multiput(0,0)(0,15){2} {\line(1,0){15}} 
		\put(15,15) {\line(1,0){10}} 
		\put(15,10) {\line(1,0){10}} 
		\put(15,5) {\line(1,0){10}} 
		\put(20,15) {\line(0,-1){10}} 
		\put(25,15) {\line(0,-1){10}} 
		}
		\thicklines
		\put(25,15) {\line(1,0){10}} 
		\put(25,10) {\line(1,0){10}} 
		%
		\put(25,15) {\line(0,-1){5}} 
		\put(30,15) {\line(0,-1){5}} 
		\put(35,15) {\line(0,-1){5}} 
		\put(0,0) {\line(1,0){5}} 
		\put(0,-5) {\line(1,0){5}} 
		\put(0,0) {\line(0,-1){5}} 
		\put(5,0) {\line(0,-1){5}} 
		{\color{dblue}
		\put(1,-4){2}
		\put(26,11){1}
		\put(31,11){3}}
		\end{picture}, \text{ and }
		T_3 = \begin{picture}(35,15)(0,2.5)
		{\color{grey}
			\multiput(0,0)(15,0){2} {\line(0,1){15}} 
		\multiput(0,0)(0,15){2} {\line(1,0){15}} 
		\put(15,15) {\line(1,0){10}} 
		\put(15,10) {\line(1,0){10}} 
		\put(15,5) {\line(1,0){10}} 
		\put(20,15) {\line(0,-1){10}} 
		\put(25,15) {\line(0,-1){10}} 
		}
		\thicklines
		\put(25,15) {\line(1,0){10}} 
		\put(25,10) {\line(1,0){10}} 
		%
		\put(25,15) {\line(0,-1){5}} 
		\put(30,15) {\line(0,-1){5}} 
		\put(35,15) {\line(0,-1){5}} 
		\put(0,0) {\line(1,0){5}} 
		\put(0,-5) {\line(1,0){5}} 
		\put(0,0) {\line(0,-1){5}} 
		\put(5,0) {\line(0,-1){5}} 
		{\color{dblue}
		\put(1,-4){1}
		\put(26,11){2}
		\put(31,11){3}}
		\end{picture}.
		\end{equation}
	}
Now consider sequences of partitions $T = (\lambda=T^{(0)}\searrow T^{(1)}\searrow \dots\searrow T^{(k)}=\mu)$ where $T^{(i)}$ is obtained from $T^{(i-1)}$ by adding a box. We can identify each  $T$ with the $\mu/\lambda$-tableau built by placing the integer $i$ in the box added at the $i^\text{th}$ step. For example, 
$$
	\begin{picture}(25,15)(0,0)
	\thicklines
		\multiput(0,0)(15,0){2} {\line(0,1){15}} 
		\multiput(0,0)(0,15){2} {\line(1,0){15}} 
		\put(15,15) {\line(1,0){10}} 
		\put(15,10) {\line(1,0){10}} 
		\put(15,5) {\line(1,0){10}} 
		\put(20,15) {\line(0,-1){10}} 
		\put(25,15) {\line(0,-1){10}} 
		\end{picture}
\searrow
	\begin{picture}(30,15)(0,0)
		\thicklines
		\multiput(0,0)(15,0){2} {\line(0,1){15}} 
		\multiput(0,0)(0,15){2} {\line(1,0){15}} 
		\put(15,15) {\line(1,0){15}} 
		\put(15,10) {\line(1,0){15}} 
		\put(15,5) {\line(1,0){10}} 
		\put(20,15) {\line(0,-1){10}} 
		\put(25,15) {\line(0,-1){10}} 
		\put(30,15) {\line(0,-1){5}} 
		\end{picture}
\searrow
	\begin{picture}(35,15)(0,0)
	\thicklines
		\multiput(0,0)(15,0){2} {\line(0,1){15}} 
		\multiput(0,0)(0,15){2} {\line(1,0){15}} 
		\put(15,15) {\line(1,0){20}} 
		\put(15,10) {\line(1,0){20}} 
		\put(15,5) {\line(1,0){10}} 
		\put(20,15) {\line(0,-1){10}} 
		\put(25,15) {\line(0,-1){10}} 
		\put(30,15) {\line(0,-1){5}} 
		\put(35,15) {\line(0,-1){5}} 
		\end{picture}
\searrow	
	\begin{picture}(35,15)(0,0)
		\thicklines
		\multiput(0,0)(15,0){2} {\line(0,1){15}} 
		\multiput(0,0)(0,15){2} {\line(1,0){15}} 
		\put(15,15) {\line(1,0){20}} 
		\put(15,10) {\line(1,0){20}} 
		\put(15,5) {\line(1,0){10}} 
		\put(20,15) {\line(0,-1){10}} 
		\put(25,15) {\line(0,-1){10}} 
		\put(30,15) {\line(0,-1){5}} 
		\put(35,15) {\line(0,-1){5}} 
		\put(0,-5) {\line(1,0){5}} 
		\put(0,0) {\line(0,-1){5}} 
		\put(5,0) {\line(0,-1){5}} 
		\end{picture}
$$
is identified with $T_1$ in \eqref{eq:tab_example}.


\begin{example}[Bratteli diagram for $\End_\fg(M\otimes V^{\otimes k})$]
\label{ex:BratMV}
Let $\fg = \fsl_n$ or $\fgl_n$, $M = L((a^p))$, and $V = L((1^1))$, and consider the example where $\cC_k = \End_\fg(M \otimes V^{\otimes {k}})$. The inclusion map in \eqref{eq:embedded_centralizers} provides a chain
\begin{equation}
\label{MV chain} \CC =  \End_\fg(M) \subseteq  \End_\fg(M\otimes V) \subseteq \cdots
\end{equation}
as in equation \eqref{chain}. By identifying classes of $\hat\cC_k$ with partitions as in Section \ref{sec:prelims}, we learn from Example \ref{ex:LR-add-a-box} that the dimensions in \eqref{multiplicity} are all 0 or 1, and that the paths in \eqref{eq:paths1} are in bijection with the set of
tableaux
\begin{equation}\label{eq:MV-tableaux}
\{ T= ((a^p) = T^{(0)} \searrow \cdots \searrow T^{(k)}) \}.
\end{equation} 
In particular, each $\cC_k$-module $\cC^\lambda$ has seminormal basis  $\{v_T\}$ indexed by the tableaux in \eqref{eq:MV-tableaux} which end at $T^{(k)} = \lambda$. Moreover, $\cC^\lambda$ is the same as $\cL^\lambda$ in \eqref{eq:centralizer_thm}, and each $v_T$ is a highest weight vector of weight $T^{(i)}$ in $\Res_{\cC_i}^{\cC_k} M \otimes V^{\otimes k}$ for each $i = 0, \dots, k$. 
\end{example}

\begin{example}[Bratteli diagram for $\End_\fg(M\otimes N \otimes V^{\otimes k})$] 
\label{ex:BratMNV}
Let $\fg = \fsl_n$ or $\fgl_n$, $M = L((a^p))$, $N=L((b^q))$, and $V = L((1^1))$, and consider the example where 
$$\cC_0 = \End_\fg(M) \quad \text{and} \quad \cC_k = \End_\fg(M \otimes N \otimes V^{\otimes {k-1}}).
$$
Just as in the previous example, these $\cC_k$ satisfy the chain 
\begin{equation}
\label{MNV chain} \CC =  \End_\fg(M) \subseteq  \End_\fg(M\otimes N) \subseteq \End_\fg(M\otimes N\otimes V) \subseteq \cdots.
\end{equation}
As in Example \ref{ex:rectangles}, if $(a^p)$ and $(b^q)$ are rectangular partitions, let $\cP=\cP((a^p),(b^q))$ be the set of partitions $\mu$ for which $L(\mu)$ appears as a submodule of  $L((a^p))\otimes L((b^q))$. In particular, each $L(\mu)$ appears with multiplicity 1. Let $\cP_0=\cP$ and define $\cP_k$ to be the set of partitions which are obtained by adding a box to an element of $\cP_{k-1}$. 

The classes in $\hat \cC_k$ are in bijection with the partitions in $\cP_{k-1}$, and the Bratteli diagram for the chain in \eqref{MNV chain} is the following oriented ranked graph:
\begin{enumerate}\item[]\begin{enumerate}
\item[Vertices:]The vertices are labeled by partitions. \begin{enumerate}
	\item[level $0$:]On level $0$, place one vertex, labeled by $(a^p)$. 
	\item[level $k>0$:]On level $k> 0$ place one vertex for each partition in $ \cP_{k-1}$. 
	\end{enumerate}
\item[Edges:]Edges connect two vertices on adjacent levels.  
	\item[]Connect the vertex on level $0$ to each of the vertices on level 1 with one edge. 
	\item[]Connect each vertex on level $k-1$ to a vertex on level $k$ if the vertex on level $k$ can be obtained by adding a box to the corresponding vertex on level $k-1$. 
\end{enumerate}
\end{enumerate}

\begin{figure}[!tp]
$$~~~~~
{\def\UNIT{.8pt}
\xymatrix{
&&
\setlength{\unitlength}{\UNIT}
\begin{picture}(25,20)(0,0)
	\put(10,23){\tiny $a$}
	\put(2,9){\tiny $p$}
	\thicklines
		\multiput(0,0)(25,0){2} {\line(0,1){20}} 
		\multiput(0,0)(0,20){2} {\line(1,0){25}} 
		\end{picture}
		\ar@/_1.5pc/[ddll]|{\color{dblue} 4a}
		\ar@/_1pc/[ddl]|{\color{dblue} 3a-p}
		\ar@/_.5pc/[dd]|(.55){\color{dblue} 2(a-p+1)}
		\ar@/^.5pc/[ddr]|{\color{dblue} 2(a-p-1)}
		\ar@/^1pc/[ddrr]|{\color{dblue} a-3p}
		\ar@/^1.5pc/[ddrrr]|{\color{dblue} -4p}
&&&
\\
~&&&&&
\\
\setlength{\unitlength}{\UNIT}
\begin{picture}(35,30)(0,-10)
	\thicklines
		\multiput(0,0)(25,0){2} {\line(0,1){20}} 
		\multiput(0,0)(0,20){2} {\line(1,0){25}} 
		\put(25,20) {\line(1,0){10}} 
		\put(25,15) {\line(1,0){10}} 
		\put(25,10) {\line(1,0){10}} 
		\put(30,20) {\line(0,-1){10}} 
		\put(35,20) {\line(0,-1){10}} 
		\end{picture}
		\ar[dd]|{ \color{dblue} -p}
		\ar@/_2pc/[ddd]|(.7){\color{dblue} a+2}
		\ar@/_3pc/[dddd]|(.75){\color{dblue} a-2}
 &
\setlength{\unitlength}{\UNIT}
\begin{picture}(35, 30)(0,-10)
	\thicklines
		\multiput(0,0)(25,0){2} {\line(0,1){20}} 
		\multiput(0,0)(0,20){2} {\line(1,0){25}} 
		\put(25,20) {\line(1,0){10}} 
		\put(25,15) {\line(1,0){10}} 
		\put(25,10) {\line(1,0){5}} 
		\put(30,20) {\line(0,-1){10}} 
		\put(35,20) {\line(0,-1){5}} 
		\put(0,-5) {\line(1,0){5}} 
		\put(0,0) {\line(0,-1){5}} 
		\put(5,0) {\line(0,-1){5}} 
		\end{picture}
		\ar[ddl]|{\color{dblue} a \!\!\phantom{|}} \ar[dd]|{\color{dblue} -p+1} 
		\ar@{->}[ddr]|(.33){\color{dblue} -p-1}
		\ar@/_2pc/[ddd]|(.7){\color{dblue} a+2}
		\ar@/_3pc/[dddd]|(.75){\color{dblue} a-2}
 &
\setlength{\unitlength}{\UNIT}
\begin{picture}(35, 30)(0,-10)
	\thicklines
		\multiput(0,0)(25,0){2} {\line(0,1){20}} 
		\multiput(0,0)(0,20){2} {\line(1,0){25}} 
		\put(25,20) {\line(1,0){10}} 
		\put(25,15) {\line(1,0){10}} 
		\put(30,20) {\line(0,-1){5}} 
		\put(35,20) {\line(0,-1){5}} 
		\put(0,-5) {\line(1,0){10}} 
		\put(0,0) {\line(0,-1){5}} 
		\put(5,0) {\line(0,-1){5}} 
		\put(10,0) {\line(0,-1){5}} 
		\end{picture}
		\ar[ddl]|(.66){\color{dblue} a-1}
		\ar[ddr]|(.33){\color{dblue} -p-1}
		\ar@/_2pc/[ddd]|(.7){\color{dblue} a+2}
		\ar@/_3pc/[dddd]|(.75){\color{dblue} -p+2}
 &
\setlength{\unitlength}{\UNIT}
\begin{picture}(30, 30)(0,-10)
	\thicklines
		\multiput(0,0)(25,0){2} {\line(0,1){20}} 
		\multiput(0,0)(0,20){2} {\line(1,0){25}} 
		\put(25,20) {\line(1,0){5}} 
		\put(25,15) {\line(1,0){5}} 
		\put(25,10) {\line(1,0){5}} 
		\put(30,20) {\line(0,-1){10}} 
		\put(0,-5) {\line(1,0){5}} 
		\put(0,-10) {\line(1,0){5}} 
		\put(0,0) {\line(0,-1){10}} 
		\put(5,0) {\line(0,-1){10}} 
     \end{picture}
		\ar[ddl]|(.66){\color{dblue} a+1}
		\ar[ddr]|(.66){\color{dblue} -p+1}
		\ar@/^2pc/[ddd]|(.7){\color{dblue} a-2}
		\ar@/^3pc/[dddd]|(.75){\color{dblue} -p-2}
 &
\setlength{\unitlength}{\UNIT}
\begin{picture}(35, 30)(0,-10)
	\thicklines
		\multiput(0,0)(25,0){2} {\line(0,1){20}} 
		\multiput(0,0)(0,20){2} {\line(1,0){25}} 
		\put(25,20) {\line(1,0){5}} 
		\put(25,15) {\line(1,0){5}} 
		\put(30,20) {\line(0,-1){5}} 
		\put(0,-5) {\line(1,0){10}} 
		\put(0,-10) {\line(1,0){5}} 
		\put(0,0) {\line(0,-1){10}} 
		\put(5,0) {\line(0,-1){10}} 
		\put(10,0) {\line(0,-1){5}} 
     \end{picture} 
		\ar[ddl]|(.33){\color{dblue} a+1}
		\ar[dd]|{\color{dblue} a-1}
		\ar[ddr]|(.6){\color{dblue} -p}
		\ar@/^2pc/[ddd]|(.7){\color{dblue} -p+2}
		\ar@/^3pc/[dddd]|(.75){\color{dblue} -p-2}
 &
\setlength{\unitlength}{\UNIT}
\begin{picture}(30, 30)(0,-10)
	\thicklines
		\multiput(0,0)(25,0){2} {\line(0,1){20}} 
		\multiput(0,0)(0,20){2} {\line(1,0){25}} 
		\put(0,-5) {\line(1,0){10}} 
		\put(0,-10) {\line(1,0){10}} 
		\put(0,0) {\line(0,-1){10}} 
		\put(5,0) {\line(0,-1){10}} 
		\put(10,0) {\line(0,-1){10}} 
     \end{picture}
		\ar[dd]|{\color{dblue} a \!\!\phantom{|}}
		\ar@/^2pc/[ddd]|(.7){\color{dblue} -p+2}
		\ar@/^3pc/[dddd]|(.75){\color{dblue} -p-2}
     \\
     ~&&&&&
     \\
\setlength{\unitlength}{\UNIT}
\begin{picture}(35, 30)(0,-10)
	\thicklines
		\multiput(0,0)(25,0){2} {\line(0,1){20}} 
		\multiput(0,0)(0,20){2} {\line(1,0){25}} 
		\put(25,20) {\line(1,0){10}} 
		\put(25,15) {\line(1,0){10}} 
		\put(25,10) {\line(1,0){10}} 
		\put(30,20) {\line(0,-1){10}} 
		\put(35,20) {\line(0,-1){10}} 
		\put(0,-5) {\line(1,0){5}} 
		\put(0,0) {\line(0,-1){5}} 
		\put(5,0) {\line(0,-1){5}} 
		%
     \end{picture} 
&
\setlength{\unitlength}{\UNIT}
\begin{picture}(35, 30)(0,-10)
	\thicklines
		\multiput(0,0)(25,0){2} {\line(0,1){20}} 
		\multiput(0,0)(0,20){2} {\line(1,0){25}} 
		\put(25,20) {\line(1,0){10}} 
		\put(25,15) {\line(1,0){10}} 
		\put(25,10) {\line(1,0){5}} 
		\put(30,20) {\line(0,-1){10}} 
		\put(35,20) {\line(0,-1){5}} 
		\put(0,-5) {\line(1,0){10}} 
		\put(0,0) {\line(0,-1){5}} 
		\put(5,0) {\line(0,-1){5}} 
		\put(10,0) {\line(0,-1){5}} 
		%
		\end{picture} 
&
\setlength{\unitlength}{\UNIT}
\begin{picture}(30, 30)(0,-10)
	\thicklines
		\multiput(0,0)(25,0){2} {\line(0,1){20}} 
		\multiput(0,0)(0,20){2} {\line(1,0){25}} 
		\put(25,20) {\line(1,0){10}} 
		\put(25,15) {\line(1,0){10}} 
		\put(25,10) {\line(1,0){5}} 
		\put(30,20) {\line(0,-1){10}} 
		\put(35,20) {\line(0,-1){5}} 
		\put(0,-5) {\line(1,0){5}} 
		\put(0,-10) {\line(1,0){5}} 
		\put(0,0) {\line(0,-1){10}} 
		\put(5,0) {\line(0,-1){10}} 
		%
		\end{picture} 
&
\setlength{\unitlength}{\UNIT}
\begin{picture}(30, 30)(0,-10)
	\thicklines
		\multiput(0,0)(25,0){2} {\line(0,1){20}} 
		\multiput(0,0)(0,20){2} {\line(1,0){25}} 
		\put(25,20) {\line(1,0){10}} 
		\put(25,15) {\line(1,0){10}} 
		\put(30,20) {\line(0,-1){5}} 
		\put(35,20) {\line(0,-1){5}} 
		\put(0,-5) {\line(1,0){10}} 
		\put(0,-10) {\line(1,0){5}} 
		\put(0,0) {\line(0,-1){10}} 
		\put(5,0) {\line(0,-1){10}} 
		\put(10,0) {\line(0,-1){5}} 
		%
		\end{picture} 
&
\setlength{\unitlength}{\UNIT}
\begin{picture}(35, 30)(0,-10)
	\thicklines
		\multiput(0,0)(25,0){2} {\line(0,1){20}} 
		\multiput(0,0)(0,20){2} {\line(1,0){25}} 
		\put(25,20) {\line(1,0){5}} 
		\put(25,15) {\line(1,0){5}} 
		\put(25,10) {\line(1,0){5}} 
		\put(30,20) {\line(0,-1){10}} 
		\put(0,-5) {\line(1,0){10}} 
		\put(0,-10) {\line(1,0){5}} 
		\put(0,0) {\line(0,-1){10}} 
		\put(5,0) {\line(0,-1){10}} 
		\put(10,0) {\line(0,-1){5}} 
		%
		\end{picture} 
&
\setlength{\unitlength}{\UNIT}
\begin{picture}(30, 30)(0,-10)
	\thicklines
		\multiput(0,0)(25,0){2} {\line(0,1){20}} 
		\multiput(0,0)(0,20){2} {\line(1,0){25}} 
		\put(25,20) {\line(1,0){5}} 
		\put(25,15) {\line(1,0){5}} 
		\put(30,20) {\line(0,-1){5}} 
		\put(0,-5) {\line(1,0){10}} 
		\put(0,-10) {\line(1,0){10}} 
		\put(0,0) {\line(0,-1){10}} 
		\put(5,0) {\line(0,-1){10}} 
		\put(10,0) {\line(0,-1){10}} 
		%
		\end{picture} 
 \\
\setlength{\unitlength}{\UNIT}
\begin{picture}(40,30)(0,-10)
	\thicklines
		\multiput(0,0)(25,0){2} {\line(0,1){20}} 
		\multiput(0,0)(0,20){2} {\line(1,0){25}} 
		\put(25,20) {\line(1,0){15}} 
		\put(25,15) {\line(1,0){15}} 
		\put(25,10) {\line(1,0){10}} 
		\put(30,20) {\line(0,-1){10}} 
		\put(35,20) {\line(0,-1){10}} 
		\put(40,20) {\line(0,-1){5}} 
		%
		\end{picture}
 &
\setlength{\unitlength}{\UNIT}
\begin{picture}(40,30)(0,-10)
	\thicklines
		\multiput(0,0)(25,0){2} {\line(0,1){20}} 
		\multiput(0,0)(0,20){2} {\line(1,0){25}} 
		\put(25,20) {\line(1,0){15}} 
		\put(25,15) {\line(1,0){15}} 
		\put(25,10) {\line(1,0){5}} 
		\put(30,20) {\line(0,-1){10}} 
		\put(35,20) {\line(0,-1){5}} 
		\put(40,20) {\line(0,-1){5}} 
		\put(0,-5) {\line(1,0){5}} 
		\put(0,0) {\line(0,-1){5}} 
		\put(5,0) {\line(0,-1){5}} 
		%
		\end{picture}
 &
\setlength{\unitlength}{\UNIT}
\begin{picture}(30,30)(0,-10)
	\thicklines
		\multiput(0,0)(25,0){2} {\line(0,1){20}} 
		\multiput(0,0)(0,20){2} {\line(1,0){25}} 
		\put(25,20) {\line(1,0){15}} 
		\put(25,15) {\line(1,0){15}} 
		\put(30,20) {\line(0,-1){5}} 
		\put(35,20) {\line(0,-1){5}} 
		\put(40,20) {\line(0,-1){5}} 
		\put(0,-5) {\line(1,0){10}} 
		\put(0,0) {\line(0,-1){5}} 
		\put(5,0) {\line(0,-1){5}} 
		\put(10,0) {\line(0,-1){5}} 
		%
		\end{picture}
 &
\setlength{\unitlength}{\UNIT}
\begin{picture}(30,30)(0,-10)
	\thicklines
		\multiput(0,0)(25,0){2} {\line(0,1){20}} 
		\multiput(0,0)(0,20){2} {\line(1,0){25}} 
		\put(25,20) {\line(1,0){5}} 
		\put(25,15) {\line(1,0){5}} 
		\put(25,10) {\line(1,0){5}} 
		\put(25,5) {\line(1,0){5}} 
		\put(30,20) {\line(0,-1){15}} 
		\put(0,-5) {\line(1,0){5}} 
		\put(0,-10) {\line(1,0){5}} 
		\put(0,0) {\line(0,-1){10}} 
		\put(5,0) {\line(0,-1){10}} 
		%
     \end{picture}
 &
\setlength{\unitlength}{\UNIT}
\begin{picture}(40,30)(0,-10)
	\thicklines
		\multiput(0,0)(25,0){2} {\line(0,1){20}} 
		\multiput(0,0)(0,20){2} {\line(1,0){25}} 
		\put(25,20) {\line(1,0){5}} 
		\put(25,15) {\line(1,0){5}} 
		\put(30,20) {\line(0,-1){5}} 
		\put(0,-5) {\line(1,0){15}} 
		\put(0,-10) {\line(1,0){5}} 
		\put(0,0) {\line(0,-1){10}} 
		\put(5,0) {\line(0,-1){10}} 
		\put(10,0) {\line(0,-1){5}} 
		\put(15,0) {\line(0,-1){5}} 
		%
     \end{picture} 
 &
\setlength{\unitlength}{\UNIT}
\begin{picture}(30,30)(0,-10)
	\thicklines
		\multiput(0,0)(25,0){2} {\line(0,1){20}} 
		\multiput(0,0)(0,20){2} {\line(1,0){25}} 
		\put(0,-5) {\line(1,0){15}} 
		\put(0,-10) {\line(1,0){10}} 
		\put(0,0) {\line(0,-1){10}} 
		\put(5,0) {\line(0,-1){10}} 
		\put(10,0) {\line(0,-1){10}} 
		\put(15,0) {\line(0,-1){5}} 
		%
\end{picture}
\\
\setlength{\unitlength}{\UNIT}
\begin{picture}(35,35)(0,-15)
	\thicklines
		\multiput(0,0)(25,0){2} {\line(0,1){20}} 
		\multiput(0,0)(0,20){2} {\line(1,0){25}} 
		\put(25,20) {\line(1,0){10}} 
		\put(25,15) {\line(1,0){10}} 
		\put(25,10) {\line(1,0){10}} 
		\put(25,5) {\line(1,0){5}} 
		\put(30,20) {\line(0,-1){15}} 
		\put(35,20) {\line(0,-1){10}} 
		%
		\end{picture}
 &
\setlength{\unitlength}{\UNIT}
\begin{picture}(35,35)(0,-15)
	\put(2,2){\tiny 6}
	\thicklines
		\multiput(0,0)(25,0){2} {\line(0,1){20}} 
		\multiput(0,0)(0,20){2} {\line(1,0){25}} 
		\put(25,20) {\line(1,0){10}} 
		\put(25,15) {\line(1,0){10}} 
		\put(25,10) {\line(1,0){5}} 
		\put(25,5) {\line(1,0){5}} 
		\put(30,20) {\line(0,-1){15}} 
		\put(35,20) {\line(0,-1){5}} 
		\put(0,-5) {\line(1,0){5}} 
		\put(0,0) {\line(0,-1){5}} 
		\put(5,0) {\line(0,-1){5}} 
		%
		\end{picture}
 &
\setlength{\unitlength}{\UNIT}
\begin{picture}(30,35)(0,-15)
	\thicklines
		\multiput(0,0)(25,0){2} {\line(0,1){20}} 
		\multiput(0,0)(0,20){2} {\line(1,0){25}} 
		\put(25,20) {\line(1,0){10}} 
		\put(25,15) {\line(1,0){10}} 
		\put(30,20) {\line(0,-1){5}} 
		\put(35,20) {\line(0,-1){5}} 
		\put(0,-5) {\line(1,0){15}} 
		\put(0,0) {\line(0,-1){5}} 
		\put(5,0) {\line(0,-1){5}} 
		\put(10,0) {\line(0,-1){5}} 
		\put(15,0) {\line(0,-1){5}} 
		%
		\end{picture}
 &
\setlength{\unitlength}{\UNIT}
\begin{picture}(30,35)(0,-15)
	\thicklines
		\multiput(0,0)(25,0){2} {\line(0,1){20}} 
		\multiput(0,0)(0,20){2} {\line(1,0){25}} 
		\put(25,20) {\line(1,0){5}} 
		\put(25,15) {\line(1,0){5}} 
		\put(25,10) {\line(1,0){5}} 
		\put(30,20) {\line(0,-1){10}} 
		\put(0,-5) {\line(1,0){5}} 
		\put(0,-10) {\line(1,0){5}} 
		\put(0,-15) {\line(1,0){5}} 
		\put(0,0) {\line(0,-1){15}} 
		\put(5,0) {\line(0,-1){15}} 
		%
     \end{picture}
 &
\setlength{\unitlength}{\UNIT}
\begin{picture}(35,35)(0,-15)
	\thicklines
		\multiput(0,0)(25,0){2} {\line(0,1){20}} 
		\multiput(0,0)(0,20){2} {\line(1,0){25}} 
		\put(25,20) {\line(1,0){5}} 
		\put(25,15) {\line(1,0){5}} 
		\put(30,20) {\line(0,-1){5}} 
		\put(0,-5) {\line(1,0){10}} 
		\put(0,-10) {\line(1,0){5}} 
		\put(0,-15) {\line(1,0){5}} 
		\put(0,0) {\line(0,-1){15}} 
		\put(5,0) {\line(0,-1){15}} 
		\put(10,0) {\line(0,-1){5}} 
		%
     \end{picture} 
 &
\setlength{\unitlength}{\UNIT}
\begin{picture}(30,35)(0,-15)
	\thicklines
		\multiput(0,0)(25,0){2} {\line(0,1){20}} 
		\multiput(0,0)(0,20){2} {\line(1,0){25}} 
		\put(0,-5) {\line(1,0){10}} 
		\put(0,-10) {\line(1,0){10}} 
		\put(0,-15) {\line(1,0){5}} 
		\put(0,0) {\line(0,-1){15}} 
		\put(5,0) {\line(0,-1){15}} 
		\put(10,0) {\line(0,-1){10}} 
		%
     \end{picture}
    }
    }
$$
\caption{Levels 0, 1, and 2 of a Bratteli diagram encoding isotypic components of $M \otimes N \otimes V$. The edges are labeled by combinatorial values given by the action of $\cH_k^\ext$ as stated later in Theorem \ref{thm:Bratteli_action}.}
\label{fig:Brat_diagram_with_contents}
\end{figure}
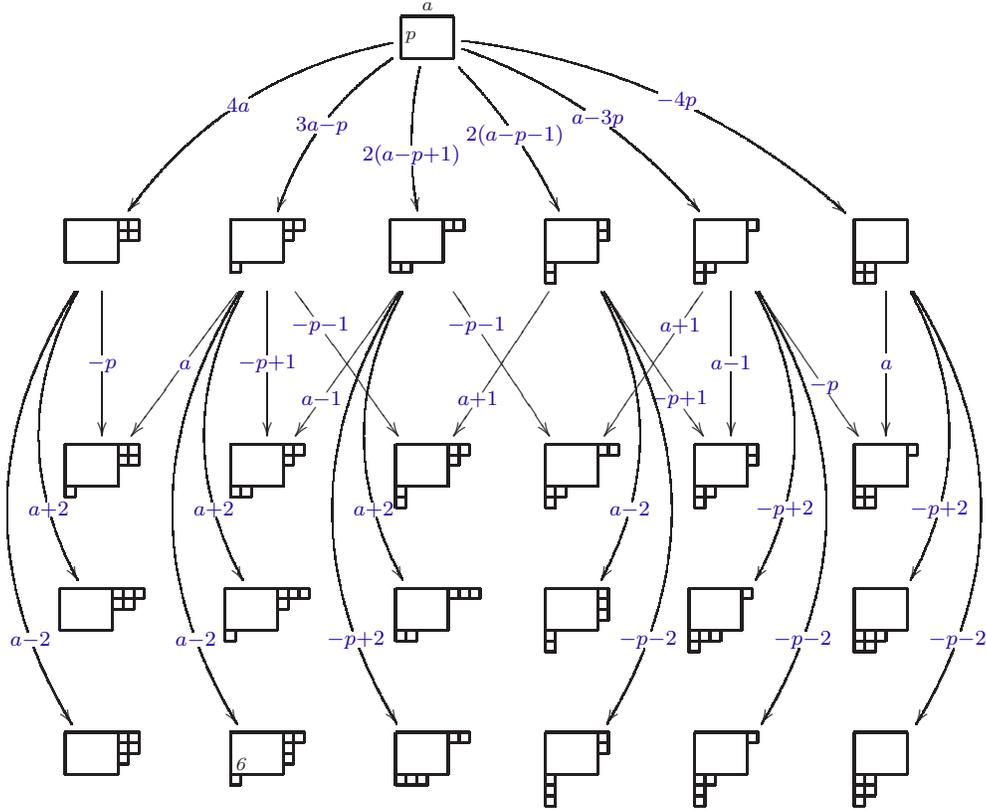

\noindent For the examples where  $a,p > 2$ and $b=q=2$, levels 0, 1, and 2 are shown of the Bratteli diagram in Figure \ref{fig:Brat_diagram_with_contents}. Notice that each of the partitions in $\cP_1$ comes from exactly one or two partitions in $\cP_0$; this happens for any choice of $a,b,p,$ and $q$ by Lemma \ref{lem:one_or_two}.

Again, we learn from Examples \ref{ex:LR-add-a-box} and \ref{ex:rectangles} that the dimensions in \eqref{multiplicity} are all 0 or 1, and that the paths in \eqref{eq:paths1} are in bijection with the set of tableaux 
	\begin{equation}\label{eq:MNV-tableaux}
	\{ T = ( T^{(1)} \searrow \cdots \searrow T^{({k})}) ~|~ T^{(i)} \in \cP_{i-1} \text{ for }i=1, \dots, k\}.\end{equation}
	In particular, for each $\lambda \in \cP_{k-1}$, the $\cC_k$-module $\cC^\lambda$ has seminormal basis  $\{v_T\}$ indexed by the tableaux in \eqref{eq:MNV-tableaux} which end at $T^{(k)} = \lambda$. Moreover, $\cC^\lambda$ is the same as $\cL^\lambda$ in \eqref{eq:centralizer_thm}, and each $v_T$ is a highest weight vector of weight $T^{(i)}$ in $\Res_{\cC_i}^{\cC_k} M \otimes N \otimes V^{\otimes k-1}$ for each $i = 0, \dots, k-1$. 

\end{example}

We now return to the representation $\Phi'$ in Theorem \ref{thm:Hecke-action}, and combinatorially determine the eigenvalues of the operators $\Phi'(x_i), \Phi'(y_i)$, and $\Phi'(z_i)$, for $i= 1, \dots, k$. If $\lambda \subset \mu$ are partitions differing by a box in column $c$ and row $r$, recall $c(\mu/\lambda) = c - r$
is the content of the box $\mu/\lambda$ in $\mu$. 

\begin{thm}~\label{thm:Bratteli_action}  Let  $\Phi': \cH_k^\ext \to \End_\fg(M \otimes N \otimes V^{\otimes k})$ be the representation in Theorem \ref{thm:Hecke-action}, with
$$c_z = \begin{cases} 0 & \text{if } \fg= \fgl_n\\
	\frac{abpq}{n}&  \text{if } \fg= \fsl_n.\end{cases}$$
 There is a basis $\{v_T\}$ of  $M\otimes N \otimes V^{\otimes k}$ indexed by standard tableaux 
$$	\{ T = ( T^{(0)} \searrow \cdots \searrow T^{(k)}) ~|~ T^{(j)} \in \cP_{j} \text{ for }j=0, \dots, k\}$$
with action 
$$\Phi'(z_i) v_T = c(T^{(i)}/T^{(i-1)})v_T, \qquad \text{for }i=1, \dots, k,$$
and $$\Phi'(z_0) v_T = \left(abq + 2 \sum_{B \in \cB_{T^{(0)}}} \left(c(B) - \half(a-p + b-q)\right)\right)v_T,$$
	where $\cB_{\lambda}$ to be the set of boxes in $\lambda$ in rows $p+1$ and below.
\end{thm}

\begin{proof}

The basis for $M \otimes N \otimes V^{\otimes k}$ produced in Example \ref{ex:BratMNV} is specifically one which satisfies
	$$v_T \in v_{T^{(i)}} \otimes V^{\otimes (k-i)} \subseteq M \otimes N \otimes V^{\otimes k} ,\quad i=0, \dots, k,$$ 
where $v_{T^{(i)}}$ is a highest weight vector of weight ${T^{(i)}}$ in $L(T^{(i)}) \subseteq M \otimes N \otimes V^{\otimes i}$.
Therefore 
\begin{align*}
\Phi'(z_i) \cdot v_T 
	&= \left( \half(\kappa_{M\otimes N,\leq i} -\kappa_{M\otimes N,\leq i-1} + \kappa_V) + ((i-1)d + c_x + c_y)\id\right)\cdot v_T\\
	&= \left( \gamma_{L(T^{(i-1)}), V} + \kappa_V  + ((i-1)d + c_x + c_y)\id\right)v_T.
\end{align*} 
By Lemma \ref{thm:gamma_contents}, \eqref{kappaV-gl} and \eqref{kappaV-sl},
$$
\gamma_{L(T^{(i-1)}), V} = \begin{cases}  
		c(T^{(i)}/T^{(i-1)}) & \text{ if } \fg=\fgl_n,\\
		c(T^{(i)}/T^{(i-1)}) - \frac{ap-bq+i-1}{n} & \text{ if } \fg=\fsl_n,\end{cases}
\quad \text{ and } \quad
\kappa_V=\begin{cases} 
		n & \text{ if } \fg=\fgl_n,\\ 
		n-\frac{1}{n} & \text{ if } \fg=\fsl_n.\end{cases}$$  So since 
		$$c_x = c_y = -\half n, \text{ and } d=0 \qquad \text{ when $\fg=\fgl_n$, and}$$
		$$c_x= \frac{ap}{n}- \frac12\left(n-\frac{1}{n}\right) ,  
	c_y =  \frac{bq}{n} - \frac{1}{2}\left(n-\frac{1}{n}\right), \text{ and }
	d= \frac{1}{n} \qquad \text{ when $\fg=\fsl_n$,}$$
we have $\Phi'(z_i) \cdot v_T = c(T^{(i)}/T^{(i-1)})v_T$ as desired. Similarly, the action of $\Phi'(z_0)$ follows from Lemma \ref{lem:z_0_on_lambda}. 
\end{proof}

\begin{example}

To illustrate, we apply Theorem \ref{thm:Bratteli_action} to the example where $a,p>2$, $b=q=2$, and $k=1$.  Returning to Figure \ref{fig:Brat_diagram_with_contents} above, we can read that there there are eighteen distinct isotypic components of $M \otimes N \otimes V$, six of which correspond to 2-dimensional $\End_\fg(M\otimes N\otimes V)$-modules and twelve of which correspond to 1-dimensional $\End_\fg(M\otimes N\otimes V)$-modules. 

The edges connecting level $0$ to level $1$ are labeled by the combinatorial formula for the action of $z_0$, and the edges connecting level $1$ to level $2$ are labeled by the content of the box added. In general, we label the edges connecting level $i$ to level $i+1$ by the content of the box added. The paths in this diagram from $(a^p)$ to $\lambda \in \cP_1$ index the basis of $\cL^\lambda$, and $\Phi'(z_1)$ and $\Phi'(z_0)$ act on those basis elements by the corresponding edge labels. 
%

\end{example}

\begin{remark}\label{rk:x_1-evalues}
Example \ref{ex:BratMV} gives a basis $\{v_{T}\}$ of  $M\otimes V^{\otimes k}$ indexed by standard tableaux 
$$	\{ T = ( (a^p) = T^{(0)} \searrow \cdots \searrow T^{(k)}) \}.$$
For every $n \in N$, there is a canonical map 
$$\iota_n: M \otimes V^{\otimes k} \stackrel{}{\hookrightarrow} M \otimes V^{\otimes k} \otimes N \cong M \otimes N \otimes V^{\otimes k}.$$
Therefore, by picking a basis $\{n_j\}$ of $N$, $\{v_T\}$ can be lifted to a basis for $M \otimes N \otimes V^{\otimes k}$
$$ \{\iota_{n_j}(v_T)\}_{j, T}.$$
A similar calculation as in Theorem \ref{thm:Bratteli_action} will produce 
$$\Phi'(x_i) \cdot \iota_{n_j} (v_T) = c( T^{(i)}/ T^{(i-1)})\iota_{n_j} (v_T).$$

Similarly, there is a basis of highest weight vectors $\{v_T\}$ for $N \otimes V^{\otimes k}$ indexed by 
standard tableaux 
$$	\{ T = ( (b^q) = T^{(0)} \searrow \cdots \searrow T^{(k)}) \}.$$ By picking a basis $\{m_j\}$ for $M$, the map
$$\iota_m: N \otimes V^{\otimes k} \stackrel{}{\hookrightarrow} M \otimes N \otimes V^{\otimes k}$$
 produces a basis 
$$\{\iota_{m_j}(v_T)\}_{j, T}$$
which satisfies
$$\Phi'(y_i) \cdot \iota_{m_j}(v_T) = c( T^{(i)}/ T^{(i-1)})\iota_{m_j}(v_T).$$

These two examples reflect the fact that the degenerate two-boundary braid algebras and Hecke algebras contain one-boundary analogs (though this isomorphic containment is left for future work).  Later, Theorem \ref{thm:seminormal-Hecke-values} will provide explicit formulas for $x_1$ (and therefore $y_1$) in terms of the basis given in Theorem \ref{thm:Bratteli_action}, but we can already ascertain the eigenvalues of $x_1$ and $y_1$, as they are stable under a change of basis. 
\end{remark}

\begin{remark}
The algebras in \eqref{chain} are also known as \emph{Gelfand-Zetlin algebras}, and the seminormal bases are (non-normalized) \emph{Gelfand-Zetlin bases}. Denote the center of $\cC_i$ by $Z(\cC_i)$. The commutative subalgebra $\cA_k \subseteq \cC_k$ generated by the subalgebras $Z(\cC_0)$, $Z(\cC_1)$, \dots $Z(\cC_k)$, is the  \emph{Gelfand-Zetlin subalgebra} in \cite{OV}. It remains for future work to show that the subalgebra of $\End(M \otimes N \otimes V^{\otimes k})$ generated by $\Phi'(z_0), \Phi'(z_1), \dots, \Phi'(z_{k})$ has large index inside of the Gelfand-Zetlin subalgebra for the chain in Example \ref{ex:BratMNV}. However, Theorem \ref{thm:Bratteli_action} is suggestive of this relationship, and future work on the center of $\cH_k^\ext$ will show that the subalgebra generated by subalgebras $Z(\cH_0^\ext)$, $Z(\cH_1^\ext)$,  \dots, $Z(\cH_k^\ext)$, is in fact the same as the subalgebra generated by $z_0, z_1, \dots, z_k$. 

\end{remark}

\section{Seminormal Representations of $\cH_k^\ext$}
\label{sec:seminormal_reps}

In Section \ref{sec:hecke-action}, we showed that a quotient of $\cH_k^\ext$ is a subalgebra 
of $\End_\fg(M \otimes N \otimes V^{\otimes k})$, when $\fg= \fgl_n$ or $\fsl_n$, $M$ and $N$ are simple $\fg$-modules indexed by rectangular partitions, and $V$ is the standard representation. Then in Section \ref{sec:seminormal_bases}, we showed that the action of the generators $z_0, \dots, z_k$ on $M \otimes N \otimes V^{\otimes k}$ is simultaneously diagonalizable with eigenvalues given by combinatorial values. In this section, we study all seminormal representations, and conclude finally in Corollary \ref{cor:punchline} that the simple $\End_\fg(M \otimes N \otimes V^{\otimes k})$-modules in $M \otimes N \otimes V^{\otimes k}$ are also simple as $\cH_k$-modules. 

This section serves as a culmination of work done so far on $\cH_k$ and $\cH_k^\ext$, and will draw on many results throughout the paper. As a guide to the reader, we will be primarily citing results from Sections \ref{sec:Hecke}, \ref{sec:prelims} and \ref{sec:seminormal_bases}. The presentation of choice for $\cH_k^\ext$ is given in Theorem \ref{thm:hecke_ext_presentation-short}; in particular, we switch from the generating set 
$$z_0, z_1, \dots, z_k, ~ x_1, \dots, x_k,~ y_1, \dots, y_k,~ t_{s_1}, \dots, t_{s_{k-1}}$$
to the generating set
$$w_0, w_1, \dots, w_k, ~ x_1, ~  t_{s_1}, \dots, t_{s_{k-1}},$$ 
where $w_i = z_i - \frac12(a-p+b-q)$. Section \ref{sec:seminormal_bases} (specifically Example \ref{ex:BratMNV}) introduces the combinatorial backbone of the modules that we study in this section; the bases for the modules in Proposition \ref{thm:seminormal-Hecke} are indexed by the same tableaux as in the Bratteli diagram for $\End_\fg(M \otimes N \otimes V^{\otimes k})$ in Example \ref{ex:BratMNV}. The specific combinatorial properties of these tableaux begin in Section \ref{sec:prelims}: Example \ref{ex:rectangles} describes the set  $\cP$ of partitions that index the simple submodules of $L((a^p)) \otimes L((p^q))$, and Figure \ref{fig:vis-of-rectangles} provides a useful illustration of the partitions in $\cP$; then Lemmas \ref{lem:one_or_two} and \ref{lem:content-distinction} tell us about the shape and symmetries of the Bratteli diagram at levels 0-2, and Figure  \ref{fig:added-boxes-for-2Ds} illustrates this symmetry. Finally, in Section \ref{sec:seminormal_bases}, Theorem \ref{thm:Bratteli_action} tells us the correct action of the $w_i$, and Remark \ref{rk:x_1-evalues} tells us how to anticipate the eigenvalues of the action of $x_1$.

\medskip 
\noindent Fix $a,b,p,q$ non-negative integers with $q \leq p$. 
Recall from Example \ref{ex:BratMNV} that $\cP_0=\cP$ is the set of partitions indexing simple submodules of $L((a^p))\otimes L((b^q))$, and $\cP_i$ is the set of partitions obtained by adding a box to any partition in $\cP_{i-1}$. Let $\cT_\lambda$ be the set of tableaux
\begin{equation}
\label{eq:indexing-tableau}\cT_\lambda = \left\{ T = (T^{(0)} \searrow  \ldots \searrow T^{(k)}=\lambda) ~|~ 
	T^{(0)} \in \cP,~ T^{(i)} \in \cP_{i}\right\}.\end{equation}
The box added to $T^{(i-1)}$ to get $T^{(i)}$ is $b_i=T^{(i)}/T^{(i-1)}$. Define \emph{shifted contents}
\begin{align*} 
	c_T(0) &= abq - (|\cB_\mu| + \half)(a-p+b-q) + 2\sum_{B \in \cB_\mu} c(B),\\
	c_T(i) &= c(T^{(i)}/T^{(i-1)}) - \half(a-p+b-q),
\end{align*} 
where $\cB_\mu$ is the set of boxes in $\mu$ in rows $p+1$ and below, as described in Lemma \ref{lem:z_0_on_lambda}. 
\begin{lemma}
\label{lem:content_list}
A tableau $T\in\cT_\lambda$ is fully determined by 
	$$ c_T(1), \dots, c_T(k), \quad \text{ and } \quad T^{(k)}.$$
 \end{lemma}
\begin{proof}
This can be shown by induction on $k$. The key observation is that the value $c_T(i)$, $i>0$, determines the diagonal on which $T^{(i)}/T^{(i-1)}$ lies. In any given partition, there is at most one removable box on any diagonal. So $c_T(k)$ and $T^{(k)}$ determines $T^{(k-1)}$. By iterating, $c_T(i)$ and $T^{(i)}$ determines $T^{(i-1)}$, so we can recover $T^{(k-1)}$, $T^{(k-2)}$, \dots,  $T^{(0)}$.   

\end{proof}

%
%

Two consecutive boxes $b_i$ and $b_{i+1}$ are in the same row or column if and only if $c(b_i) = c(b_{i+1}) \pm 1$. So 
for any $i$ for which $c_T(i) \neq c_T(i+1) \pm 1$, we can define $s_iT^{(i)}$ as the partition obtained removing $b_i$ and adding $b_{i+1}$, and so
	\begin{equation}\label{def:s_i}s_iT =  (T^{(0)}\searrow T^{(1)}\searrow \ldots \searrow T^{(i-1)}\searrow s_iT^{(i)}\searrow T^{(i+1)}\searrow \ldots\searrow T^{(k)})\end{equation}
	is the tableau constructed from $T$ by switching the order of adding the $i^\text{th}$ and $(i+1)^\text{th}$ boxes. Notice that if  $c_T(i) \neq c_T(i+1) \pm 1$, then $s_iT$ is the only tableau which varies from $T$ only at the $i^\text{th}$ position; otherwise, if  $c_T(i) = c_T(i+1) \pm 1$, then there is no such tableau. 
	
Similarly, for any $\mu \in \cP_1$, there are exactly one or two partitions $\nu \in \cP$ which differ from $\mu$ by a box by Lemma \ref{lem:one_or_two}. In other words, there are exactly one or two $\nu \in \cP$ which could be the first step in a tableau with a given shifted content list $c_T(2), \dots, c_T(k)$. Lemma \ref{lem:content-distinction} tells us that there is one when $c_T(1) =\half(\pm(a+p) \pm (b+q))$, and there are two otherwise. So if  $c_T(1) \neq \half(\pm(a+p) \pm (b+q))$ define
	\begin{equation}\label{def:s_0}s_0 T = (s_0T^{(0)} \searrow  T^{(1)}\searrow \ldots\searrow T^{(k)}),\end{equation}
	where $s_0T^{(0)}$ is the unique partition built by moving $T^{(1)}/T^{(0)}$ to its complementary position (see Remark \ref{rk:move-a-box} or Figure \ref{fig:added-boxes-for-2Ds}). Since 
		$c(T^{(1)}/s_0T^{(0)}) = a-p+b-q - 2 c(T^{(1)}/T^{(0)})$,
we have 
	\begin{equation}
	\label{eq:symmetric_contents}
		c_{s_0T}(1) = - c_{T}(1) .
	\end{equation} 


\begin{prop} \label{thm:seminormal-Hecke}
Fix $\lambda \in \cP_k$ and define
$$\cH^\lambda = \mathrm{span}_\CC\{ ~ {v}_T ~|~ T \in \cT_\lambda~\}$$ 
as a vector space with basis indexed by all tableaux from any $\mu \in \cP$ to $\lambda$. 
Define an action of $\cH_k^{\ext}$ by 
\begin{align*}
	w_i \cdot {v}_T &= c_T(i) {v}_T,
	&\text{ for $0 \leq i \leq k$}\\
	t_{s_i} \cdot {v}_T &=  [t_i]_{T,T} {v}_T + [t_i]_{T,s_iT} v_{s_i T},
	& \text{ for $1 \leq i \leq k-1$}\\
	x_{1} \cdot {v}_T &= [x_1]_{T,T} {v}_T +  [x_1]_{T,s_0T} {v}_{s_0T}
\end{align*}
where $[t_i]_{T,s_iT} = 0$ if and only if $c_T(i) = c_T(i+1) \pm 1$, and $ [x_1]_{T,s_0T}=0 $ if and only if  
$c_T(1) = \half(\pm(a+p) \pm (b+q))$. 
Then $\cH^\lambda$ is a simple $\cH^\mathrm{ext}_k$-module with respect to this action if and only if
\begin{enumerate}
\item
$	\lbrack t_i \rbrack_{T,T} = 1/(c_T(i+1) - c_T(i))$,
\item
$\displaystyle	\lbrack x_1 \rbrack_{T,T}   = 
			\frac{(a-p) c_T(1) + c^2_T(1) +   \left(\frac{(a+p)+(b+q)}{2} \right)\left(\frac{(a+p)-(b+q)}{2}\right)}{2 c_T(1) } $,
\item Commutation:
 $$[t_i]_{s_jT,s_is_jT} [t_j]_{T,s_jT}=[t_i]_{T,s_iT} [t_j]_{s_iT,s_js_iT},\qquad \textrm{for } j \neq i \pm 1,$$
 $$[t_i]_{s_0T,s_is_0T} [x_1]_{T,s_0T}=[t_i]_{T,s_iT} [x_1]_{s_iT,s_0s_iT},\qquad \textrm{for } i>1,$$
\item Involutions:
	$$[t_{i}]_{T,s_iT}[t_{i}]_{s_iT,T} = 1 - ([t_{i}]_{T,T})^2,$$
\item Quadratic relation:
	\begin{align*}
	[x_1]_{T,s_0T}[x_1]_{s_0T,T} &= \begin{array}{l}
		 		-\frac{1}{(2 c_T(1))^2}\left(c_T(1) + \frac{(a+p) + (b+q)}{2} \right)\left(c_T(1)- \frac{(a+p) - (b+q)}{2} \right)\\
					\qquad \qquad\cdot\left(c_T(1) - \frac{(a+p) + (b+q)}{2} \right)\left(c_T(1) + \frac{(a+p) - (b+q)}{2}  \right) ,
					\end{array}
			\end{align*}
\item Braid relations:
$$	[t_{i}]_{T,s_iT} 
	[t_{{i+1}}]_{s_iT,s_{i+1}s_iT}
	[t_{i}]_{s_{i+1}s_iT,s_i s_{i+1}s_iT} 
	=
	 [t_{i+1}]_{T,s_{i+1}T} 
	 [t_{i}]_{s_{i+1}T,s_{i}s_{i+1}T} 
	 [t_{i+1}]_{s_{i}s_{i+1}T,s_i s_{i+1}s_iT}, $$
\begin{align*}
[x_1]_{s_1 T,s_0s_1 T}[x_1]_{s_1 s_0 s_1T,s_0s_1 s_0 s_1 T}&[t_1]_{T,s_1T}[t_1]_{s_0 s_1T,s_1s_0 s_1T}  \\&= [x_1]_{T,s_0T} [x_1]_{s_1s_0T,s_0s_1s_0T} [t_1]_{s_0T,s_1s_0T}[t_1]_{s_0s_1s_0T,s_1s_0s_1s_0T}.\end{align*}

\end{enumerate}

\end{prop}

Before we provide a proof of this proposition, we will give a nice example of such a seminormal representation. 
	
\begin{thm}\label{thm:seminormal-Hecke-values}
 Define an action of $\cH^\mathrm{ext}_k$ on $\cH^\lambda$ by
\begin{align*}
	w_i \cdot {v}_T &= c_T(i) {v}_T,&
	\text{ for $0 \leq i \leq k$}\\
	t_{s_i} \cdot {v}_T &=  [t_i]_{T,T} {v}_T + [t_i]_{T,s_iT} v_{s_i T},&
	 \text{ for $1 \leq i \leq k-1$,}\\
	x_{1} \cdot {v}_T &= [x_1]_{T,T} {v}_T +  [x_1]_{T,s_0T} {v}_{s_0T},
\end{align*}
and 
\begin{displaymath}
	\lbrack t_i \rbrack_{T,S}   =  \left\{
		\begin{array}{ll} 
			\sqrt{1 -  \lbrack t_i \rbrack^2_{T,T}} & \textrm{if $S \neq T$,}\\ 
			~ & ~ \\
			1/(c_T(i+1) - c_T(i)) & \textrm{if $S=T$,} \nonumber
		\end{array} \right.
\end{displaymath}
\begin{displaymath}
	\lbrack x_1 \rbrack_{T,S}   =  \left\{
		\begin{array}{ll} 
			 \sqrt{\begin{array}{l}
		 		-\frac{1}{(2 c_T(1))^2}\left(c_T(1) + \frac{(a+p) + (b+q)}{2} \right)\left(c_T(1)- \frac{(a+p) - (b+q)}{2} \right)\\
					\qquad \qquad\cdot\left(c_T(1) - \frac{(a+p) + (b+q)}{2} \right)\left(c_T(1) + \frac{(a+p) - (b+q)}{2}  \right) 
					\end{array}} & \textrm{if $S \neq T$,}\\ 
			~ & ~ \\
			\displaystyle \frac{(a-p) c_T(1) + c^2_T(1) +   \left(\frac{(a+p)+(b+q)}{2} \right)\left(\frac{(a+p)-(b+q)}{2}\right)}{2 c_T(1) } & \textrm{if $S=T$.} \nonumber
		\end{array} \right.
\end{displaymath}
With this action, $\cH^\lambda$ is a simple $\cH^\mathrm{ext}_k$-module. 
\end{thm}
\begin{proof}
The values for $[t_{i}]_{T,T}$ and  $[x_1]_{T,T}$ are pulled directly from Theorem \ref{thm:seminormal-Hecke}, so we need only check criteria 3-6: Commutation, Quadratic relation, and Braid relations. We will verify these using the fact that  $\lbrack x_1 \rbrack_{T,S}$ and $[t_{i}]_{T,S}$ for $S \neq T$ are functions of shifted contents $c_T(j)$. 
\paragraph{\textbf{Commutation:}}  For $j \neq i \pm 1$, 
$c_T(i) = c_{s_jT}(i)$, $c_T(i+1) = c_{s_jT}(i+1)$,
$c_T(j) = c_{s_iT}(j)$, and $c_T(j+1) = c_{s_iT}(j+1)$,
so
$$[t_i]_{s_jT,s_is_jT}=[t_i]_{T,s_iT} \quad \text{and} \quad [t_j]_{T,s_jT} = [t_j]_{s_iT,s_js_iT}.$$
Similarly, for $i>1$, $c_T(i) = c_{s_0T}(i)$ and $c_T(i+1) = c_{s_0T}(i+1)$, so 
	$[t_{i}]_{s_0T,s_is_0T} = [t_{i}]_{T,s_iT},$
	and $c_T(1) = c_{s_iT}(1)$, so 
	$ [x_1]_{T,s_0T} = [x_1]_{s_iT,s_0s_iT}.$
Thus criteria 3 is satisfied. 
\paragraph{\textbf{Quadratic Relation:}} By equation \eqref{eq:symmetric_contents}, 
	$\lbrack x_1 \rbrack_{T,s_0T} = \lbrack x_1 \rbrack_{T,s_0T},$
	so criteria 4 is satisfied.
\paragraph{\textbf{Braid relations:}} For the first braid relation, let $A = c_T(i)$, $B = c_T(i+1)$, and $C =c_T(i+2)$. Either both sides of the equality $$
	[t_{i}]_{T,s_iT} 
	[t_{{i+1}}]_{s_iT,s_{i+1}s_iT}
	[t_{i}]_{s_{i+1}s_iT,s_i s_{i+1}s_iT} 
	=
	 [t_{i+1}]_{T,s_{i+1}T} 
	 [t_{i}]_{s_{i+1}T,s_{i}s_{i+1}T} 
	 [t_{i+1}]_{s_{i}s_{i+1}T,s_i s_{i+1}s_iT} 
$$
are zero, or the six tableaux involved sit in a subgraph of the Bratteli diagram depicted in Figure \ref{fig:sub_Brat1}.
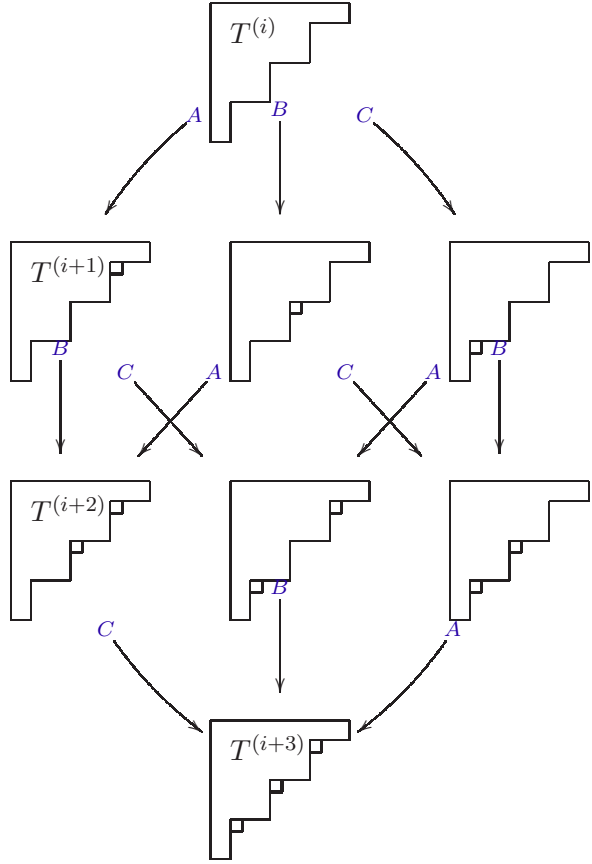
\begin{figure}
\caption{A generic subgraph of the Bratteli diagram in levels $i+1$ through $i+4$.}
\label{fig:sub_Brat1}
$$\def\UNIT{1.5pt}
\xymatrix{
&\setlength{\unitlength}{\UNIT}
\begin{picture}(35, 40)(0,0)
		\put(5,25) {$T^{(i)}$}
		\put(0,0) {\line(1,0){5}} 
		\put(5,10) {\line(1,0){10}} 
		\put(15,20) {\line(1,0){10}} 
		\put(25,30) {\line(1,0){10}} 
		\put(0,35) {\line(1,0){35}} 
		\put(0,0) {\line(0,1){35}} 
		\put(5,0) {\line(0,1){10}} 
		\put(15,10) {\line(0,1){10}} 
		\put(25,20) {\line(0,1){10}} 
		\put(35,30) {\line(0,1){5}} 
		%
		\end{picture} 
		\ar@/_1pc/[dl]|(.3){\phantom{|}\!\! \color{dblue} A} 
		\ar[d]|(.2){\phantom{|}\!\! \color{dblue} B}
		\ar@/^1pc/[dr]|(.3){\phantom{|}\!\! \color{dblue} C}
& \\
\setlength{\unitlength}{\UNIT}
\begin{picture}(35, 40)(-5,0)
		\put(5,25) {$T^{(i+1)}$}
		\put(0,0) {\line(1,0){5}} 
		\put(5,10) {\line(1,0){10}} 
		\put(15,20) {\line(1,0){10}} 
		\put(25,30) {\line(1,0){10}} 
		\put(0,35) {\line(1,0){35}} 
		\put(0,0) {\line(0,1){35}} 
		\put(5,0) {\line(0,1){10}} 
		\put(15,10) {\line(0,1){10}} 
		\put(25,20) {\line(0,1){10}} 
		\put(35,30) {\line(0,1){5}} 
		\thicklines
		\put(25,27) {\line(1,0){3}} 
		\put(28,27) {\line(0,1){3}} 
		\end{picture} 
		\ar[d]|(.2){\phantom{|}\!\! \color{dblue} B}
		\ar@/^0pc/[dr]|(.3){\phantom{|}\!\! \color{dblue} C}&\setlength{\unitlength}{\UNIT}
\begin{picture}(35, 40)(-5,0)
		\put(0,0) {\line(1,0){5}} 
		\put(5,10) {\line(1,0){10}} 
		\put(15,20) {\line(1,0){10}} 
		\put(25,30) {\line(1,0){10}} 
		\put(0,35) {\line(1,0){35}} 
		\put(0,0) {\line(0,1){35}} 
		\put(5,0) {\line(0,1){10}} 
		\put(15,10) {\line(0,1){10}} 
		\put(25,20) {\line(0,1){10}} 
		\put(35,30) {\line(0,1){5}} 
		\thicklines
		%
		\put(15,17) {\line(1,0){3}} 
		\put(18,17) {\line(0,1){3}} 
		%
		\end{picture} 
		\ar@/_0pc/[dl]|(.3){\phantom{|}\!\! \color{dblue} A} 
		\ar@/^0pc/[dr]|(.3){\phantom{|}\!\! \color{dblue} C}
&\setlength{\unitlength}{\UNIT}
\begin{picture}(35, 40)(-5,0)
		\put(0,0) {\line(1,0){5}} 
		\put(5,10) {\line(1,0){10}} 
		\put(15,20) {\line(1,0){10}} 
		\put(25,30) {\line(1,0){10}} 
		\put(0,35) {\line(1,0){35}} 
		\put(0,0) {\line(0,1){35}} 
		\put(5,0) {\line(0,1){10}} 
		\put(15,10) {\line(0,1){10}} 
		\put(25,20) {\line(0,1){10}} 
		\put(35,30) {\line(0,1){5}} 
		\thicklines
		\put(5,7) {\line(1,0){3}} 
		\put(8,7) {\line(0,1){3}} 
		%
		\end{picture} 
		\ar@/_0pc/[dl]|(.3){\phantom{|}\!\! \color{dblue} A} 
		\ar[d]|(.2){\phantom{|}\!\! \color{dblue} B}
\\
\setlength{\unitlength}{\UNIT}
\begin{picture}(35, 40)(-5,0)
		\put(5,25) {$T^{(i+2)}$}
		\put(0,0) {\line(1,0){5}} 
		\put(5,10) {\line(1,0){10}} 
		\put(15,20) {\line(1,0){10}} 
		\put(25,30) {\line(1,0){10}} 
		\put(0,35) {\line(1,0){35}} 
		\put(0,0) {\line(0,1){35}} 
		\put(5,0) {\line(0,1){10}} 
		\put(15,10) {\line(0,1){10}} 
		\put(25,20) {\line(0,1){10}} 
		\put(35,30) {\line(0,1){5}} 
		\thicklines
		%
		\put(15,17) {\line(1,0){3}} 
		\put(18,17) {\line(0,1){3}} 
		\put(25,27) {\line(1,0){3}} 
		\put(28,27) {\line(0,1){3}} 
		\end{picture} 
		\ar@/_1pc/[dr]|(.3){\phantom{|}\!\! \color{dblue} C}
&
\setlength{\unitlength}{\UNIT}
\begin{picture}(35, 40)(-5,0)
		\put(0,0) {\line(1,0){5}} 
		\put(5,10) {\line(1,0){10}} 
		\put(15,20) {\line(1,0){10}} 
		\put(25,30) {\line(1,0){10}} 
		\put(0,35) {\line(1,0){35}} 
		\put(0,0) {\line(0,1){35}} 
		\put(5,0) {\line(0,1){10}} 
		\put(15,10) {\line(0,1){10}} 
		\put(25,20) {\line(0,1){10}} 
		\put(35,30) {\line(0,1){5}} 
		\thicklines
		\put(5,7) {\line(1,0){3}} 
		\put(8,7) {\line(0,1){3}} 
		%
		%
		\put(25,27) {\line(1,0){3}} 
		\put(28,27) {\line(0,1){3}} 
		\end{picture} 
		\ar[d]|(.2){\phantom{|}\!\! \color{dblue} B}
&
\setlength{\unitlength}{\UNIT}
\begin{picture}(35, 40)(-5,0)
		\put(0,0) {\line(1,0){5}} 
		\put(5,10) {\line(1,0){10}} 
		\put(15,20) {\line(1,0){10}} 
		\put(25,30) {\line(1,0){10}} 
		\put(0,35) {\line(1,0){35}} 
		\put(0,0) {\line(0,1){35}} 
		\put(5,0) {\line(0,1){10}} 
		\put(15,10) {\line(0,1){10}} 
		\put(25,20) {\line(0,1){10}} 
		\put(35,30) {\line(0,1){5}} 
		\thicklines
		\put(5,7) {\line(1,0){3}} 
		\put(8,7) {\line(0,1){3}} 
		\put(15,17) {\line(1,0){3}} 
		\put(18,17) {\line(0,1){3}} 
		%
		\end{picture} 
		\ar@/^1pc/[dl]|(.3){\phantom{|}\!\! \color{dblue} A} 
\\
&\setlength{\unitlength}{\UNIT}
\begin{picture}(35, 40)(0,0)
		\put(5,25) {$T^{(i+3)}$}
		\put(0,0) {\line(1,0){5}} 
		\put(5,10) {\line(1,0){10}} 
		\put(15,20) {\line(1,0){10}} 
		\put(25,30) {\line(1,0){10}} 
		\put(0,35) {\line(1,0){35}} 
		\put(0,0) {\line(0,1){35}} 
		\put(5,0) {\line(0,1){10}} 
		\put(15,10) {\line(0,1){10}} 
		\put(25,20) {\line(0,1){10}} 
		\put(35,30) {\line(0,1){5}} 
		\thicklines
		\put(5,7) {\line(1,0){3}} 
		\put(8,7) {\line(0,1){3}} 
		\put(15,17) {\line(1,0){3}} 
		\put(18,17) {\line(0,1){3}} 
		\put(25,27) {\line(1,0){3}} 
		\put(28,27) {\line(0,1){3}} 
		\end{picture} 
&
}
$$
\end{figure}
  This encodes the fact that for whichever of these $S$ exist, their shifted contents follow the pattern in \ref{tab:length-3-contents}, and one can use these values to check that the first braid relation is satisfied.
\begin{equation}\label{tab:length-3-contents}
\begin{array}{r|llllll}
	S\to		& T 		& s_iT 	& s_{i+1}T& s_i s_{i+1} T	& s_{i+1}s_i  T & s_i s_{i+1} s_i  T \\
	\hline
	c_S(i)	&A		&B		&A		&C			&B			&C		\\
	c_S(i+1)	&B		&A		&C		&A			&C			&B		\\
	c_S(i+2)	&C		&C		&B		&B			&A			&A	\\
\end{array}
\end{equation}
 

For the second braid relation,  let $A = c_T(1)$ and $B = c_T(2)$. So either both sides of the equality 
\begin{align*}[x_1]_{s_1 T,s_0s_1 T}[x_1]_{s_1 s_0 s_1T,s_0s_1 s_0 s_1 T}&[t_1]_{T,s_1T}[t_1]_{s_0 s_1T,s_1s_0 s_1T}  \\&= [x_1]_{T,s_0T} [x_1]_{s_1s_0T,s_0s_1s_0T} [t_1]_{s_0T,s_1s_0T}[t_1]_{s_0s_1s_0T,s_1s_0s_1s_0T}\end{align*}
are zero, or the eight tableaux involved sit in a subgraph of the Bratteli diagram depicted in Figure \ref{fig:sub_Brat12}.
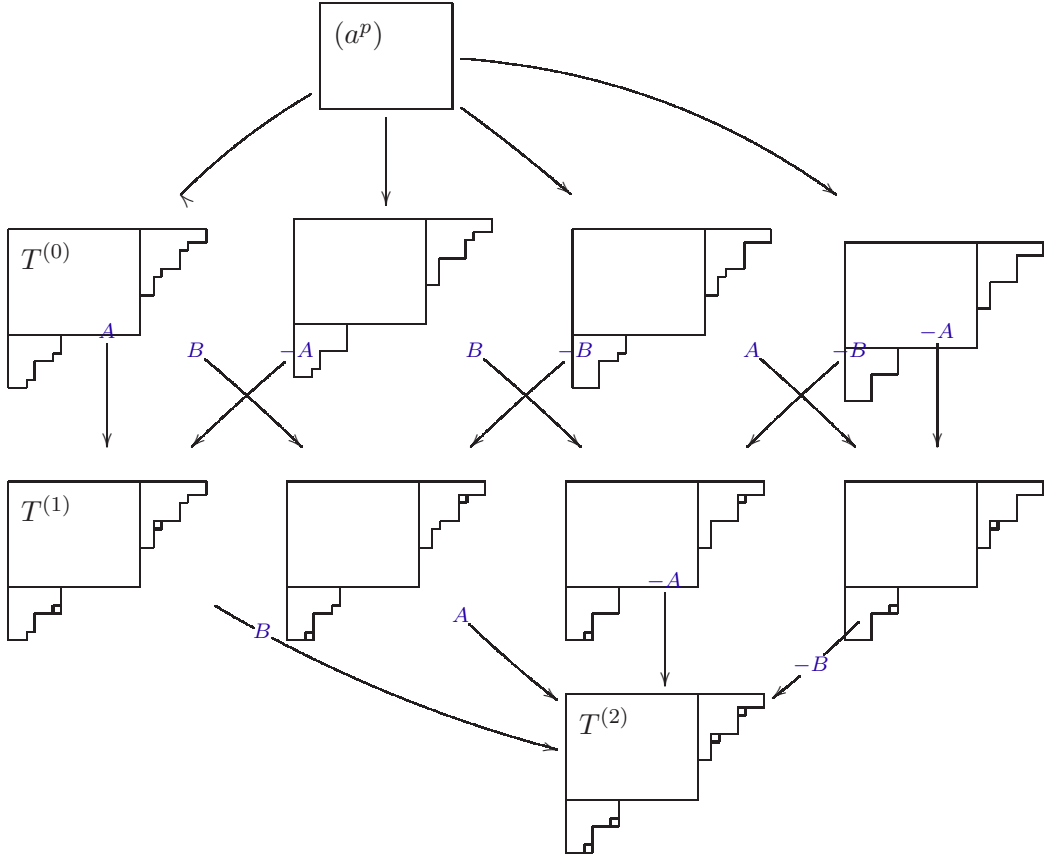
\begin{figure}
\caption{A generic subgraph of the Bratteli diagram in levels $0$ through $3$.}
\label{fig:sub_Brat12}
{\def\UNIT{1pt}$$\xymatrix{
&\setlength{\unitlength}{\UNIT}
	\begin{picture}(50, 40)(0,0)
		\multiput(0,0)(50,0){2} {\line(0,1){40}}
		\multiput(0,0)(0,40){2} {\line(1,0){50}}
		\put(5,25){$(a^p)$}
		\end{picture}
		\ar@/_.8pc/[dl] \ar[d] \ar@/^.3pc/[dr] \ar@/^2pc/[drr] 
& &\\
\setlength{\unitlength}{\UNIT}
\begin{picture}(75, 60)(0,-10)
		\put(5,25) {$T^{(0)}$}
		\multiput(0,0)(50,0){2} {\line(0,1){40}}
		\multiput(0,0)(0,40){2} {\line(1,0){50}}
		\put(50,40) {\line(1,0){25}} 
		\put(75,40) {\line(0,-1){5}}
		\put(75,35) {\line(-1,0){7}}
		\put(68,35) {\line(0,-1){3}}
		\put(68,32) {\line(-1,0){3}}
		\put(65,32) {\line(0,-1){7}}
		\put(65,25) {\line(-1,-0){7}}
		\put(58,25) {\line(0,-1){3}}
		\put(58,22) {\line(-1,0){3}}
		\put(55,22) {\line(0,-1){7}}
		\put(55,15) {\line(-1,0){5}}
		\put(0,0) {\line(0,-1){20}}
		\put(0,-20) {\line(1,0){7}}
		\put(7,-20) {\line(0,1){3}}
		\put(7,-17) {\line(1,0){3}}
		\put(10,-17) {\line(0,1){7}}
		\put(10,-10) {\line(1,0){7}}
		\put(17,-10) {\line(0,1){3}}
		\put(17,-7) {\line(1,0){3}}
		\put(20,-7) {\line(0,1){7}}
		\end{picture} 
		\ar[d]|(.2){\phantom{|}\!\! \color{dblue} A}	
		\ar@/^.3pc/[dr]|(.3){\phantom{|}\!\! \color{dblue} B}
&\setlength{\unitlength}{\UNIT}
\begin{picture}(70, 52)(0,-10)
		\multiput(0,0)(50,0){2} {\line(0,1){40}}
		\multiput(0,0)(0,40){2} {\line(1,0){50}}
		\put(50,40) {\line(1,0){25}} 
		\put(75,40) {\line(0,-1){5}}
		\put(75,35) {\line(-1,0){7}}
		\put(68,35) {\line(0,-1){3}}
		\put(68,32) {\line(-1,0){3}}
		\put(65,32) {\line(0,-1){7}}
		\put(65,25) {\line(-1,-0){10}}
		\put(55,25) {\line(0,-1){10}}
		\put(55,15) {\line(-1,0){5}}
		\put(0,0) {\line(0,-1){20}}
		\put(0,-20) {\line(1,0){7}}
		\put(7,-20) {\line(0,1){3}}
		\put(7,-17) {\line(1,0){3}}
		\put(10,-17) {\line(0,1){7}}
		\put(10,-10) {\line(1,0){10}}
		\put(20,-10) {\line(0,1){10}}
		\end{picture} 
		\ar@/_.3pc/[dl]|(.3){\phantom{|}\!\! \color{dblue} -A} 
		\ar@/^.3pc/[dr]|(.3){\phantom{|}\!\! \color{dblue} B}
&\setlength{\unitlength}{\UNIT}
\begin{picture}(70, 60)(0,-10)
		\multiput(0,0)(50,0){2} {\line(0,1){40}}
		\multiput(0,0)(0,40){2} {\line(1,0){50}}
		\put(50,40) {\line(1,0){25}} 
		\put(75,40) {\line(0,-1){5}}
		\put(75,35) {\line(-1,0){10}}
		\put(65,35) {\line(0,-1){10}}
		\put(65,25) {\line(-1,-0){7}}
		\put(58,25) {\line(0,-1){3}}
		\put(58,22) {\line(-1,0){3}}
		\put(55,22) {\line(0,-1){7}}
		\put(55,15) {\line(-1,0){5}}
		\put(0,0) {\line(0,-1){20}}
		\put(0,-20) {\line(1,0){10}}
		\put(10,-20) {\line(0,1){10}}
		\put(10,-10) {\line(1,0){7}}
		\put(17,-10) {\line(0,1){3}}
		\put(17,-7) {\line(1,0){3}}
		\put(20,-7) {\line(0,1){7}}
		\end{picture} 
		\ar@/_.3pc/[dl]|(.3){\phantom{|}\!\! \color{dblue} -B} 
		\ar@/^.3pc/[dr]|(.3){\phantom{|}\!\! \color{dblue} A}
&\setlength{\unitlength}{\UNIT}
\begin{picture}(70, 70)(0,-10)
		\multiput(0,0)(50,0){2} {\line(0,1){40}}
		\multiput(0,0)(0,40){2} {\line(1,0){50}}
		\put(50,40) {\line(1,0){25}} 
		\put(75,40) {\line(0,-1){5}}
		\put(75,35) {\line(-1,0){10}}
		\put(65,35) {\line(0,-1){10}}
		\put(65,25) {\line(-1,-0){10}}
		\put(55,25) {\line(0,-1){10}}
		\put(55,15) {\line(-1,0){5}}
		\put(0,0) {\line(0,-1){20}}
		\put(0,-20) {\line(1,0){10}}
		\put(10,-20) {\line(0,1){10}}
		\put(10,-10) {\line(1,0){10}}
		\put(20,-10) {\line(0,1){10}}
		\end{picture} 
		\ar@/_.3pc/[dl]|(.3){\phantom{|}\!\! \color{dblue} -B} 
		\ar[d]|(.2){\phantom{|}\!\! \color{dblue} -A}
\\\setlength{\unitlength}{\UNIT}
\begin{picture}(75, 60)(0,-10)
		\put(5,25) {$T^{(1)}$}
		\multiput(0,0)(50,0){2} {\line(0,1){40}}
		\multiput(0,0)(0,40){2} {\line(1,0){50}}
		\put(50,40) {\line(1,0){25}} 
		\put(75,40) {\line(0,-1){5}}
		\put(75,35) {\line(-1,0){7}}
		\put(68,35) {\line(0,-1){3}}
		\put(68,32) {\line(-1,0){3}}
		\put(65,32) {\line(0,-1){7}}
		\put(65,25) {\line(-1,-0){10}}
		\put(55,25) {\line(0,-1){10}}
		{\thicklines
		\put(58,25) {\line(0,-1){3}}
		\put(58,22) {\line(-1,0){3}}
		}
		\put(55,15) {\line(-1,0){5}}
		\put(0,0) {\line(0,-1){20}}
		\put(0,-20) {\line(1,0){7}}
		\put(7,-20) {\line(0,1){3}}
		\put(7,-17) {\line(1,0){3}}
		\put(10,-17) {\line(0,1){7}}
		\put(10,-10) {\line(1,0){10}}
		\put(20,-10) {\line(0,1){10}}
		{\thicklines
		\put(17,-10) {\line(0,1){3}}
		\put(17,-7) {\line(1,0){3}}
		}
		\end{picture} 
		\ar@/_1.1pc/[drr]|(.3){\phantom{|}\!\! \color{dblue} B}
&\setlength{\unitlength}{\UNIT}
\begin{picture}(75, 60)(0,-10)
		\multiput(0,0)(50,0){2} {\line(0,1){40}}
		\multiput(0,0)(0,40){2} {\line(1,0){50}}
		\put(50,40) {\line(1,0){25}} 
		\put(75,40) {\line(0,-1){5}}
		\put(75,35) {\line(-1,0){10}}
		\put(65,35) {\line(0,-1){10}}
		{\thicklines
		\put(68,35) {\line(0,-1){3}}
		\put(68,32) {\line(-1,0){3}}
		}
		\put(65,25) {\line(-1,-0){7}}
		\put(58,25) {\line(0,-1){3}}
		\put(58,22) {\line(-1,0){3}}
		\put(55,22) {\line(0,-1){7}}
		\put(55,15) {\line(-1,0){5}}
		\put(0,0) {\line(0,-1){20}}
		\put(0,-20) {\line(1,0){10}}
		\put(10,-20) {\line(0,1){10}}
		{\thicklines
		\put(7,-20) {\line(0,1){3}}
		\put(7,-17) {\line(1,0){3}}
		}
		\put(10,-10) {\line(1,0){7}}
		\put(17,-10) {\line(0,1){3}}
		\put(17,-7) {\line(1,0){3}}
		\put(20,-7) {\line(0,1){7}}
		\end{picture} 
		\ar@/_.5pc/[dr]|(.3){\phantom{|}\!\! \color{dblue} A}
&\setlength{\unitlength}{\UNIT}
\begin{picture}(75, 60)(0,-10)
		\multiput(0,0)(50,0){2} {\line(0,1){40}}
		\multiput(0,0)(0,40){2} {\line(1,0){50}}
		\put(50,40) {\line(1,0){25}} 
		\put(75,40) {\line(0,-1){5}}
		\put(75,35) {\line(-1,0){10}}
		\put(65,35) {\line(0,-1){10}}
		{\thicklines
		\put(68,35) {\line(0,-1){3}}
		\put(68,32) {\line(-1,0){3}}
		}
		\put(65,25) {\line(-1,-0){10}}
		\put(55,25) {\line(0,-1){10}}
		\put(55,15) {\line(-1,0){5}}
		\put(0,0) {\line(0,-1){20}}
		\put(0,-20) {\line(1,0){10}}
		\put(10,-20) {\line(0,1){10}}
		{\thicklines
		\put(7,-20) {\line(0,1){3}}
		\put(7,-17) {\line(1,0){3}}
		}
		\put(10,-10) {\line(1,0){10}}
		\put(20,-10) {\line(0,1){10}}
		\end{picture} 		
		\ar[d]|(.2){\phantom{|}\!\! \color{dblue} -A}
&\setlength{\unitlength}{\UNIT}
\begin{picture}(70, 60)(0,-10)
		\multiput(0,0)(50,0){2} {\line(0,1){40}}
		\multiput(0,0)(0,40){2} {\line(1,0){50}}
		\put(50,40) {\line(1,0){25}} 
		\put(75,40) {\line(0,-1){5}}
		\put(75,35) {\line(-1,0){10}}
		\put(65,35) {\line(0,-1){10}}
		\put(65,25) {\line(-1,-0){10}}
		\put(55,25) {\line(0,-1){10}}
		{\thicklines
		\put(58,25) {\line(0,-1){3}}
		\put(58,22) {\line(-1,0){3}}
		}
		\put(55,15) {\line(-1,0){5}}
		\put(0,0) {\line(0,-1){20}}
		\put(0,-20) {\line(1,0){10}}
		\put(10,-20) {\line(0,1){10}}
		{\thicklines
		\put(17,-10) {\line(0,1){3}}
		\put(17,-7) {\line(1,0){3}}}
		\put(10,-10) {\line(1,0){10}}
		\put(20,-10) {\line(0,1){10}}
		\end{picture} 
		\ar@/^.5pc/[dl]|{\phantom{|}\!\! \color{dblue} -B}
\\
&&\setlength{\unitlength}{\UNIT}\begin{picture}(75, 60)(0,-20)
		\put(5,25) {$T^{(2)}$}
		\multiput(0,0)(50,0){2} {\line(0,1){40}}
		\multiput(0,0)(0,40){2} {\line(1,0){50}}
		\put(50,40) {\line(1,0){25}} 
		\put(75,40) {\line(0,-1){5}}
		\put(75,35) {\line(-1,0){10}}
		\put(65,35) {\line(0,-1){10}}
		{\thicklines
		\put(68,35) {\line(0,-1){3}}
		\put(68,32) {\line(-1,0){3}}}
		\put(65,25) {\line(-1,-0){10}}
		\put(55,25) {\line(0,-1){10}}
		{\thicklines
		\put(58,25) {\line(0,-1){3}}
		\put(58,22) {\line(-1,0){3}}}
		\put(55,15) {\line(-1,0){5}}
		\put(0,0) {\line(0,-1){20}}
		\put(0,-20) {\line(1,0){10}}
		\put(10,-20) {\line(0,1){10}}
		{\thicklines
		\put(7,-20) {\line(0,1){3}}
		\put(7,-17) {\line(1,0){3}}}
		\put(10,-10) {\line(1,0){10}}
		\put(20,-10) {\line(0,1){10}}
		{\thicklines
		\put(17,-10) {\line(0,1){3}}
		\put(17,-7) {\line(1,0){3}}}
		\end{picture}
&
 }
$$}\end{figure}
  This encodes the fact that for whichever of these $S$ exist, their shifted contents follow the pattern in \ref{tab:length-4-contents}, and one can use these values to check that the first braid relation is satisfied
\begin{equation}\label{tab:length-4-contents}
\begin{array}{rrrrrrrrr}
	S\to		& T 		& s_0T 	& s_1T& s_0s_1 T& s_1s_0  T & s_0s_1s_0  T& s_1s_0s_1  T& s_0s_1s_0s_1  T \\
	\hline
	c_S(1)	&A		&-A	&B	    &-B	     &B		&-B		&A			&-A	\\
	c_S(2)	&B		&B		&A	    &A	     &-A		&-A		&-B		&-B	
\end{array}
\end{equation}
Thus criteria 6 is satisfied, concluding the proof of Theorem \ref{thm:seminormal-Hecke-values}.
\end{proof}

\begin{proof}[Proof of Proposition \ref{thm:seminormal-Hecke}]
We prove Proposition  \ref{thm:seminormal-Hecke} in two parts. In Part 1, we check that the relations in the presentation of $\cH_k^\ext$ given in Theorem \ref{thm:hecke_ext_presentation-short}  hold, showing that $\cH^\lambda$ is a $\cH_k^\ext$-module. In Part 2, we verify that $\cH^\lambda$ is simple.

\paragraph{\textbf{Part 1: $\cH^\lambda$ is a $\cH_k^\ext$-module.}}

By \eqref{eq:comm-ww}, the elements $w_0, w_1, \dots, w_k$ generate a commutative subalgebra of $\cH^\mathrm{ext}_k$, so we begin by fixing the diagonal action as stated above, 
\begin{align*}
	w_0  {v}_T &= \left( abq + 2\sum_{B \in \cB_\lambda} \left( c(B) - \half(a-p+b-q)\right) \right) v_T,\\
	w_i  {v}_T &= c_T(i) {v}_T,\quad \text{ for $1 \leq i \leq k$}.\\
\end{align*}

Now write 
	$$t_{s_i} v_T = \sum_{S \in \cT_\lambda} [t_{i}]_{T,S} v_S 
	\quad \text{and} \quad 
	x_1 v_T = \sum_{S \in \cT_\lambda} [x_1]_{T,S} v_S,$$ 
	where $\cT_\lambda$ is the set of tableaux \eqref{eq:indexing-tableau} and $[t_{i}]_{T,S}, [x_1]_{T,S} \in \CC$. 
	
\begin{enumerate}
\item[]\textbf{Claim 1:} Relations $t_{s_i}^2=1$, \eqref{eq:comm-tw}, and \eqref{eq:twist-tw} are satisfied if and only if
$$t_{s_i} v_T =  [t_{i}]_{T,T} v_T + [t_{i}]_{T,s_iT} v_{s_iT}, \quad \text{for $i = 1, \dots, k-1$,}$$
$$[t_{i}]_{T,T} = \frac{1}{c_T(i+1) - c_T(i)} \quad \text{ and } \quad  [t_{i}]_{T,S}[t_{i}]_{S,T} = 1 - ([t_{i}]_{T,T})^2.$$

\noindent\emph{Proof:} The first commutation relation \eqref{eq:comm-tw}, $t_{s_i} w_j = w_j t_{s_i}$ for $j \ne i, i+1$, implies 
\begin{align*}
	t_{s_i} w_j  v_T 	&=  \sum_{S \in \cT_\lambda }  c_T(j) [t_{i}]_{T,S} v_S \\
	= w_j t_{s_i}   v_T 	&= \sum_{S \in \cT_\lambda}   c_S(j) [t_{i}]_{T,S} v_S.
\end{align*}
So for each $S$, either 
\begin{equation}\label{eq:aligned_contents_1}
[t_{i}]_{T,S} = 0 \qquad \text{ or } \qquad  c_T(j) =  c_S(j) \quad \text{for all $j \ne i, i+1$.}
\end{equation}
The first twisting relation \eqref{eq:twist-tw}, together with relation $t_{s_i}^2=1$, require $$t_{s_i} w_i - w_{i+1} t_{s_i} = - 1 = w_it_{s_i}  - t_{s_i} w_{i+1},$$ i.e., 
\begin{align*}
	(t_{s_i} w_i - w_{i+1} t_{s_i})  v_T 	&= \sum_{S\in \cT_\lambda}  (c_T(i) -  c_S(i+1)) [t_{i}]_{T,S} v_S\\
							&=   - v_T\\
	=(w_i t_{s_i}  - t_{s_i} w_{i+1})  v_T 	&= \sum_{S \in \cT_\lambda}  (c_S(i) -  c_T(i+1)) [t_{i}]_{T,S} v_S.
\end{align*}
So 
		\begin{equation}
		[t_{i}]_{T,T} = \frac{1}{c_T(i+1) - c_T(i)}
		\end{equation} 
and for $S \neq T$, either 
	\begin{equation}\label{eq:aligned_contents_2}
	c_S(i+1) = c_T(i) \quad \text{ and } \quad c_S(i) = c_T(i+1)
		\qquad \text{ or } \qquad [t_{i}]_{T,S}=0.
		\end{equation}
By Lemma \ref{lem:content_list}, equations \eqref{eq:aligned_contents_1} and \eqref{eq:aligned_contents_2} tell us	
$$t_{s_i} v_T =  [t_{i}]_{T,T} v_T + [t_{i}]_{T,s_iT} v_{s_iT}, \quad \text{for $i = 1, \dots, k-1$,}$$
where $[t_{i}]_{T,s_iT}=0$ if $c_T(i) = c_T(i+1) \pm 1$. Finally, the involution relation $t_{s_i}^2=1$ implies 
	$$[t_{i}]_{T,T} = - [t_{i}]_{s_iT,s_iT} \quad \text{ and }
	[t_{i}]_{T,S}[t_{i}]_{S,T} = 1 - ([t_{i}]_{T,T})^2.$$
The first is implied by $[t_{i}]_{T,T} = \frac{1}{c_T(i+1) - c_T(i)}$, but the second places a new condition on coefficients. 
This concludes the proof of Claim 1.


\item[]\textbf{Claim 2:} Relation \eqref{eq:comm-xw} is satisfied if and only if 
$$
x_1v_T = [x_1]_{T,T} v_T +  [x_1]_{T,s_0T} v_{s_0T}, \quad \text{where $[x_1]_{T,s_0T} = 0$ if $c_T(1) = \pm(a+p) \pm(b+q)$.}$$
Furthermore, \eqref{eq:quadratic}, \eqref{eq:twist-xw0}, and \eqref{eq:xw} are additionally satisfied if and only if
$$\lbrack x_1 \rbrack_{T,T}   = 
			\frac{(a-p) c_T(1) + c^2_T(1) +   \left(\frac{(a+p)+(b+q)}{2} \right)\left(\frac{(a+p)-(b+q)}{2}\right)}{2 c_T(1) }$$
			and 
			\begin{align*}
	[x_1]_{T,s_0T}[x_1]_{s_0T,T} &= \begin{array}{l}
		 		-\frac{1}{(2 c_T(1))^2}\left(c_T(1) + \frac{(a+p) + (b+q)}{2} \right)\left(c_T(1)- \frac{(a+p) - (b+q)}{2} \right)\\
					\qquad \qquad\cdot\left(c_T(1) - \frac{(a+p) + (b+q)}{2} \right)\left(c_T(1) + \frac{(a+p) - (b+q)}{2}  \right) .
					\end{array}
			\end{align*}

\noindent\emph{Proof:} 
The relation $x_1 w_i = w_i x_1$ for $i > 1$ implies
	\begin{align*}x_1 w_i v_T &= \sum_{S \in T_\lambda} c_T(i) [x_1]_{T,S} v_S\\
	=w_ix_1 v_T &= \sum_{S \in T_\lambda} c_S(i) [x_1]_{T,S} v_S.
	\end{align*}
So by Lemmas \ref{lem:one_or_two} and \ref{lem:content_list}, 
\begin{equation}\label{eq:aligned_contents_x}
x_1v_T = [x_1]_{T,T} v_T +  [x_1]_{T,s_0T} v_{s_0T}, 
	\end{equation}
	where $[x_1]_{T,s_0T} = 0$ if $c_T(1) = \pm(a+p) \pm(b+q)$.

Now let $K= \left(\frac{a+p + b+ q}{2}\right) \left(\frac{a+p - (b+q)}{2}\right)$, so the third twisting relation \eqref{eq:xw},
$$	x_1 w_1 =  -w_1 x_1+ (a-p)w_1 + w_1^2 + K,$$
says
	\begin{align*}
		(x_1 w_1 + w_1x_1)v_T &= (c_T(1) + c_T(1)) [x_1]_{T,T} v_T +   (c_T(1) + c_{s_0T}(1))[x_1]_{T,s_0T} v_{s_0T} \\
					&=  2 c_T(1)[x_1]_{T,T} v_T\\
		=\left((a-p)w_1 + w_1^2 + K\right)v_T	&= \left((a-p)c_T(1) + (c_T(1))^2 + K\right) v_T.
	\end{align*}
So
	\begin{equation}\label{eq:xT,T}
	[x_1]_{T,T} = \frac{\left((a-p)c_T(1) + (c_T(1))^2 + K\right)}{ 2 c_T(1)}.
	\end{equation}
If $S= s_0T$ exists, then the quadratic relation \eqref{eq:quadratic} implies
	\begin{align*} 
		x_1^2 v_T &= ([x_1]_{T,T}^2 + [x_1]_{T,S}[x_1]_{S,T}) v_T \\
			&\quad + ([x_1]_{T,T}[x_1]_{T,S} + [x_1]_{T,S} [x_1]_{S,S}) v_S\\
		=(a-p)x_1 + ap 
			&= ((a-p)[x_1]_{T,T} + ap)v_T + (a-p)[x_1]_{T,S} v_S.
	\end{align*}
We could conclude $([x_1]_{T,T}[x_1]_{T,S} + [x_1]_{T,S} [x_1]_{S,S}) = (a-p)[x_1]_{T,S} $ from \eqref{eq:xT,T}, so this simply tells us that 
	\begin{align*}
	 [x_1]_{T,S}[x_1]_{S,T} &= -[x_1]_{T,T}^2  + (a-p)[x_1]_{T,T} + ap\\
		&= -\left( \frac{\left((a-p)c_T(1) + c^2_T(1) + K\right)}{ 2 c_T(1)}\right)^2 \\
		& \qquad \qquad + (a-p) \left(\frac{\left((a-p)c_T(1) + c^2_T(1) + K\right)}{ 2 c_T(1)}\right) + ap\\
		&=-\frac{1}{4c^2_T(1)} \left(c_T(1) + \frac{(a+p) + (b+q)}{2} \right)\left(c_T(1)- \frac{(a+p) - (b+q)}{2} \right)\\
			& \qquad \qquad \qquad \qquad \left(c_T(1) - \frac{(a+p) + (b+q)}{2} \right)\left(c_T(1) + \frac{(a+p) - (b+q)}{2}\right).
	\end{align*}

Finally, the second twisting relation \eqref{eq:twist-xw0} implies 
	\begin{align*}
		&x_1(w_0 + w_1)v_T = (c_T(0) + c_T(1))[x_1]_{T,T}v_T + (c_T(0) + c_T(1))[x_1]_{T,s_0T} v_{s_0T}\\
		&= (w_0 + w_1)x_1v_T =  (c_T(0) + c_T(1))[x_1]_{T,T}v_T + (c_{s_0T}(0) + c_{s_0T}(1))[x_1]_{T,s_0T} v_{s_0T}.
	\end{align*}
So we require  $$[x_1]_{T,s_0T}=[x_1]_{T,s_0T}=0 \quad \text{ or } \quad  c_T(0) + c_T(1) = c_{s_0T}(0) + c_{s_0T}(1).$$
Recall from \eqref{eq:symmetric_contents} that if $v_{s_0T}$ exists, then $c_{s_0T}(1) = - c_{T}(1) $. So this requirement is equivalent to
	$$ [x_1]_{T,s_0T}=0 \quad \text{ or } \quad c_T(1) = \half( c_{s_0T}(0) - c_T(0) ),$$
and is therefore a consequence of the construction in Lemmas \ref{lem:move-a-box} and \ref{lem:z_0_on_lambda}. 
This concludes the proof of Claim 2.\\

\item[]\textbf{Claim 3:} The second relation in \eqref{eq:S_braid} and relation \eqref{eq:comm-xt} are satisfied if and only if
 $$[t_i]_{s_jT,s_is_jT} [t_j]_{T,s_jT}=[t_i]_{T,s_iT} [t_j]_{s_iT,s_js_iT},\qquad \textrm{for } j \neq i \pm 1,$$
 and
 $$[t_i]_{s_0T,s_is_0T} [x_1]_{T,s_0T}=[t_i]_{T,s_iT} [x_1]_{s_iT,s_0s_iT},\qquad \textrm{for } i>1,$$
respectively. 

\noindent \emph{Proof:} 
For $j \neq i \pm 1$, the second relation in \eqref{eq:S_braid} implies
	\begin{align*}
		t_{s_i}t_{s_j} v_T&= [t_i]_{T,T}[t_j]_{T,T}v_{T}
					+[t_i]_{T,s_iT}[t_j]_{T,T}v_{s_iT}\\
					&\qquad
					+[t_i]_{s_jT,s_jT}[t_j]_{T,s_jT}v_{s_jT}
					+[t_i]_{s_jT,s_is_jT}[t_j]_{T,s_jT}v_{s_is_jT}\\
		=t_{s_j}t_{s_i} v_T&= [t_i]_{T,T}[t_j]_{T,T}v_{T}
					+[t_i]_{T,s_iT}[t_j]_{s_iT,s_iT}v_{s_iT}\\
					&\qquad
					+[t_i]_{T,T}[t_j]_{T,s_iT}v_{s_jT}
					+[t_i]_{T,s_iT}[t_j]_{s_iT,s_js_iT}v_{s_js_iT}.
	\end{align*}
If $s_iT$ and $s_jT$ exist, we already know $[t_j]_{T,T} = [t_j]_{s_iT,s_iT}$ and $[t_i]_{s_jT,s_jT}=[t_i]_{T,T}$ because 
	$c_T(j) = c_{s_iT}(j)$ and $c_T(i) = c_{s_jT}(i)$ for $j \neq i \pm 1$. However, since $s_is_jT = s_js_iT$, we gain the requirement 
	$$[t_i]_{s_jT,s_is_jT}[t_j]_{T,s_jT} = [t_i]_{T,s_iT}[t_j]_{s_iT,s_js_iT}.$$

Similarly, for $i>1$, relation \eqref{eq:comm-xt}  implies
	\begin{align*}
		t_{s_i} x_1v_T &= [t_i]_{T,T} [x_1]_{T,T} v_{T} + [t_i]_{s_0T,s_0T} [x_1]_{T,s_0T} v_{s_0T} \\
					&\qquad + [t_i]_{T,s_iT} [x_1]_{T,T} v_{s_iT} + [t_i]_{s_0T,s_is_0T} [x_1]_{T,s_0T} v_{s_is_0T}\\
		= x_1 t_{s_i} v_T 
					&= [t_i]_{T,T} [x_1]_{T,T} v_{T} + [t_i]_{T,T} [x_1]_{T,s_0T} v_{s_0T} \\
					&\qquad + [t_i]_{T,s_iT} [x_1]_{s_iT,s_iT} v_{s_iT} + [t_i]_{T,s_iT} [x_1]_{s_iT,s_0s_iT} v_{s_is_0T}
	\end{align*}
since $s_0s_iT = s_is_0T$ for $i>1$. If $s_0T$ and $s_iT$ exist, we already require that
	$$[t_i]_{s_0T,s_0T} =  [t_i]_{T,T}  \qquad \text{and}\qquad  [x_1]_{T,T} = [x_1]_{s_iT,s_iT},$$
	 since  $c_T(i) = c_{s_0T}(i)$, $c_T(i+1) = c_{s_0T}(i+1)$, and $c_T(1) = c_{s_iT}(1)$.  However, given 
$s_0T$,  $s_iT,$ and $s_0s_iT$ exist, we gain the requirement
	\begin{equation}
	 [t_i]_{s_0T,s_is_0T} [x_1]_{T,s_0T}=[t_i]_{T,s_iT} [x_1]_{s_iT,s_0s_iT},
	\end{equation}
concluding the proof of Claim 3. \\

\item[]
	\textbf{Claim 4:} The braid relation \eqref{eq:S_braid}
is satisfied if and only if 
	\begin{align*}	
	[t_{i}]_{T,s_iT} &
	[t_{{i+1}}]_{s_iT,s_{i+1}s_iT}
	[t_{i}]_{s_{i+1}s_iT,s_i s_{i+1}s_iT} \\
	=&
	 [t_{i+1}]_{T,s_{i+1}T} 
	 [t_{i}]_{s_{i+1}T,s_{i}s_{i+1}T} 
	 [t_{i+1}]_{s_{i}s_{i+1}T,s_i s_{i+1}s_iT}. \end{align*}

\emph{Proof:}
If $v_S$ exists for $S= s_iT$, $s_{i+1}T$, $s_i s_{i+1} T$, $s_{i+1}s_i  T$, $s_i s_{i+1} s_i  T$, then
\begin{align*}
	t_{s_i}t_{s_{i+1}}t_{s_i} v_T  	
		&= 
			\left(
				[t_{i}]^2_{T,T}
				[t_{{i+1}}]_{T,T} 
				+ 
				[t_{i}]_{T,s_iT}
				[t_{{i+1}}]_{s_iT,s_iT}
				[t_{i}]_{s_{i}T,T} 
				\right)v_T\\
		& \qquad + \left(
				[t_{i}]_{T,T}
				[t_{{i+1}}]_{T,T} 
				[t_{i}]_{T,s_iT} 
				+ 
				[t_{i}]_{T,s_iT}
				[t_{{i+1}}]_{s_iT,s_iT} 
				[t_{i}]_{s_{i}T,s_{i}T} 
				\right) v_{s_iT} \\
		& \qquad	+ \left(
				[t_{i}]_{T,T}
				[t_{{i+1}}]_{T,s_{i+1}T} 
				[t_{i}]_{s_{i+1}T,s_{i+1}T} 
				\right)v_{s_{i+1}T} \\
		& \qquad	+ \left(
				[t_{i}]_{T,T}
				[t_{{i+1}}]_{T,s_{i+1}T}
				[t_{i}]_{s_{i+1}T,s_is_{i+1}T} 
				\right)v_{s_is_{i+1}T}\\
		&\qquad 	 +\left(
				[t_{i}]_{T,s_iT} 
				[t_{{i+1}}]_{s_iT,s_{i+1}s_iT} 
				[t_{i}]_{s_{i+1}s_iT,s_{i+1}s_iT} 
				\right)v_{s_{i+1}s_iT}  \\
		& \qquad	+ \left(
				[t_{i}]_{T,s_iT} 
				[t_{{i+1}}]_{s_iT,s_{i+1}s_iT}
				[t_{i}]_{s_{i+1}s_iT,s_i s_{i+1}s_iT} 
				\right)v_{s_i s_{i+1}s_iT}
\end{align*}
because $s_is_iT = T$. 
Similarly, 
\begin{align*}
	t_{s_{i+1}}t_{s_{i}}t_{s_{i+1}} v_T 
		&= 
			\left([t_{i+1}]^2_{T,T}[t_{{i}}]_{T,T} + [t_{i+1}]_{T,s_{i+1}T}[t_{{i}}]_{s_{i+1}T,s_{i+1}T}[t_{i+1}]_{s_{i+1}T,T} \right)v_T\\
		& \qquad + \big([t_{i+1}]_{T,T}[t_{{i}}]_{T,T} [t_{i+1}]_{T,s_{i+1}T} \\
		&\qquad \qquad + [t_{i+1}]_{T,s_{i+1}T}[t_{{i}}]_{s_{i+1}T,s_{i+1}T} [t_{i+1}]_{s_{i+1}T,s_{i+1}T} \big) v_{s_{i+1}T} \\
		& \qquad	+ [t_{i+1}]_{T,T}[t_{{i}}]_{T,s_{i}T} [t_{i+1}]_{s_{i}T,s_{i}T} v_{s_{i}T} \\
		& \qquad	+ [t_{i+1}]_{T,T}[t_{{i}}]_{T,s_{i}T}[t_{i+1}]_{s_{i}T,s_{i+1}s_{i}T} v_{s_{i+1}s_{i}T}\\
		&\qquad 	 +[t_{i+1}]_{T,s_{i+1}T} [t_{{i}}]_{s_{i+1}T,s_{i}s_{i+1}T} [t_{i+1}]_{s_{i}s_{i+1}T,s_{i}s_{i+1}T} v_{s_{i}s_{i+1}T}  \\
		& \qquad	+ [t_{i+1}]_{T,s_{i+1}T}  [t_{i}]_{s_{i+1}T,s_{i}s_{i+1}T} [t_{i+1}]_{s_{i}s_{i+1}T,s_i s_{i+1}s_iT} v_{s_i s_{i+1}s_iT}.
\end{align*}
To check the identity $t_{s_{i}}t_{s_{i+1}}t_{s_{i}} v_T  = t_{s_{i+1}}t_{s_{i}}t_{s_{i+1}} v_T$, we show that each coefficient in $t_{s_{i}}t_{s_{i+1}}t_{s_{i}} v_T -  t_{s_{i+1}}t_{s_{i}}t_{s_{i+1}} v_T$ is 0, noting that if some $S$ does not exist, the result is trivial.

Let $A = c_T(i)$, $B = c_T(i+1)$, and $C =c_T(i+2)$. By definition, for whichever of these $S$ exist, their shifted contents are given by the table in \eqref{tab:length-3-contents}.
%
So, by using the condition that 
	$[t_i]_{T,T} = 1/(c_T(i+1) - c_T(i))$
to simplify the above expansion, we find that the coefficients on each $v_S$, for $S = T$,  $s_iT$, $s_{i+1}T$, $s_i s_{i+1} T$, $s_{i+1}s_i  T$, is 0. The remaining term, 
\begin{align*}
	([t_{i}]_{T,s_iT} 
	[t_{{i+1}}]_{s_iT,s_{i+1}s_iT}&
	[t_{i}]_{s_{i+1}s_iT,s_i s_{i+1}s_iT}\\
	-
	 [t_{i+1}]_{T,s_{i+1}T} &
	 [t_{i}]_{s_{i+1}T,s_{i}s_{i+1}T} 
	 [t_{i+1}]_{s_{i}s_{i+1}T,s_i s_{i+1}s_iT})v_{s_is_{i+1}s_iT} 
\end{align*}
cannot be reduced using the determined values, and so we add the assumption that this coefficient is 0. This concludes the proof of Claim 4. \\
\item[] \textbf{Claim 5:} The braid relation \eqref{eq:braid 4} is satisfied if and only if
\begin{align*}
[x_1]_{s_1 T,s_0s_1 T}&[x_1]_{s_1 s_0 s_1T,s_0s_1 s_0 s_1 T}[t_1]_{T,s_1T}[t_1]_{s_0 s_1T,s_1s_0 s_1T}  \\&= [x_1]_{T,s_0T} [x_1]_{s_1s_0T,s_0s_1s_0T} [t_1]_{s_0T,s_1s_0T}[t_1]_{s_0s_1s_0T,s_1s_0s_1s_0T}.\end{align*}

\emph{Proof:}
Let $a_T = [x_1]_{T,T}$, $b_T = [x_1]_{T,s_0T}$, $d_T = [t_1]_{T,T}$,  $e_T = [t_1]_{T,s_1T}$. So 
\begin{align*}
	x_1 t_{s_1} &v_T 
		= a_{T} d_{T} v_{T} + b_{T} d_{T} v_{s_0 T} + a_{s_1 T} e_{T} v_{s_1T} + b_{s_1 T} e_{T} v_{s_0 s_1 T}, \qquad \text{and}\\
	x_1 t_{s_1}& x_1 t_{s_1} v_T  \\
	&= (a_{T}^2 d_{T}^2 + b_{T}  b_{s_0 T} d_{T}d_{s_0 T} + a_{T} a_{s_1 T} e_{T} e_{s_1T} )  v_{T}\\
	&\quad + 
		(a_{T} b_{T} d_{T}^2 + a_{s_0 T}b_{T} d_{T}  d_{s_0 T}  + a_{s_1 T} b_{T}e_{T}  e_{s_1T}) v_{s_0 T}\\
	&\quad + 
		(a_{T}a_{s_1 T}  d_{T} e_{T} + a_{s_1 T}^2d_{s_1T} e_{T}+ b_{s_1 T}b_{s_0 s_1 T} d_{s_0 s_1 T}e_{T}  ) v_{s_1 T}\\
	&\quad + 
		(a_{T} b_{s_1 T} d_{T} e_{T} + a_{s_1 T} b_{s_1T} d_{s_1T}e_{T} +a_{s_0 s_1 T}b_{s_1 T} e_{T} d_{s_0 s_1 T} ) v_{s_0 s_1T}\\
	&\quad + 
		(a_{s_1 s_0 T}b_{T} d_{T}  e_{s_0 T}) v_{s_1s_0 T}\\
	&\quad + 
		(b_{s_1 s_0 T} b_{T} d_{T} e_{s_0 T} ) v_{s_0s_1s_0 T}\\
	&\quad + 
		(a_{s_1 s_0 s_1 T}b_{s_1 T} e_{T} e_{s_0 s_1 T} ) v_{s_1s_0s_1 T}\\
	&\quad + 
		(b_{s_1 T}b_{s_1 s_0 s_1 T} e_{T} e_{s_0 s_1 T} ) v_{s_0s_1s_0s_1 T},
\end{align*}
and so
\begin{align*}
	(x_1 t_{s_1}& x_1 t_{s_1} +x_1 t_{s_1}) v_T \\
	&= (a_{T}^2 d_{T}^2 + b_{T}  b_{s_0 T} d_{T}d_{s_0 T} + a_{T} a_{s_1 T} e_{T} e_{s_1T} + a_{T} d_{T}  )  v_{T}\\
	&\quad + 
		(a_{T} b_{T} d_{T}^2 + a_{s_0 T}b_{T} d_{T}  d_{s_0 T}  + a_{s_1 T} b_{T}e_{T}  e_{s_1T} + b_{T} d_{T}) v_{s_0 T}\\
	&\quad + 
		(a_{T}a_{s_1 T}  d_{T} e_{T} + a_{s_1 T}^2d_{s_1T} e_{T}+ b_{s_1 T}b_{s_0 s_1 T} d_{s_0 s_1 T}e_{T}  + a_{s_1 T} e_{T} ) v_{s_1 T}\\
	&\quad + 
		(a_{T} b_{s_1 T} d_{T} e_{T} + a_{s_1 T} b_{s_1T} d_{s_1T}e_{T} +a_{s_0 s_1 T}b_{s_1 T} e_{T} d_{s_0 s_1 T} + b_{s_1 T} e_{T} ) v_{s_0 s_1T}\\
	&\quad + 
		(a_{s_1 s_0 T}b_{T} d_{T}  e_{s_0 T}) v_{s_1s_0 T}\\
	&\quad + 
		(b_{s_1 s_0 T} b_{T} d_{T} e_{s_0 T} ) v_{s_0s_1s_0 T}\\
	&\quad + 
		(a_{s_1 s_0 s_1 T}b_{s_1 T} e_{T} e_{s_0 s_1 T} ) v_{s_1s_0s_1 T}\\
	&\quad + 
		(b_{s_1 T}b_{s_1 s_0 s_1 T} e_{T} e_{s_0 s_1 T} ) v_{s_0s_1s_0s_1 T}.
\end{align*}
Similarly, 
since 
 $s_0s_1s_0s_1 T = s_1s_0s_1s_0 T$,
\begin{align*}
(t_{s_1} x_1 &t_{s_1} x_1 +t_{s_1}x_1 ) v_T \\
		&=( a_{T}^2 d_{T}^2 +b_{T}  b_{s_0T} d_{T}d_{s_0T} +a_{ T}  a_{s_1 T} e_{T}e_{s_1T} +a_{T} d_{T})v_T\\
		&\quad 
		+( a_{T}b_{T} d_{T} d_{s_0T} +  a_{s_0T} b_{T} d_{s_0T}^2 +a_{s_1s_0 T} b_{T} e_{s_0T}e_{s_1s_0T} + b_{T} d_{s_0T})v_{s_0T}\\
		&\quad 
		+(a_{T}^2 d_{T}e_{T}+b_{T} b_{s_0T}d_{s_0T}  e_{T} +a_{ T}  a_{s_1T} d_{s_1T}e_{T} + a_{ T} e_{T})v_{s_1T}\\
		&\quad 
		+( a_{ T} b_{s_1T} d_{s_0s_1T} e_{T} )v_{s_0s_1T}\\
		&\quad 
		+(a_{T} b_{T} d_{T}e_{s_0T}+a_{s_0 T}b_{T} d_{s_0T}  e_{s_0T}+a_{s_1s_0T}b_{T} d_{s_1s_0T}e_{s_0T} + b_{T} e_{s_0T})v_{s_1s_0T}\\
		&\quad 
		+( b_{T}b_{s_1s_0T}d_{s_0s_1s_0T} e_{s_0T})v_{s_0s_1s_0T}\\
		&\quad 
		+( a_{ T} b_{s_1T} e_{T} e_{s_0s_1T} )v_{s_1s_0s_1T}\\
		&\quad 
		+(b_{T}b_{s_1s_0T}  e_{s_0T}e_{s_0s_1s_0T} )v_{s_0s_1s_0s_1T}.
\end{align*}

Let $A = c_T(1)$ and $B = c_T(2)$. By definition, for whichever of these $S$ exist, their shifted contents are given by the table in \eqref{tab:length-4-contents}. 
Thus the values of $a_S$ and $d_S$ are given by
\begin{align*}
&\begin{array}{r|cccc}
	S\to		& T 		& s_0T 	& s_1T& s_0s_1 T \\
	\hline
	a_S	&a_T			
		&- a_{T} + (a-p)		
		&a_{s_1T}  		
		& - a_{s_1T} + (a-p)	
		\\
	d_S	&\hbox{$\frac{1}{B - A}$}	
		&\hbox{$\frac{1}{B + A}$}	
		&-d_T					
		&d_{s_0T}				
		\\
\end{array}
\\
&\begin{array}{r|cccc}
	S\to		& s_1s_0  T & s_0s_1s_0  T& s_1s_0s_1  T& s_0s_1s_0s_1  T \\
	\hline
	a_S			
		&a_{s_1T}		
		&- a_{s_1T} + (a-p)	
		&a_T			
		&- a_{T} + (a-p)		
		\\
	d_S	
		&-d_{s_0T}				
		&d_T					
		&-d_{s_0T}				
		&-d_{T}					
		\\
\end{array}
\end{align*}
Furthermore recall that 
 $b_T b_{s_0 T} = -a_T^2 + (a-p)a_T + ap$ and $e_{s_1T}e_T = 1 - d_T^2$. Using these values, we can simplify the expansion of 
 $$\big((x_1 t_{s_1} x_1 t_{s_1} +x_1 t_{s_1}) - (t_{s_1} x_1 t_{s_1}x_1  +t_{s_1}x_1 )\big) v_T$$
 to find that the coefficients of $v_S$ for $S= T$, $s_0T$, $s_1T$, $s_0s_1 T$, $s_1s_0  T$, $s_0s_1s_0  T$,  and $s_1s_0s_1T$ are 0. The remaining term, 
\begin{align*}
\big([x_1]_{s_1 T,s_0s_1 T}&[x_1]_{s_1 s_0 s_1T,s_0s_1 s_0 s_1 T}[t_1]_{T,s_1T}[t_1]_{s_0 s_1T,s_1s_0 s_1T}  \\&- [x_1]_{T,s_0T} [x_1]_{s_1s_0T,s_0s_1s_0T} [t_1]_{s_0T,s_1s_0T}[t_1]_{s_0s_1s_0T,s_1s_0s_1s_0T}\big)v_{s_0s_1s_0s_1  T},\end{align*}
cannot be reduced using the determined values, and so we add the assumption that this coefficient is 0. This concludes the proof of Claim 5.

\end{enumerate}
This concludes Part 1, showing that $\cH^\lambda$ is a $\cH_k^\ext$-module. \\

\paragraph{\textbf{Part 2: $\cH^\lambda$ is simple.}} ~

\noindent We first show that any nontrivial submodule of $\cH^\lambda$ contains some basis element $v_T$. We then prove that any basis element $v_T$ generates  $\cH^\lambda$, and conclude that  $\cH^\lambda$ contains no nontrivial proper submodules.

\begin{itemize}\item[]\textbf{Claim 1:} If  $0 \neq v \in \cH^\lambda$, then $\cH_k^\ext v$ contains some element of the basis $v_T$. \\
\emph{Proof.}
For any $S \in \cT_\lambda$, let
$$W_S = (w_1 - c_S(1))^2 + (w_2 - c_S(2))^2 +  \cdots + (w_k - c_S(k))^2.$$
By Lemma \ref{lem:content_list}, 
	$$W_S v_T = \left(\sum_{i = 1}^k (c_T(i) - c_S(i))^2\right) v_T = 0 \quad \text{ if and only if } \quad  T=S.$$
Therefore, if
	$$\mathrm{Pr}_T =\prod_{\begin{picture}(10,18)(9,0)\put(0,10){\footnotesize$S \in \cT_\lambda$}\put(0,0){\footnotesize $S \neq T$} \end{picture}}  \left(\frac{W_S}{\sum_{i = 1}^k (c_T(i) - c_S(i))^2}\right) \quad \text{ then } \quad \mathrm{Pr}_T v_S = \delta_{ST} v_T.$$
Write $$v = \sum_{S \in \cT_\lambda} d_S v_S, \qquad d_S \in \CC.$$ 
Since $v \neq 0$, there is some $d_T \neq 0$, and so 
	$v_T = \frac{1}{d_T} \mathrm{Pr}_T v \in \cH_k^\ext v,$
concluding the proof of Claim 1. 
\end{itemize}



If $c_T(1) \neq \pm\half((a+p) \pm (b+q))$, then $[x_1]_{T,s_0T} \neq 0$.  
Define the operator $\sigma_0$ on the basis $\{v_T\}_{T \in \cT_\lambda}$ of $\cH^\lambda$ by
	\begin{equation}\sigma_0 v_T = \begin{cases} 
				0 & \text{ if } c_T(1) = \pm\half((a+p) \pm (b+q)),\\
				\frac{1}{[x_1]_{T,s_0T}}\left(x_1 - [x_1]_{T,T}\right)v_T & \text{ otherwise,}
			\end{cases}\end{equation}
and extend linearly. Though $\sigma_0$ is not formally an element of $\cH_k^\ext$, it defines operator on $\cH^\lambda$ via $\cH_k^\ext$, i.e. $\sigma_0 v_T \in \cH_k^\ext v_T$. Therefore if $v_{s_0T}$ exists, then
\begin{align}
\sigma_0 v_T &= \frac{1}{[x_1]_{T,s_0T}}\left(x_1 - [x_1]_{T,T}\right)v_T \nonumber\\
	&=  \frac{1}{[x_1]_{T,s_0T}}\left( [x_1]_{T,T}v_T + [x_1]_{T,s_0T}v_{s_0T} 
			-[x_1]_{T,T}v_T \right)  \nonumber\\
 	&= v_{s_0T}, \nonumber
	\end{align}
and so $v_{s_0T} \in \cH^\ext_k v_T$. 

Similarly, if $c_T(i+1) \neq c_T(i) \pm 1$, then $[t_i]_{T,s_iT} \neq 0$. 
Define the operator $\sigma_i$, $i = 1, \dots, k-1$, on the basis $\{v_T\}_{T \in \cT_\lambda}$ of $\cH^\lambda$ by
	\begin{equation}\sigma_iv_T = \begin{cases} 
				0 & \text{ if } c_T(i+1) = c_T(i) \pm 1,\\
				\frac{1}{[t_i]_{T,s_iT}}\left(t_{s_i} - [t_i]_{T,T}\right)  v_T& \text{ otherwise}
			\end{cases}\end{equation}
and extend linearly. Again, $\sigma_i$ is not formally an element of $\cH_k^\ext$, but rather defines an operator on $\cH^\lambda$ via $\cH_k^\ext$.
So if $v_{s_i T} $ exists, we have
\begin{align}
\sigma_i v_T &= \frac{1}{[t_i]_{T,s_iT}}\left(t_{s_i} - [t_i]_{T,T}\right)v_T \nonumber\\
	&= \frac{1}{[t_i]_{T,s_iT}}\left([t_{i}]_{T,T}v_T + [t_i]_{T,s_iT}v_{s_iT}-[t_i]_{T,T}v_T\right)  \nonumber\\
 	&= v_{s_iT}, \nonumber
	\end{align}
and so $v_{s_iT} \in \cH^\ext_k v_T$.

Recall from Section \ref{sec:seminormal_bases} that we can view every tableau either as a sequence of partitions, as we have been doing, or as a skew shape filled with integers $1, \cdots, k$ with increasing rows and columns. Viewing $T$ as a standard filling now, consider the placement of labels $i$ and $i+1$. If they are adjacent (in row or column), then $c_T(i+1) = c_T(i) \pm 1$, and so $s_iT$ does not exist. However, if labels $i$ and $i+1$ are nonadjacent, then $s_iT$ is gotten from $T$ by switching $i$ and $i+1$. For example, 
	\begin{equation}\label{eq:s2-pic}
	 \phantom{\left| \begin{picture}(1,18)(0,7)\end{picture}\right|} 
	s_2 \cdot
	\def\UNIT{1.7pt}
	\setlength{\unitlength}{\UNIT}	
		\phantom{\left| \begin{picture}(0,15)(0,0) \end{picture}\right.}
		\begin{picture}(35,12.5)(0,0)
		{\color{grey}
		\multiput(0,0)(20,0){2} {\line(0,1){15}} 
		\multiput(0,0)(0,15){2} {\line(1,0){20}} 
		\put(20,15) {\line(1,0){5}} 
		\put(20,10) {\line(1,0){5}} 
		\put(25,15) {\line(0,-1){5}} 
		\put(0,0) {\line(1,0){10}} 
		\put(0,-5) {\line(1,0){10}} 
		\put(0,-10) {\line(1,0){5}} 
		\put(0,0) {\line(0,-1){10}} 
		\put(5,0) {\line(0,-1){10}} 
		\put(10,0) {\line(0,-1){5}} 
		}
		\thicklines
		\put(25,15) {\line(1,0){10}} 
		\put(25,10) {\line(1,0){10}} 
		\put(25,15) {\line(0,-1){5}} 
		\put(30,15) {\line(0,-1){5}} 
		\put(35,15) {\line(0,-1){5}} 
		\put(10,0) {\line(1,0){5}} 
		\put(5,-5) {\line(1,0){10}} 
		\put(5,-10) {\line(1,0){10}} 
		\put(5,-5) {\line(0,-1){5}} 
		\put(10,0) {\line(0,-1){10}} 
		\put(15,0) {\line(0,-1){10}} 
		\put(26,11){\scriptsize 3} 
		\put(31,11){\scriptsize 4} 
		\put(11,-4){\scriptsize 1} 
		\put(6,-9){\scriptsize 2} 
		\put(11,-9){\scriptsize 5} 
		\end{picture}
		=
		\begin{picture}(35,12.5)(0,0)
		{\color{grey}
		\multiput(0,0)(20,0){2} {\line(0,1){15}} 
		\multiput(0,0)(0,15){2} {\line(1,0){20}} 
		\put(20,15) {\line(1,0){5}} 
		\put(20,10) {\line(1,0){5}} 
		\put(25,15) {\line(0,-1){5}} 
		\put(0,0) {\line(1,0){10}} 
		\put(0,-5) {\line(1,0){10}} 
		\put(0,-10) {\line(1,0){5}} 
		\put(0,0) {\line(0,-1){10}} 
		\put(5,0) {\line(0,-1){10}} 
		\put(10,0) {\line(0,-1){5}} 
		}
		\thicklines
		\put(25,15) {\line(1,0){10}} 
		\put(25,10) {\line(1,0){10}} 
		\put(25,15) {\line(0,-1){5}} 
		\put(30,15) {\line(0,-1){5}} 
		\put(35,15) {\line(0,-1){5}} 
		\put(10,0) {\line(1,0){5}} 
		\put(5,-5) {\line(1,0){10}} 
		\put(5,-10) {\line(1,0){10}} 
		\put(5,-5) {\line(0,-1){5}} 
		\put(10,0) {\line(0,-1){10}} 
		\put(15,0) {\line(0,-1){10}} 
		{\color{dblue}
		\put(26,11){\scriptsize 2} 
		\put(31,11){\scriptsize 4} 
		\put(11,-4){\scriptsize 1} 
		\put(6,-9){\scriptsize 3} 
		\put(11,-9){\color{dred}\scriptsize 5} 
		}
		\end{picture}.
		\phantom{\left| \begin{picture}(1,12.5)(0,5)\end{picture}\right|}
\end{equation}
\def\row{\mathrm{row}}
	Define the tableau $\row(T)$ as the filling of $\lambda/T^{(0)}$ built by placing values $1, \dots, k$ left to right, top to bottom, consecutively (this tableau only depends on the shape of the first and last partitions in $T$).

\begin{itemize}\item[]\textbf{Claim 2:} For any tableau $T \in \cT_\lambda$ and any submodule $U \subseteq \cH^\lambda$, 
	$$v_T \in U \quad \text{ if and only if } \quad v_{\row(T)} \in U.$$
\emph{Proof.}
For any $T$, the following process allows us to construct $\row(T)$ by applying a series of $s_i$ moves to $T$. 
	\begin{enumerate}[1.]
	\item Reading left to right, top to bottom, find the first box which has a different filling from $\row(T)$. Let $j$ be the filling in this box and let $i$ be the box immediately before it.
	\item Notice $j-1$ is not placed in any boxes north (east or west) or directly west of $j$, since those boxes are filled with $1, \dots, i$. Therefore, $j-1$ and $j$ can be switched by applying $s_{j-1}$. 
	\item If $s_{j-1}T = \row(T)$, we are done. Otherwise, begin again at step 1 with $s_{j-1}T$. 
	\end{enumerate}
Let $\mathrm{w} = s_{i_\ell} \dots s_{i_2}s_{i_1}$ be the word generated by this process (where $s_{i_1}$ is the first transposition applied, and so on). 
In the example begun in \eqref{eq:s2-pic}, this process unfolds as follows. 
{
\def\UNIT{1.7pt}
\begin{equation}\label{eq:Trow}
	\setlength{\unitlength}{\UNIT}	
		\phantom{\left| \begin{picture}(0,15)(0,0) \end{picture}\right.}
		\begin{picture}(35,12.5)(0,0)
		{\color{grey}
		\multiput(0,0)(20,0){2} {\line(0,1){15}} 
		\multiput(0,0)(0,15){2} {\line(1,0){20}} 
		\put(20,15) {\line(1,0){5}} 
		\put(20,10) {\line(1,0){5}} 
		\put(25,15) {\line(0,-1){5}} 
		\put(0,0) {\line(1,0){10}} 
		\put(0,-5) {\line(1,0){10}} 
		\put(0,-10) {\line(1,0){5}} 
		\put(0,0) {\line(0,-1){10}} 
		\put(5,0) {\line(0,-1){10}} 
		\put(10,0) {\line(0,-1){5}} 
		}
		\put(5,7){$T$}
		\thicklines
		\put(25,15) {\line(1,0){10}} 
		\put(25,10) {\line(1,0){10}} 
		\put(25,15) {\line(0,-1){5}} 
		\put(30,15) {\line(0,-1){5}} 
		\put(35,15) {\line(0,-1){5}} 
		\put(10,0) {\line(1,0){5}} 
		\put(5,-5) {\line(1,0){10}} 
		\put(5,-10) {\line(1,0){10}} 
		\put(5,-5) {\line(0,-1){5}} 
		\put(10,0) {\line(0,-1){10}} 
		\put(15,0) {\line(0,-1){10}} 
		{\color{dblue}
		\put(26,11){\scriptsize 3} 
		\put(31,11){\scriptsize 4} 
		\put(11,-4){\scriptsize 1} 
		\put(6,-9){\scriptsize 2} 
		\put(11,-9){\color{dred}\scriptsize 5} 
		}
		\end{picture}
		\xrightarrow{s_2}
		\begin{picture}(35,12.5)(0,0)
		{\color{grey}
		\multiput(0,0)(20,0){2} {\line(0,1){15}} 
		\multiput(0,0)(0,15){2} {\line(1,0){20}} 
		\put(20,15) {\line(1,0){5}} 
		\put(20,10) {\line(1,0){5}} 
		\put(25,15) {\line(0,-1){5}} 
		\put(0,0) {\line(1,0){10}} 
		\put(0,-5) {\line(1,0){10}} 
		\put(0,-10) {\line(1,0){5}} 
		\put(0,0) {\line(0,-1){10}} 
		\put(5,0) {\line(0,-1){10}} 
		\put(10,0) {\line(0,-1){5}} 
		}
		\thicklines
		\put(25,15) {\line(1,0){10}} 
		\put(25,10) {\line(1,0){10}} 
		\put(25,15) {\line(0,-1){5}} 
		\put(30,15) {\line(0,-1){5}} 
		\put(35,15) {\line(0,-1){5}} 
		\put(10,0) {\line(1,0){5}} 
		\put(5,-5) {\line(1,0){10}} 
		\put(5,-10) {\line(1,0){10}} 
		\put(5,-5) {\line(0,-1){5}} 
		\put(10,0) {\line(0,-1){10}} 
		\put(15,0) {\line(0,-1){10}} 
		{\color{dblue}
		\put(26,11){\scriptsize 2} 
		\put(31,11){\scriptsize 4} 
		\put(11,-4){\scriptsize 1} 
		\put(6,-9){\scriptsize 3} 
		\put(11,-9){\color{dred}\scriptsize 5} 
		}
		\end{picture}
		\xrightarrow{s_1}
		\begin{picture}(35,12.5)(0,0)
		{\color{grey}
		\multiput(0,0)(20,0){2} {\line(0,1){15}} 
		\multiput(0,0)(0,15){2} {\line(1,0){20}} 
		\put(20,15) {\line(1,0){5}} 
		\put(20,10) {\line(1,0){5}} 
		\put(25,15) {\line(0,-1){5}} 
		\put(0,0) {\line(1,0){10}} 
		\put(0,-5) {\line(1,0){10}} 
		\put(0,-10) {\line(1,0){5}} 
		\put(0,0) {\line(0,-1){10}} 
		\put(5,0) {\line(0,-1){10}} 
		\put(10,0) {\line(0,-1){5}} 
		}
		\thicklines
		\put(25,15) {\line(1,0){10}} 
		\put(25,10) {\line(1,0){10}} 
		\put(25,15) {\line(0,-1){5}} 
		\put(30,15) {\line(0,-1){5}} 
		\put(35,15) {\line(0,-1){5}} 
		\put(10,0) {\line(1,0){5}} 
		\put(5,-5) {\line(1,0){10}} 
		\put(5,-10) {\line(1,0){10}} 
		\put(5,-5) {\line(0,-1){5}} 
		\put(10,0) {\line(0,-1){10}} 
		\put(15,0) {\line(0,-1){10}} 
		{\color{dblue}
		\put(26,11){\color{dred}\scriptsize 1} 
		\put(31,11){\scriptsize 4} 
		\put(11,-4){\scriptsize 2} 
		\put(6,-9){\scriptsize 3} 
		\put(11,-9){\color{dred}\scriptsize 5} 
		}
		\end{picture}
		\xrightarrow{s_3}
		\begin{picture}(35,12.5)(0,0)
		{\color{grey}
		\multiput(0,0)(20,0){2} {\line(0,1){15}} 
		\multiput(0,0)(0,15){2} {\line(1,0){20}} 
		\put(20,15) {\line(1,0){5}} 
		\put(20,10) {\line(1,0){5}} 
		\put(25,15) {\line(0,-1){5}} 
		\put(0,0) {\line(1,0){10}} 
		\put(0,-5) {\line(1,0){10}} 
		\put(0,-10) {\line(1,0){5}} 
		\put(0,0) {\line(0,-1){10}} 
		\put(5,0) {\line(0,-1){10}} 
		\put(10,0) {\line(0,-1){5}} 
		}
		\thicklines
		\put(25,15) {\line(1,0){10}} 
		\put(25,10) {\line(1,0){10}} 
		\put(25,15) {\line(0,-1){5}} 
		\put(30,15) {\line(0,-1){5}} 
		\put(35,15) {\line(0,-1){5}} 
		\put(10,0) {\line(1,0){5}} 
		\put(5,-5) {\line(1,0){10}} 
		\put(5,-10) {\line(1,0){10}} 
		\put(5,-5) {\line(0,-1){5}} 
		\put(10,0) {\line(0,-1){10}} 
		\put(15,0) {\line(0,-1){10}} 
		{\color{dblue}
		\put(26,11){\color{dred}\scriptsize 1} 
		\put(31,11){\scriptsize 3} 
		\put(11,-4){\scriptsize 2} 
		\put(6,-9){\color{dred}\scriptsize 4} 
		\put(11,-9){\color{dred}\scriptsize 5} 
		}
		\end{picture}
		\xrightarrow{s_2}
		\begin{picture}(35,12.5)(0,0)
		{\color{grey}
		\multiput(0,0)(20,0){2} {\line(0,1){15}} 
		\multiput(0,0)(0,15){2} {\line(1,0){20}} 
		\put(20,15) {\line(1,0){5}} 
		\put(20,10) {\line(1,0){5}} 
		\put(25,15) {\line(0,-1){5}} 
		\put(0,0) {\line(1,0){10}} 
		\put(0,-5) {\line(1,0){10}} 
		\put(0,-10) {\line(1,0){5}} 
		\put(0,0) {\line(0,-1){10}} 
		\put(5,0) {\line(0,-1){10}} 
		\put(10,0) {\line(0,-1){5}} 
		}
	\put(3,7){$\row(T)$}
		\thicklines
		\put(25,15) {\line(1,0){10}} 
		\put(25,10) {\line(1,0){10}} 
		\put(25,15) {\line(0,-1){5}} 
		\put(30,15) {\line(0,-1){5}} 
		\put(35,15) {\line(0,-1){5}} 
		\put(10,0) {\line(1,0){5}} 
		\put(5,-5) {\line(1,0){10}} 
		\put(5,-10) {\line(1,0){10}} 
		\put(5,-5) {\line(0,-1){5}} 
		\put(10,0) {\line(0,-1){10}} 
		\put(15,0) {\line(0,-1){10}} 
		{\color{dblue}
		\put(26,11){\color{dred}\scriptsize 1} 
		\put(31,11){\color{dred}\scriptsize 2} 
		\put(11,-4){\color{dred}\scriptsize 3} 
		\put(6,-9){\color{dred}\scriptsize 4} 
		\put(11,-9){\color{dred}\scriptsize 5} 
		}
		\end{picture}
	\phantom{\left| \begin{picture}(0,15)(0,0) \end{picture}\right.}
\end{equation}
		}
	So $\mathrm{w} = s_2s_3s_1s_2$ and  $s_2s_3s_1s_2T=\row(T)$.

If $\mathrm{w}T = s_{i_\ell} \dots s_{i_2}s_{i_1}T = \row(T)$, then 
$\sigma_{i_\ell} \dots \sigma_{i_2} \sigma_{i_1}v_T = v_{\row(T)}$
and so $v_{\row(T)} \in \cH^\ext_k v_T$. 
We can apply the same process to find 
 $\mathrm{w}^{-1}\row(T) = T,$
implying
$\sigma_{i_1}\sigma_{i_2} \dots \sigma_{i_\ell}v_{\row(T)} = v_{T}$
and so $v_{T} \in \cH^\ext_k v_{\row(T)}$. This concludes the proof of Claim 2. 
	\end{itemize}
	
Recall from Lemma \ref{lem:z_0_on_lambda} that if $\mu \in \cP$, then $\cB_\mu$ is the set of boxes in rows $p+1$ and below in $\mu$. If $\lambda$ is a partition containing $\mu$, let  $\cB_\mu^\lambda$ be the set of boxes $(i,j)$ in $\cB_\mu$ for which box $(a+b+1-i, p+q+1-j)$ is also in $\lambda$. The criteria in \eqref{eq:rectangle_reqs} imply that the shape obtained by moving each of the boxes $(i,j) \in \cB_\mu^\lambda$ to their complementary position $(a+b+1-i, p+q+1-j)$ gives another partition in $\cP$, which we denote $(\lambda/\mu)^{\max}$. For example, if 
	\begin{equation}\def\UNIT{1pt}
		\begin{matrix}
		\lambda = 
		 \setlength{\unitlength}{\UNIT}	
		\begin{picture}(40,15)(0,2.5)
		\thicklines
		{
			\multiput(0,0)(25,0){2} {\line(0,1){20}} 
		\multiput(0,0)(0,20){2} {\line(1,0){25}} 
		\put(25,20) {\line(1,0){5}} 
		\put(25,15) {\line(1,0){5}} 
		\put(30,20) {\line(0,-1){5}} 
		\put(0,-5) {\line(1,0){15}} 
		\put(0,-10) {\line(1,0){15}} 
		\put(0,0) {\line(0,-1){10}} 
		\put(5,0) {\line(0,-1){10}} 
		\put(10,0) {\line(0,-1){10}} 
		\put(15,0) {\line(0,-1){10}} 
		}
		\put(30,20) {\line(1,0){10}} 
		\put(30,15) {\line(1,0){10}} 
		\put(30,20) {\line(0,-1){5}} 
		\put(35,20) {\line(0,-1){5}} 
		\put(40,20) {\line(0,-1){5}} 
		\end{picture}
		\quad \text{ and } \quad 
		T^{(0)} = 		\begin{picture}(40,15)(0,2.5)
		\thicklines
		{
			\multiput(0,0)(25,0){2} {\line(0,1){20}} 
		\multiput(0,0)(0,20){2} {\line(1,0){25}} 
		\put(25,20) {\line(1,0){5}} 
		\put(25,15) {\line(1,0){5}} 
		\put(30,20) {\line(0,-1){5}} 
		\put(0,-5) {\line(1,0){10}} 
		\put(0,-10) {\line(1,0){5}} 
		\put(0,0) {\line(0,-1){10}} 
		\put(5,0) {\line(0,-1){10}} 
		\put(10,0) {\line(0,-1){5}} 
		}
		\end{picture}
		\\~\\
		\text{then }\quad 
		 \cB_{T^{(0)}}^\lambda = \left\{ (1, p+2) \right\}
		\quad\text{ and } \quad
		(\lambda/\mu)^{\max} = \begin{picture}(40,15)(0,2.5)
		\thicklines
		{
			\multiput(0,0)(25,0){2} {\line(0,1){20}} 
		\multiput(0,0)(0,20){2} {\line(1,0){25}} 
		\put(25,20) {\line(1,0){10}} 
		\put(25,15) {\line(1,0){10}} 
		\put(30,20) {\line(0,-1){5}} 
		\put(35,20) {\line(0,-1){5}} 
		\put(0,-5) {\line(1,0){10}} 
		\put(0,0) {\line(0,-1){5}} 
		\put(5,0) {\line(0,-1){5}} 
		\put(10,0) {\line(0,-1){5}} 
		}
		\end{picture}.
		\end{matrix}\label{ex:lambda_max}\end{equation}

Moreover, since $(\lambda/\mu)^{\max} \in \lambda$, there is a tableau 
	$S = ((\lambda/\mu)^{\max} = S^{(0)}\searrow \dots \searrow S^{(k)} = \lambda)$
from $(\lambda/\mu)^{\max}$ to $\lambda$. Define $T_\lambda = \row(S)$,
i.e.\ $T_\lambda$ is the unique tableau in $\cT_\lambda$ with $T^{(0)}$ highest in lexicographical order and with fillings reading left to right, top to bottom. From the example in \eqref{ex:lambda_max}, $T_\lambda$ is the last tableau pictured below in \eqref{eq:Tlambda}.
\begin{itemize}\item[]\textbf{Claim 3:} For any tableau $T \in \cT_\lambda$ and submodule $U \subseteq \cH^\lambda$, 
	$$v_T \in U \quad \text{ if and only if } \quad  v_{T_\lambda} \in U.$$ 
\emph{Proof.} 
The following process allows us to construct $T_\lambda$ from $T$ through a series of $s_i$:
	\begin{enumerate}[1.]
		\item[0.] Use the process in Claim 2 to move $T$ to $T'=\row(T)$. 
		\item Reading left to right, top to bottom, find the last box $(i,j)$ in $\cB^\lambda_{T'^{(0)}}$.
		\item The box in position $(a+b+1-i, p+q+1-j)$ is filled with a 1. Therefore, we can construct a new tableau 
			$S = (S^{(0)}\searrow \dots\searrow S^{(k)} = \lambda) \in \cT_\lambda$,
			where $S^{(0)}$ is built from $T'^{(0)}$ by moving box $(a+b+1-i, p+q+1-j)$ to $(i,j)$, and $S^{(i)} = T'^{(i)}$ for $i = 1, \dots, k$. The resulting filling will have a 1 in box $(i,j)$ and $2, \dots, k$ identical to $T'$. This new tableau $S$ is equal to $s_0T'$ (see the description of \eqref{def:s_0}). 
		\item Use the process in Claim 2 to move to $\row(s_0T')$.
		\item If $\row(s_0T') = T_\lambda$, then we are done. If not, return to step 1 with $\row(s_0T')$. 
	\end{enumerate}
Let $\mathrm{w} = s_{i_\ell} \dots s_{i_2}s_{i_1}$ be the word generated by this process (where $s_{i_1}$ is the first transposition applied, and so on). Continuing from the example in \eqref{eq:Trow}
 this process proceeds as follows.
{
\def\UNIT{1.7pt}
\begin{equation}\label{eq:Tlambda}
	\setlength{\unitlength}{\UNIT}	
		\phantom{\left| \begin{picture}(0,15)(0,0) \end{picture}\right.}
		\begin{picture}(35,12.5)(0,0)
		{\color{grey}
		\multiput(0,0)(20,0){2} {\line(0,1){15}} 
		\multiput(0,0)(0,15){2} {\line(1,0){20}} 
		\put(20,15) {\line(1,0){5}} 
		\put(20,10) {\line(1,0){5}} 
		\put(25,15) {\line(0,-1){5}} 
		\put(0,0) {\line(1,0){10}} 
		\put(0,-5) {\line(1,0){10}} 
		\put(0,-10) {\line(1,0){5}} 
		\put(0,0) {\line(0,-1){10}} 
		\put(5,0) {\line(0,-1){10}} 
		\put(10,0) {\line(0,-1){5}} 
		}
		\put(5,7){$T$}
		\thicklines
		\put(25,15) {\line(1,0){10}} 
		\put(25,10) {\line(1,0){10}} 
		\put(25,15) {\line(0,-1){5}} 
		\put(30,15) {\line(0,-1){5}} 
		\put(35,15) {\line(0,-1){5}} 
		\put(10,0) {\line(1,0){5}} 
		\put(5,-5) {\line(1,0){10}} 
		\put(5,-10) {\line(1,0){10}} 
		\put(5,-5) {\line(0,-1){5}} 
		\put(10,0) {\line(0,-1){10}} 
		\put(15,0) {\line(0,-1){10}} 
		{\color{dblue}
		\put(26,11){\scriptsize 3} 
		\put(31,11){\scriptsize 4} 
		\put(11,-4){\scriptsize 1} 
		\put(6,-9){\scriptsize 2} 
		\put(11,-9){\color{dred}\scriptsize 5} 
		}
		\end{picture}
		\to \cdots \to 
		\begin{picture}(35,12.5)(0,0)
		{\color{grey}
		\multiput(0,0)(20,0){2} {\line(0,1){15}} 
		\multiput(0,0)(0,15){2} {\line(1,0){20}} 
		\put(20,15) {\line(1,0){5}} 
		\put(20,10) {\line(1,0){5}} 
		\put(25,15) {\line(0,-1){5}} 
		\put(0,0) {\line(1,0){10}} 
		\put(0,-5) {\line(1,0){10}} 
		\put(0,-10) {\line(1,0){5}} 
		\put(0,0) {\line(0,-1){10}} 
		\put(5,0) {\line(0,-1){10}} 
		\put(10,0) {\line(0,-1){5}} 
		}
	\put(3,7){$\row(T)$}
		\thicklines
		\put(25,15) {\line(1,0){10}} 
		\put(25,10) {\line(1,0){10}} 
		\put(25,15) {\line(0,-1){5}} 
		\put(30,15) {\line(0,-1){5}} 
		\put(35,15) {\line(0,-1){5}} 
		\put(10,0) {\line(1,0){5}} 
		\put(5,-5) {\line(1,0){10}} 
		\put(5,-10) {\line(1,0){10}} 
		\put(5,-5) {\line(0,-1){5}} 
		\put(10,0) {\line(0,-1){10}} 
		\put(15,0) {\line(0,-1){10}} 
		{\color{dblue}
		\put(26,11){\color{dred}\scriptsize 1} 
		\put(31,11){\color{dred}\scriptsize 2} 
		\put(11,-4){\color{dred}\scriptsize 3} 
		\put(6,-9){\color{dred}\scriptsize 4} 
		\put(11,-9){\color{dred}\scriptsize 5} 
		}
		\put(1,-9){$*$}
		\end{picture}
		\xrightarrow{s_0}
		\begin{picture}(35,12.5)(0,0)
		{\color{grey}
		\multiput(0,0)(20,0){2} {\line(0,1){15}} 
		\multiput(0,0)(0,15){2} {\line(1,0){20}} 
		\put(20,15) {\line(1,0){10}} 
		\put(20,10) {\line(1,0){10}} 
		\put(25,15) {\line(0,-1){5}} 
		\put(0,0) {\line(0,-1){5}} 
		\put(5,0) {\line(0,-1){5}} 
		}
		\thicklines
		\put(30,15) {\line(1,0){5}} 
		\put(30,10) {\line(1,0){5}} 
		\put(30,15) {\line(0,-1){5}} 
		\put(35,15) {\line(0,-1){5}} 
		\put(10,0) {\line(1,0){5}} 
		\put(0,-5) {\line(1,0){15}} 
		\put(0,-10) {\line(1,0){15}} 
		\put(0,-5) {\line(0,-1){5}} 
		\put(5,-5) {\line(0,-1){5}} 
		\put(10,0) {\line(0,-1){10}} 
		\put(15,0) {\line(0,-1){10}} 
		\put(26,11){$*$} 
		{\color{dblue}
		\put(31,11){\scriptsize 2} 
		\put(11,-4){\scriptsize 3} 
		\put(1,-9){\scriptsize 1} 
		\put(6,-9){\color{dred}\scriptsize 4} 
		\put(11,-9){\color{dred}\scriptsize 5}} 
		\end{picture}
		\xrightarrow{s_1}
		\begin{picture}(35,12.5)(0,0)
		{\color{grey}
		\multiput(0,0)(20,0){2} {\line(0,1){15}} 
		\multiput(0,0)(0,15){2} {\line(1,0){20}} 
		\put(20,15) {\line(1,0){10}} 
		\put(20,10) {\line(1,0){10}} 
		\put(25,15) {\line(0,-1){5}} 
		\put(0,0) {\line(0,-1){5}} 
		\put(5,0) {\line(0,-1){5}} 
		}
		\thicklines
		\put(30,15) {\line(1,0){5}} 
		\put(30,10) {\line(1,0){5}} 
		\put(30,15) {\line(0,-1){5}} 
		\put(35,15) {\line(0,-1){5}} 
		\put(10,0) {\line(1,0){5}} 
		\put(0,-5) {\line(1,0){15}} 
		\put(0,-10) {\line(1,0){15}} 
		\put(0,-5) {\line(0,-1){5}} 
		\put(5,-5) {\line(0,-1){5}} 
		\put(10,0) {\line(0,-1){10}} 
		\put(15,0) {\line(0,-1){10}} 
			%
		{\color{dblue}
		\put(31,11){\color{dred}\scriptsize 1} 
		\put(11,-4){\scriptsize 3} 
		\put(1,-9){\scriptsize 2} 
		\put(6,-9){\color{dred}\scriptsize 4} 
		\put(11,-9){\color{dred}\scriptsize 5}} 
		\end{picture}
		\xrightarrow{s_2}
		\begin{picture}(35,12.5)(0,0)
		{\color{grey}
		\multiput(0,0)(20,0){2} {\line(0,1){15}} 
		\multiput(0,0)(0,15){2} {\line(1,0){20}} 
		\put(20,15) {\line(1,0){10}} 
		\put(20,10) {\line(1,0){10}} 
		\put(25,15) {\line(0,-1){5}} 
		\put(0,0) {\line(0,-1){5}} 
		\put(5,0) {\line(0,-1){5}} 
		}
		\put(2,7){$T_\lambda$}
		\thicklines
		\put(30,15) {\line(1,0){5}} 
		\put(30,10) {\line(1,0){5}} 
		\put(30,15) {\line(0,-1){5}} 
		\put(35,15) {\line(0,-1){5}} 
		\put(10,0) {\line(1,0){5}} 
		\put(0,-5) {\line(1,0){15}} 
		\put(0,-10) {\line(1,0){15}} 
		\put(0,-5) {\line(0,-1){5}} 
		\put(5,-5) {\line(0,-1){5}} 
		\put(10,0) {\line(0,-1){10}} 
		\put(15,0) {\line(0,-1){10}} 
		{\color{dblue}
		\put(31,11){\color{dred}\scriptsize 1} 
		\put(11,-4){\color{dred}\scriptsize 2} 
		\put(1,-9){\color{dred}\scriptsize 3} 
		\put(6,-9){\color{dred}\scriptsize 4} 
		\put(11,-9){\color{dred}\scriptsize 5}} 
		\end{picture}
	\phantom{\left| \begin{picture}(0,15)(0,0) \end{picture}\right.}
\end{equation}
		}
So $\mathrm{w}= s_2s_1s_0s_2s_3s_1s_2$, and $\mathrm{w}T = T_\lambda$.

If $\mathrm{w}T = s_{i_\ell} \dots s_{i_2}s_{i_1}T = T_\lambda$, then 
$\sigma_{i_\ell} \dots \sigma_{i_2} \sigma_{i_1}v_T = v_{T_\lambda}$
and so $v_{T^ \lambda} \in \cH^\ext_k v_T$. 
We can apply the same process to find 
 $\mathrm{w}^{-1}T^ \lambda = s_{i_1} s_{i_2} \dots s_{i_\ell}T_\lambda = T,$
implying
$$\sigma_{i_1}\sigma_{i_2} \dots \sigma_{i_\ell}v_{T_\lambda} = v_{T}$$
and so $v_{T} \in \cH^\ext_k v_{T_\lambda}$. This concludes the proof of Claim 3. 
\end{itemize}


By Claim 1, any nonzero submodule $U \subseteq \cH^\lambda$ contains some basis vector $v_T$. By Claim 3, $U$ therefore contains $v_{T_\lambda}$, and consequently contains all basis vectors $v_T$ of $\cH^\lambda$. Thus, $U = \cH^\lambda$ and so $\cH^\lambda$ is simple. 
This concludes Part 2, and therefore completes the proof of Proposition \ref{thm:seminormal-Hecke}.
\end{proof}

\begin{remark}
\label{rk:seminormal-Hecke}
We have shown slightly more than was stated in Proposition \ref{thm:seminormal-Hecke}. Namely, if $\cH^\lambda$ is a $\cH_k^{\ext}$-module with basis indexed by $T \in \cT_\lambda$ and $w_i  {v}_T = c_T(i) {v}_T$ for $0 = 1, \dots, k$, then 
\begin{enumerate}
	\item $t_{s_i}  {v}_T =  [t_i]_{T,T} {v}_T + [t_i]_{T,s_iT} v_{s_i T}$ and $x_{1}  {v}_T = [x_1]_{T,T} {v}_T +  [x_1]_{T,s_0T} {v}_{s_0T},$
where $[t_i]_{T,s_iT} = 0$ if and only if $c_T(i) = c_T(i+1) \pm 1$, and $ [x_1]_{T,s_0T}=0 $ if and only if  
$c_T(1) = \half(\pm(a+p) \pm (b+q))$,
	\item $[x_1]_{T,S}$ and $[t_i]_{T,S}$ satisfy items (1)-(6) of Proposition \ref{thm:seminormal-Hecke}, and
	\item  $\cH^\lambda$ is simple as an $\cH_k^{\ext}$-module.
	\end{enumerate}
What is more is that the proof that $\cH^\lambda$ is simple (Part 2) relies only on the action of $\cH_k$, and so
	$\Res^{\cH_k^{\ext}}_{\cH_k}\left(\cH^\lambda\right)$ is simple. 
\end{remark}

\begin{cor}\label{cor:punchline}
In the setting of Theorem \ref{thm:Bratteli_action}, 
	$$\Res^{\End_\fg(M \otimes N \otimes V^{\otimes k})}_{\Phi'(\cH_k^\ext)} (\cL^\mu)
	 \quad \text{ and } \quad
	  \Res^{\End_\fg(M \otimes N \otimes V^{\otimes k})}_{\Phi'(\cH_k)} (\cL^\mu)$$
	are simple $\cH_k^\ext$- and $\cH_k$-modules, respectively.
\end{cor}
\begin{proof}
By Theorem \ref{thm:Bratteli_action} any simple $\End_\fg(M \otimes N \otimes V^{\otimes k})$-module $\cL^\mu \subseteq M \otimes N \otimes V^{\otimes k}$ has basis $\{ v_T ~|~ T \in \cT_\mu\}$ on which $w_i$ acts via $\Phi'$ by $w_i  v_T = c_T(i)v_T$. The restatement of Proposition \ref{thm:seminormal-Hecke} in Remark \ref{rk:seminormal-Hecke} implies $\cL^\mu$ is simple as both a $\cH_k^\ext$-module and a $\cH_k$-module. 
\end{proof}



\end{document}